\documentclass[oneside]{amsart}

\pdfoutput=1
\usepackage{color}
\usepackage{changepage}
\usepackage[hyphens]{url}
\usepackage{hyperref}
\usepackage[hyphenbreaks]{breakurl}
\usepackage{caption}
\usepackage{subcaption}
\usepackage{mathdots}
 \usepackage{relsize}

\usepackage[
    backend=biber,
    style=alphabetic,
    sorting=nyt
]{biblatex}
\usepackage{amsmath}
\usepackage{amssymb}
\usepackage{IEEEtrantools}
\usepackage{booktabs}
\usepackage{tikz-cd}
\usetikzlibrary{matrix,arrows,decorations.pathmorphing}

\usepackage{enumitem}
\newlist{step}{enumerate}{1}
\setlist[step]{label=\textbf{Step \arabic*.}}

\usepackage{makecell}
\usepackage{diagbox}
\usepackage{colortbl}

\usepackage{pdflscape}
\usepackage{mathtools}
\usepackage{graphicx}
\usepackage{amsfonts}

\usepackage{nicematrix}
\usepackage{tikz}
\usepackage{stmaryrd}
\usepackage{bm}
\usepackage{amsthm}

\newcommand{\CC}{\mathbb{C}} 
\newcommand{\ZZ}{\mathbb{Z}} 
\newcommand{\PP}{\mathbb{P}} 
\newcommand{\KK}{\mathbb{K}} 

\newcommand{\gl}[1]{\mathrm{GL}(#1)}
\newcommand{\charac}{\mathrm{char}}
\newcommand{\fano}[2]{\mathrm{\textbf{F}}_{#1}(#2)}  
\newcommand{\SD}[2][]{\mathrm{SD}\ifx&#1&_#2^#2\else_#2^#1\fi} 
\newcommand{\AD}[2][]{\mathrm{Pf}\ifx&#1&_#2^#2\else_#2^#1\fi} 
\newcommand{\sym}{\mathcal{S}} 
\newcommand{\alt}{\mathcal{A}} 
\newcommand{\rect}{\mathcal{M}} 
\newcommand{\flag}[2]{\mathrm{Fl}_{#1}(#2)}
\newcommand{\comp}[1]{\mathrm{\textbf{C}}_{\kappa(#1)}(#1)}
\newcommand{\compk}[2]{\mathrm{\textbf{C}}_{#1}(#2)}
\newcommand{\ucompk}[1]{\mathrm{\textbf{C}}_{#1}}
\newcommand{\rowcol}{row-column}
\newcommand{\ann}{\mathrm{Ann}\,}

\newcommand*{\Scale}[2][4]{\scalebox{0.9}{$#2$}}
\newcommand{\outmat}[1]{\scalebox{0.9}{$#1$}}

\DeclareMathOperator{\Hom}{Hom}

\DeclareMathOperator{\Gr}{Gr} 
\DeclareMathOperator{\Sym}{Sym} 
\DeclareMathOperator{\spn}{span} 

\numberwithin{equation}{section}
\newtheorem{theorem}[equation]{Theorem}
\newtheorem{proposition}[equation]{Proposition}
\newtheorem{lemma}[equation]{Lemma}
\newtheorem{conjecture}[equation]{Conjecture}
\newtheorem{corollary}[equation]{Corollary}

\theoremstyle{definition}
\newtheorem{definition}[equation]{Definition}
\newtheorem{remark}[equation]{Remark}
\newtheorem{question}[equation]{Question}
\newtheorem{problem}[equation]{Problem}
\newtheorem{example}[equation]{Example}

\newtheorem{notation}[equation]{Notation}

\allowdisplaybreaks

\title[Fano schemes of symmetric matrices of bounded rank]{Fano schemes of symmetric matrices of bounded rank}
\author{Ahmad Mokhtar}
\address{MACSYS Headquarters, The University of Melbourne, Baldwin Spencer Building, Parkville, Victoria, 3010, Australia}
\email{\href{mailto:ahmad.mokhtar@unimelb.edu.au}{ahmad.mokhtar@unimelb.edu.au}}

\date{}

\begin{document}

    \begin{abstract}
        We study the geometry of the Fano schemes $\fano{k}{\SD[r]{n}}$ of the projective variety $\SD[r]{n}$ defined by the $r\times r$ minors of a symmetric $n\times n$ matrix filled with indeterminates.
        These schemes are fine moduli spaces parameterizing $(k+1)$-dimensional linear spaces of $n\times n$ symmetric matrices of rank less than $r$.
        We prove that the schemes $\fano{k}{\SD[r]{n}}$ can have generically non-reduced components, and characterize their irreducibility, connectedness, and smoothness.
        Our approach to connectedness also applies to Fano schemes of rectangular matrices as well as alternating matrices and answers a question of Ilten and Chan.
        Furthermore, we give a complete description of $\fano{1}{\SD[r]{n}}$ and show that when $r=n$, the Fano schemes of lines have the expected dimension.
        As an application, we provide geometric arguments for several previous results concerning spaces of symmetric matrices of bounded rank.
    \end{abstract}

    \maketitle

    \section{Introduction}\label{sec:introduction}

    \subsection{Motivation}\label{subsec:motivation}
    Let $\KK$ be an algebraically closed field of characteristic not equal to 2.
    For an integer $n\geq 3$, consider the $\KK$-vector space $\sym_n$ of symmetric $n \times n$ matrices with entries in $\KK$.
    Fixing an integer $3\leq r\leq n$, we are interested in studying those subspaces $Q\subseteq \sym_n$ such that every matrix in $Q$ has rank less than $r$.

    Denote by $\SD[r]{n}$ the closed subscheme of $\PP_{\KK}^{\binom{n+1}{2}-1}$ cut out by the $r\times r$ minors of the generic symmetric matrix
    \begin{equation}
        \label{eq:gen-sym-matrix}
        \begin{pmatrix}
            x_{1,1} & x_{1,2} & \cdots & x_{1,n} \\
            x_{1,2} & x_{2,2} &        & x_{2,n} \\
            \vdots  &         & \ddots & \vdots  \\
            x_{1,n} & x_{2,n} & \cdots & x_{n,n}
        \end{pmatrix},
    \end{equation}
    filled with $\binom{n+1}{2}$ independent forms $x_{i,j}$.
    The variety $\SD[r]{n}$ is irreducible and has dimension $n(r-1)-\binom{r-1}{2}-1$ (see~\cite[Theorem 1, Proposition 6.1]{kutz74}).
    In this article we study the \textit{Fano scheme} $\fano{k}{\SD[r]{n}}$ that parameterizes those $k$-planes in $\PP_{\KK}^{\binom{n+1}{2}-1}$ that lie on $\SD[r]{n}$.
    This is a subscheme of the Grassmannian $\Gr\left(k+1,\sym_n\right)$ whose closed points correspond to $(k+1)$-dimensional subspaces $Q\subseteq \sym_n$ containing matrices of rank less than $r$.
    See~\cite[Chapter 6]{3264} for details on Fano schemes.

    Our original motivation for studying the Fano schemes $\fano{k}{\SD[r]{n}}$ was to generalize the work of Ilten and Chan~\cite{fano-nathan-chan} on Fano schemes $\fano{k}{D_{m,n}^r}$ of rectangular matrices of bounded rank to the symmetric setting.
    Their systematic use of the dimensions of tangent spaces to $\fano{k}{D_{m,n}^r}$ at select points characterizes the irreducibility and smoothness of $\fano{k}{D_{m,n}^r}$, and gives results on the components and connectedness of these schemes.
    Their techniques will prove helpful in our study of $\fano{k}{\SD[r]{n}}$, although we will encounter new challenges in the symmetric case.
    The difficulty is due to the existence of generically non-reduced components which prevents us from leveraging the full power of tangent space dimensions.

    The next reason to study $\fano{k}{\SD[r]{n}}$ is to provide geometric proofs for several previous results concerning subspaces of $\sym_n$ of bounded rank.
    One motivation for these results is a classical theorem due to Kronecker and Weierstrass that gives a canonical form for a pencil of singular matrices (see~\cite[Chapter XII]{gantmacher}).
    Even though the classification of subspaces of $\sym_n$ of bounded rank in linear dimension at least 3 is open, there are partial results.
    One is Meshulam's theorem~\cite[Theorem 1]{meshulam-max-dim} giving the maximum dimension of a subspace of symmetric (and alternating) matrices with rank less than $r$.
    In Proposition~\ref{prop:max-dim-non-empty} we give a short proof of this theorem using $\fano{k}{\SD[r]{n}}$.
    Actually, our method treats subspaces of symmetric, alternating and rectangular matrices of bounded rank simultaneously and gives a unified geometric proof of Meshulam's theorem for symmetric and alternating matrices and of Flanders' theorem~\cite[Theorem 1]{flanders} for rectangular matrices.

    We will also reprove (in Corollary~\ref{cor:loewy-radwan}) Loewy and Radwan's theorem~\cite[Theorem 6.1]{loewy-radwan} that classifies subspaces of $\sym_n$ of bounded rank having the maximum dimension in Meshulam's result.
    Loewy and Radwan's proof uses elementary methods in linear algebra and graph theory and has the advantage that it imposes few restrictions on the field $\KK$.
    We present a concise and more insightful geometric proof for an algebraically closed field $\KK$ by studying the Fano schemes $\fano{k}{\SD[r]{n}}$.

    Our third motivation concerns the classification of trilinear forms, i.e.\ elements $\omega\in \KK^{k+1}\otimes \KK^m\otimes \KK^n$ up to the action of $\gl{k+1}\times\gl{m}\times\gl{n}$.
    Classifying the orbits in general is an open problem but particular cases have been studied (see e.g.\ ~\cite{parfenov-orbits} and ~\cite{moduli-333}).
    A closely related problem is the classification of elements in $\KK^{k+1}\otimes\Sym^2(\KK^n)$, or equivalently $\Hom\left((\KK^{k+1})^*, \Sym^2(\KK^n)\right)$, up to $\gl{k+1}\times\gl{n}$.
    With respect to the standard bases, an element $\omega:(\KK^{k+1})^*\to \Sym^2(\KK^n)$ is an $n\times n$ matrix filled with linear forms in the coordinates of $(\KK^{k+1})^*$.
    We may assume, after possibly quotienting by the kernel, that $\omega$ is injective.
    Up to $\gl{k+1}$, $\omega$ then corresponds to the pair $(k+1,Q)$, where $Q$ is the image of $\omega$, a $(k+1)$-dimensional subspace of $n\times n$ symmetric matrices.
    The Fano scheme $\fano{k}{\SD[r]{n}}$ is the moduli space of those $Q$ having elements with rank less than $r$.

    Symmetric matrices are also directly related to quadrics.
    A subspace $Q\subseteq \sym_n$ gives a linear family of quadrics.
    The linear family $Q$ is called singular if every quadric in $Q$ is given by a singular matrix, otherwise $Q$ is non-singular.
    Pencils of quadrics have been extensively studied.
    A theorem due to Weierstrass and Segre classifies the $\gl{n}$-orbits of non-singular pencils of quadrics via a combinatorial object called the Segre symbol (see~\cite[Theorem 1.1]{sturmfels-pencil}).
    In~\cite{mezzetti} the author studies singular pencils of quadrics of constant rank.
    Orbits of nets of conics have been classified in~\cite{nets-of-conics}.
    Our Fano schemes parameterize linear spaces (not necessarily pencils or nets) of quadrics of bounded rank.
    We hope that this study leads to a better understanding of these linear spaces.

    Finally, this study provides a family of examples for general Hilbert schemes.
    For a projective scheme $Y\subseteq \PP_{\KK}^n$ and a polynomial $f(t)\in\KK[t]$, the general Hilbert scheme $\mathrm{Hilb}_f(Y)$ is the fine moduli space parameterizing those closed subschemes of $Y$ with Hilbert polynomial $f$.
    When $Y=\PP_{\KK}^n$, we have the classical Hilbert scheme $\mathrm{Hilb}_f(\PP_{\KK}^n)$ which is known to be connected~\cite{hartshorne-connectedness}.
    But for general Hilbert schemes, connectedness is not guaranteed.
    The Fano schemes $\fano{k}{Y}$ are particular examples of general Hilbert schemes.
    In~\cite{fano-nathan-chan}, the authors show that the Fano schemes $\fano{k}{D_{m,n}^r}$ parameterizing $k$-planes of the space of $m\times n$ matrices of rank less than $r$ is not always connected.
    We will see that $\fano{k}{\SD[r]{n}}$ is also not always connected (see Corollary~\ref{cor:connectedness}) and we will even encounter generically non-reduced components (see Proposition~\ref{prop:when-gen-non-reduced}).
    Our investigation provides a family of examples of general Hilbert schemes that arise naturally but exhibit such geometric pathologies.

    \subsection{Main results}\label{subsec:main-results}
    Let $\rect_{m,n}$, respectively $\alt_n$ be the vector space of all $m\times n$ matrices, respectively $n\times n$ alternating matrices over $\KK$.
    As before, $\sym_n$ denotes the vector space of $n\times n$ symmetric matrices over $\KK$.
    An easy example of a subspace of matrices with rank less than $r$ is given by compression spaces.

    \begin{definition}[Standard compression space]
        \label{def:standard-compression-space}
        Fix $2\leq r\leq m\leq n$ and let $0\leq s\leq r-1$.
        The \textit{standard} rectangular $s$-compression space is the subspace of $\rect_{m,n}$ of all matrices of the block form
        \begin{equation}
            \label{eq:standard-rect-zero-pattern}
            \begin{pNiceArray}{c|cc}[first-row,first-col,margin]
                & \outmat{r-s-1} & & \outmat{s+1+n-r}\\
                & & & \\
                \outmat{s} & * & &* \\
                &&&\\
                \hline
                & &&   \\
                \outmat{m-s}&*&&0\\
                &&   &\\
            \end{pNiceArray}.
        \end{equation}
        Next, fix $3\leq r\leq n$ and let $0\leq s\leq \lfloor\frac{r-1}{2}\rfloor$ (we take $r$ to be even when considering alternating matrices).
        The \textit{standard} symmetric, respectively alternating, $s$-compression space is the subspace of $\sym_n$, respectively of $\alt_n$, of all symmetric, respectively alternating, matrices of the block form
        \begin{equation}
            \label{eq:standard-zero-pattern}
            \begin{pNiceArray}{ccccc}[first-row,first-col]
                &&s&&\outmat{r-2s-1}&\outmat{s+1+n-r}\\
                &&&&& \\
                \outmat{s}&&*&&*&* \\
                & &  &&& \\
                \outmat{r-2s-1}&&*&&*&0\\
                & &&  && \\
                \outmat{s+1+n-r}&&*&&0&0\\
                &&&&&\\
            \end{pNiceArray}.
        \end{equation}

        \phantomsection\label{nota:kappa}
        Define $\kappa(s,m,n,r)$, respectively $\kappa(s,n,r)$, to be the projective dimension of the standard rectangular, respectively symmetric or alternating, $s$-compression space
        \begin{equation} \label{eq:kappa}
            \begin{aligned}
                \kappa(s,m,n,r) &=mn-(s+1+n-r)(m-s)-1 \quad \text{for rectangular matrices},\\
                \kappa(s,n,r) &=
                \begin{cases}
                s(s+1+n-r)+\binom{r-s}{2}-1 & \text{for symmetric matrices},\\
                s(s+1+n-r)+\binom{r-s-1}{2}-1 & \text{for alternating matrices}.
                \end{cases}
            \end{aligned}
        \end{equation}
        For simplicity, we use the notation $\kappa(s)$ when $m, n, r$ are known from the context.
    \end{definition}
    \begin{remark}
    The integer $\kappa(s,m,n,r)+1$ is the number of asterisks in~\eqref{eq:standard-rect-zero-pattern}.
    Similarly, for symmetric, respectively alternating, matrices, $\kappa(s,n,r)+1$ is the number of asterisks on and above, respectively above, the diagonal in~\eqref{eq:standard-zero-pattern}.
    \end{remark}

    The standard compression spaces all have rank less than $r$.
    The general linear group $\gl{m}\times\gl{n}$ acts on $\rect_{m,n}$ via $B_{1} M B_2^t$ for $(B_1,B_2)\in\gl{m}\times\gl{n}$ and $M\in \rect_{m,n}$.
    Similarly $\gl{n}$ acts on $\sym_n$ and $\alt_n$ by $BMB^t$ for $B\in\gl{n}$ and $M\in \sym_n$ or $\alt_n$.
    The matrix $BMB^t$ is obtained by applying to $M$ a sequence of \textit{row-column} operations.
    A row-column operation is an elementary row operation followed by the same operation on the columns of the matrix.
    Spaces of matrices in the orbit of standard compression spaces are known as compression spaces, first introduced in~\cite{atkinson-lloyd}.

    \begin{definition}[Compression space]
        Let $0\leq s\leq r-1$.
        A rectangular $s$-\textit{compression space} is a subspace of a $\left(\gl{m}\times\gl{n}\right)$-translate of the standard rectangular $s$-compression space.
        When $m=n$, a \textit{symmetric} (or \textit{alternating}) $s$-compression space is a rectangular $s$-compression space consisting of symmetric (or alternating) matrices.
    \end{definition}

    Note that the standard symmetric, respectively alternating, $s$-compression space is a symmetric, respectively alternating, $s$-compression space, but a symmetric or alternating $s$-compression space is not necessarily a subspace of a $\gl{n}$-translate of the standard symmetric or alternating $s$-compression space (see Corollary~\ref{cor:all-s-comp-spaces}).

    Our first result is a complete characterization of the connectedness of $\fano{k}{\SD[r]{n}}$.
    Our approach also treats the connectedness of Fano schemes $\fano{k}{D_{m,n}^r}$ of rectangular $m\times n$ matrices of rank less than $r$, and Fano schemes $\fano{k}{\AD[r]{n}}$ of alternating $n\times n$ matrices of rank less than $r$  (see §\ref{subsec:the-three-varieties} for exact definitions of $D_{m,n}^r$ and $\AD[r]{n}$).

    \begin{theorem}
        \label{thm:connectedness}
        Let $X_k$ be one of the schemes $\fano{k}{\SD[r]{n}}$, $\fano{k}{\AD[r]{n}}$ or $\fano{k}{D_{m,n}^r}$ where $r\leq m\leq n$ and $k$ are integers such that $X_k$ is non-empty.
        The compression spaces form a closed subscheme $\ucompk{k}$ of $X_k$ and the natural map $\pi_0(\ucompk{k})\rightarrow\pi_0\left(X_k\right)$ on the sets of connected components of $\ucompk{k}$ and $X_k$  is bijective.
        In particular, $X_k$ is connected if and only if $\ucompk{k}$ is.
    \end{theorem}

    We will prove this in §\ref{sec:connectedness} and then we will describe the connectedness of the compression subscheme $\ucompk{k}$ in terms of the connectedness of a certain graph and obtain the following characterization of the connectedness of $X_k$.

    \begin{corollary}
        \label{cor:connectedness}
        Let $m$, $n$, $r$, $k$, and $X_k$ be as in Theorem~\ref{thm:connectedness}.
        Define
        \[
            \mathbf{V}_k=\left\{ s\in\ZZ \, | \, 0\leq s\leq s_{\max},\,\mathrm{and}\, k\leq\kappa(s) \right\},
        \]
        where
        \[
            s_{\max}=
            \begin{cases}
                \left\lfloor \frac{r-1}{2} \right\rfloor & \text{if } X=\fano{k}{\mathrm{SD}_n^r}\, \text{ or }\, \fano{k}{\mathrm{Pf}_n^r},\\
                r-1 & \text{if } X=\fano{k}{D_{m,n}^r}.
            \end{cases}
        \]
        Assume $s_1,s_2,\ldots,s_{\ell}$ are the distinct elements of $\mathbf{V}_k$ in ascending order.
        Then $X_k$ is disconnected if and only if
        \begin{enumerate}
            \item $\ell\geq 2$, and
            \item At least two terms of the sequence
            \[
                g(\{s_1,s_2\}),g(\{s_2,s_3\}),\ldots,g(\{s_{\ell-1},s_{\ell}\}),g(\{s_1,s_{\ell}\})
            \]
            are less than $k$, where $g(\{s,s'\})=\kappa(s')-(s+1+n-r)(s'-s)$ for $s<s'$.
        \end{enumerate}
    \end{corollary}

    Corollary~\ref{cor:connectedness} completes the partial characterization of the connectedness of $\fano{k}{D_{m,n}^r}$ by Ilten and Chan~\cite[Theorem 5.3]{fano-nathan-chan}.
    In loc.\ sit.\ the authors ask if there exist $m,n,r,k$ such that $\fano{k}{D_{m,n}^r}$ is connected but its compression subscheme is not~\cite[Question 8.2]{fano-nathan-chan}.
    Theorem~\ref{thm:connectedness} answers this in the negative.

    The next result (proved in §\ref{sec:fano-scheme-of-lines}) characterizes the irreducibility of $\fano{k}{\SD[r]{n}}$.
    We note that by Meshulam's theorem~\cite[Theorem 1]{meshulam-max-dim}, the scheme $\fano{k}{\SD[r]{n}}$ is non-empty if and only if $k\leq\max\left\{ \kappa(0),\kappa(s_{\max}) \right\}$ where $s_{\max}=\lfloor\frac{r-1}{2}\rfloor$.
    We will give a short proof of this fact in §\ref{subsec:fixed-points}.

    \begin{theorem}
        \label{thm:irreducibility}
        The Fano scheme $\fano{k}{\SD[r]{n}}$ is irreducible if and only if
        \[
            \max\left\{\kappa(1),\kappa(s_{\max})\right\}< k\leq\kappa(0),
        \]
        or
        \[
            \max\left\{\kappa(0),\kappa(s_{\max}-1)\right\}< k\leq\kappa(s_{\max}).
        \]
    \end{theorem}

    In the particular case of the Fano scheme of lines $\fano{1}{\SD[r]{n}}$ we can describe all the irreducible components.
    We say a subspace of $\sym_n$ is a \textit{nested} symmetric $s$-compression space if it is a subspace of a $\gl{n}$-translate of the standard symmetric $s$-compression space in~\eqref{eq:standard-zero-pattern} (see Definition~\ref{def:nested-comp-space}).

    We will see that the nested symmetric $s$-compression spaces form a closed locus of $\fano{k}{\SD[r]{n}}$.
    We denote by $\compk{k}{s}$ the corresponding closed subscheme with the induced reduced subscheme structure.
    We will prove the following result for symmetric matrices in §\ref{sec:fano-scheme-of-lines} as parts of Theorem~\ref{thm:fano-lines-components-complete} and Corollary~\ref{cor:dim-compk}.
    \begin{theorem}[Fano scheme of lines]
        \label{thm:fano-lines-components}
        Let $n$ and $r$ be integers with $3\leq r\leq n$.
        The scheme $\fano{1}{\SD[r]{n}}$ has $\lfloor\frac{r-1}{2}\rfloor+1$ irreducible components.
        These components taken with the induced reduced subscheme structure are $\compk{1}{s}$ for $s=0, 1,\ldots,\lfloor\frac{r-1}{2}\rfloor$ and $\bigcap_s \compk{1}{s}\neq\emptyset$.
        Moreover, $\dim\compk{1}{s}=nr+(s-1)(n-r)-5$, and in particular when $r=n$, the components are equidimensional and of the expected dimension.
    \end{theorem}

    It is known that if the Fano scheme $\fano{k}{Y}$ of a general hypersurface $Y$ of $\PP^{N}$ of degree $d$ is non-empty, then it has the \textit{expected dimension} $(k+1)(N-k)-\binom{k+d}{k}$.
    When the Fano scheme of a hypersurface, e.g. $\fano{k}{\SD[n]{n}}$, has the expected dimension, it is given by the vanishing of a section of a certain vector bundle on the Grassmannian.
    In particular, its class in the Chow ring of the Grassmannian can be determined using tools from intersection theory (see~\cite[Proposition 6.4]{3264}).

    Theorem~\ref{thm:fano-lines-components} was already known in the special case $r=n$ (see~\cite[§7]{waterhouse}).
    For $r<n$, De Terán, Dmytryshyn, and Dopico proved that over $\CC$ the space of pencils of $n\times n$ symmetric matrices of rank less than $r$ is the union of closures, in Euclidean topology, of $\lfloor\frac{r-1}{2}\rfloor+1$ explicitly defined congruence bundles where a congruence bundle is a union of $\gl{n}$-orbits.
    They showed that in Euclidean topology, these bundle closures are distinct, and they computed the codimension of each closure~\cite[Theorem 3.2, and §5]{dopico}, which gives the dimensions in Theorem~\ref{thm:fano-lines-components}.
    Moreover, they conjectured~\cite[Remark 4.2]{dopico} that the closures, in Zariski topology, of these $\lfloor\frac{r-1}{2}\rfloor+1$ congruence bundles are also distinct and give the irreducible components of the moduli space.
    Although one can verify the equality of the closures in Zariski and Euclidean topologies of each congruence bundle (see our discussion before Theorem~\ref{thm:fano-lines-components-complete}), Theorem~\ref{thm:fano-lines-components} provides an alternative verification of that conjecture for an algebraically closed field $\KK$ with $\charac(\KK)\neq 2$.

    For $k>1$, a full characterization of the irreducible components of $\fano{k}{\SD[r]{n}}$ is difficult since non-compression components appear.
    We suspect that for every $k$, $\compk{k}{s}$ is equal to a component of $\fano{k}{\SD[r]{n}}$ taken with the induced reduced subscheme structure, but we are able to prove this only for certain values of $s$ and $k$, summarized in Table~\ref{tab:summary-symmetric-components}.
    One obstacle in proving $\compk{k}{s}$ gives a component is that the scheme $\fano{k}{\SD[r]{n}}$  may be non-reduced generically on that component as the next result shows.

    \begin{proposition}
        \label{prop:when-gen-non-reduced}
        Let $n,r,k$ and $s$ be integers with $3\leq r\leq n$ and $0\leq s\leq \lfloor\frac{r-1}{2}\rfloor$.
        Assume the closed subscheme $\compk{k}{s}$ of $\fano{k}{\SD[r]{n}}$ is equal to an irreducible component of $\fano{k}{\SD[r]{n}}$ taken with the induced reduced subscheme structure.
        Then $\fano{k}{\SD[r]{n}}$ is reduced generically on that component if and only if $r$ is odd and $s=\frac{r-1}{2}$.
    \end{proposition}

    We will prove Proposition~\ref{prop:when-gen-non-reduced} in §\ref{sec:smoothness}.
    In particular, $\fano{1}{\SD[r]{n}}$ is generically non-reduced on every component (except one component when $r$ is odd).
    Obviously the scheme $\fano{k}{\SD[r]{n}}$ is not smooth when it is generically non-reduced on a component.
    However, when $r$ is odd, the scheme can be smooth.
    In §\ref{subsec:smoothness}, we will prove the following theorem that fully characterizes when $\fano{k}{\SD[r]{n}}$ is smooth.

    \begin{theorem}\label{thm:characterize-smoothness}
    The scheme $\fano{k}{\SD[r]{n}}$ is smooth if and only if $r$ is odd and
    \begin{equation}\label{eq:smoothness-criterion}
    \max\left\{ \kappa(0),\kappa\left(\frac{r-3}{2}\right) \right\} < k \leq \kappa\left(\frac{r-1}{2}\right).
    \end{equation}
    \end{theorem}

    Table~\ref{tab:summary-of-results} summarizes the results of Corollary~\ref{cor:connectedness}, Theorem~\ref{thm:irreducibility}, and Theorem~\ref{thm:characterize-smoothness} for $\fano{k}{\SD[n]{n}}$ for $n\leq 11$.
    We also demonstrate our results by giving a complete description of the Fano schemes of singular plane conics in Example~\ref{ex:3by3matrix}.
    \begin{table}
        \centering
        \scalebox{0.9}{
        \begin{tabular}{c c c c c c c c c c}
            \toprule
            $n$ & 3 & 4 & 5 & 6 & 7 & 8 & 9 & 10 & 11\\
            \midrule
            Non-empty iff $k\leq$ & 2 & 5 & 9 & 14 & 20 & 27 & 35 & 44 & 54\\
            \addlinespace
            Irreducible iff $k=$ & -- & 5 & 9 & 12--14 & 18--20 & 23--27 & 30--35 & 38--44 & 47--54\\
            \addlinespace
            Smooth iff $k=$ & -- & -- & --& --&-- &--&--&--&--\\
            \addlinespace
            Connected iff $k\leq$ & 1 & 3 & 6 & 9 & 14 & 18 & 24 & 30 & 36\\
            or $k=$ & & 5 & 9 & 12--14 & 18--20 & 23--27 & 30--35 & 35,36,38--44 & 45,47--54\\
            \bottomrule
        \end{tabular}}
        \caption{Properties of the Fano scheme $\fano{k}{\SD[n]{n}}$ for $n\leq 11$}
        \label{tab:summary-of-results}
    \end{table}

    \subsection{Related work and organization}\label{subsec:related-work-and-organization}

    The study of linear spaces of matrices with bounded rank has attracted attention in the last sixty years.
    Dieudonné~\cite[Theorem 1]{dieudonne} computed the maximum dimension a subspace of singular square matrices can have and showed subspaces with maximal dimension are compression spaces.
    Flanders generalized this result to rectangular matrices of bounded rank and Atkinson, Lloyd, Beasley and de Seguins Pazzis~\cite{atkinson-lloyd, beasley2, pazzi-revisited} improved Flanders' classification to dimensions near the maximal one.
    For classifications of linear spaces of matrices of low rank, see~\cite{atkinson-rank3,landsberg-rank4}.

    For symmetric matrices, Meshulam~\cite[Theorem 1]{meshulam-max-dim} proved the equivalent of Flanders' result and determined the maximum dimension of a subspace of $\sym_n$ of bounded rank.
    Loewy and Radwan~\cite[Theorem 6.1]{loewy-radwan} showed subspaces having the maximum dimension are compression spaces and de Seguins Pazzis~\cite[Theorem 1.3]{pazzis-classification} generalized this classification to subspaces near the maximal dimension.
    For a study of subspaces of $\sym_n$ of \textit{fixed} rank, see e.g.~\cite{degen-landsberg}.

    Fano schemes have been studied in numerous contexts.
    Fano~\cite{fano1904} studied the family of lines on a cubic hypersurface having a finite number of singularities.
    Altman and Kleiman~\cite{altman-fano} gave a modern treatment of these moduli spaces.
    There have also been extensive studies on Fano schemes of hypersurfaces~\cite{barth-fano,langer-fano} and complete intersections~\cite{debarre-fano, ilten-kelly}.

    In~\cite{fano-nathan-chan}, the authors study the Fano schemes $\fano{k}{D_{m,n}^r}$ and $\fano{k}{P_{m,n}^r}$ where $D_{m,n}^r$ and $P_{m,n}^r$ are subschemes of $\PP_{\KK}^{mn-1}$ given by the vanishing of the $r\times r$ determinants, respectively permanents, of the generic $m\times n$ matrix.
    They give results on irreducibility, smoothness, connectedness and components of these schemes.

    This paper is organized as follows.
    In §\ref{sec:preliminaries} we review preliminary results on compression spaces, ending it with a short proof of Flanders' and Meshulam's theorems.
    In §\ref{sec:connectedness} we characterize the connectedness of the schemes $\fano{k}{\SD[r]{n}}$, $\fano{k}{\AD[r]{n}}$, and $\fano{k}{D_{m,n}^r}$ and prove Theorem~\ref{thm:connectedness} and Corollary~\ref{cor:connectedness}.
    In §\ref{sec:on-the-components} we will prove results on irreducible components of $\fano{k}{\SD[r]{n}}$ and give a geometric proof for Loewy and Radwan's theorem.
    In §\ref{sec:fano-scheme-of-lines} we will prove Theorem~\ref{thm:fano-lines-components} on the Fano scheme of lines and characterize the irreducibility of $\fano{k}{\SD[r]{n}}$ by proving Theorem~\ref{thm:irreducibility}.
    In §\ref{sec:smoothness} we discuss the existence of generically non-reduced components, prove Proposition~\ref{prop:when-gen-non-reduced}, and characterize the smoothness of $\fano{k}{\SD[r]{n}}$ by proving Theorem~\ref{thm:characterize-smoothness}.

    \subsection*{Acknowledgements}
    This paper was written on the unceded traditional territories of the Coast Salish peoples of the Tsleil-Waututh, Squamish, Kwikwetlem, and Musqueam nations.
    I would like to thank Michel Brion and Allen Knutson for helpful discussions and Nathan Ilten for suggesting this problem and being a great mentor to me.
    The author was supported by NSERC and the ARC Centre for the Mathematical Analysis of Cellular Systems (MACSYS).

    \section{Preliminaries}\label{sec:preliminaries}

    \subsection{The varieties $\SD[r]{n}$, $\AD[r]{n}$ and $D_{m,n}^r$}\label{subsec:the-three-varieties}
    In §\ref{subsec:motivation}, we defined the variety $\SD[r]{n}$ of symmetric matrices of bounded rank.
    This paper focuses on the Fano schemes $\fano{k}{\SD[r]{n}}$.
    However, in §\ref{sec:connectedness}, we will treat the connectedness of these schemes together with that of Fano schemes of alternating matrices $\fano{k}{\AD[r]{n}}$ and rectangular matrices $\fano{k}{D_{m,n}^r}$ of bounded rank.
    Here we define $\AD[r]{n}, D_{m,n}^r$ and related terminology.
    For alternating and rectangular matrices, $\KK$ can be of any characteristic (even 2).

    An \textit{alternating} matrix is a skew-symmetric matrix with zeros on the diagonal.
    As before, $\alt_n$ denotes the vector space of $n\times n$ alternating matrices over $\KK$.
    For integers $2<r\leq n$ with $r$ even, we define $\AD[r]{n}$ to be the subscheme of $\PP^{\binom{n}{2}-1}$ defined by the vanishing of the $r\times r$ Pfaffians of the generic alternating matrix in Figure~\ref{fig:gen-alter-matrix}.

    \begin{remark}
        The rank of an alternating matrix over a field is always even.
        For an even integer $r$, the $r\times r$ Pfaffians of the matrix in Figure~\ref{fig:gen-alter-matrix} generate a prime ideal (see~\cite[Theorem 12]{kleppe}) that cuts out the locus of alternating matrices of rank at most $r-2$ (see~\cite[Corollary 2.6]{buchsbaum}).
        We choose Pfaffians over minors because the ideal generated by minors of size $r$ (or $r-1$) is not prime.
        Moreover, in the study of $\fano{k}{\AD[r]{n}}$ we consider $r>2$ since $\AD[2]{n}$ is empty.
        Similarly, for $\fano{k}{\SD[r]{n}}$, the case $r=2$ is trivial as a subspace of symmetric matrices of rank at most 1 is necessarily the span of $u^t u$ for a row vector $u\in\KK^n$.
    \end{remark}

    Next, fix integers $2\leq r\leq m\leq n$.
    We denoted by $\rect_{m,n}$ the vector space of all $m\times n$ matrices over $\KK$.
    Define $D_{m,n}^r$ to be the subscheme of $\PP^{mn-1}$ given by the vanishing of the $r\times r$ minors of the generic $m\times n$ matrix in Figure~\ref{fig:generic-m-by-n}.

    \begin{figure}
        \centering
        \begin{subfigure}{.49\textwidth}
            \centering
            $
                \begin{pmatrix}
                    0        & x_{1,2}  & \cdots & x_{1,n} \\
                    -x_{1,2} & 0        & \cdots & x_{2,n} \\
                    \vdots   &  \vdots        & \ddots & \vdots  \\
                    -x_{1,n} & -x_{2,n} & \cdots & 0
                \end{pmatrix}
            $
            \caption{Alternating}\label{fig:gen-alter-matrix}
        \end{subfigure}
        \begin{subfigure}{.49\textwidth}
            \centering
            $
            \begin{pmatrix}
                x_{1,1} & x_{1,2} & \cdots & x_{1,n} \\
                x_{2,1} & x_{2,2} & \cdots & x_{2,n} \\
                \vdots  &         & \ddots & \vdots  \\
                x_{m,1} & x_{m,2} & \cdots & x_{m,n}
            \end{pmatrix}
            $
            \caption{Rectangular}\label{fig:generic-m-by-n}
        \end{subfigure}
        \caption{Generic alternating and rectangular matrices}
    \end{figure}

    In §\ref{subsec:main-results}, we defined the action of $\gl{n}$ (respectively $\gl{m}\times\gl{n}$) on $\sym_n$ and $\alt_n$ (respectively on $\rect_{m,n}$).
    This induces a corresponding action on the Fano schemes $\fano{k}{\SD[r]{n}}$ and $\fano{k}{\AD[r]{n}}$ (respectively $\fano{k}{D_{m,n}^r}$).
    A linear space of matrices $Q$ can be represented non-uniquely as a matrix with entries that are linear forms by picking a basis in $Q$ (see Example~\ref{ex:0-proper-comp-not-1-proper-comp}).
    We will freely switch between the subspace and the matrix perspective and use the same notation $Q$ to denote a chosen matrix representation.
    For instance, ``the subspace $BQB^t$'' and ``the matrix $BQB^t$'' are both valid phrases, the latter denoting the matrix obtained by  acting on a previously chosen matrix that represents $Q$.

    \subsection{Bilinear forms}\label{subsec:Bilinear-forms}
    An $n\times n$ matrix $A$ over $\KK$ defines a bilinear form $\beta_A:\KK^n\times \KK^n\to \KK$ via $\beta_A(u,w)=u A w^t$ for row vectors $u,w\in\KK^n$.
    A change of basis on $\KK^n$ corresponds to changing the matrix $A$ to $BAB^t$ for some $B\in\gl{n}$.
    We say $u,w\in \KK^n$ are orthogonal with respect to $A$ and write $u\perp_A w$ if $\beta_A(u,w)=0$.
    For a set $Q$ of matrices, $u\perp_Q w$ means $u\perp_A w$ for every $A\in Q$.
    Similarly, for $U,W\subseteq \KK^n$, by $U\perp_Q W$ we mean $u\perp_Q w$ for every $u\in U, w\in W$.
    The next proposition characterizes compression spaces in terms of bilinear forms.
    \begin{proposition}
        \label{prop:characterize-comp-space}
        A linear space of $n\times n$ matrices $Q$ over $\KK$ is an $s$-compression space if and only if there exist subspaces $U,W\subseteq \KK^n$ of dimensions $s+1+n-r$ and $n-s$ respectively, such that $U\perp_Q W$.
    \end{proposition}
    \begin{proof}
        A matrix $A\in Q$ defines a map $L_A:\KK^n\to (\KK^n)^*$ by $v\mapsto (x\mapsto\beta_A(x,v))$ for every $v\in \KK^n$.
        The matrix of $L_A$ with respect to the standard basis on $\KK^n$ and the dual basis on $(\KK^n)^*$ is $A$.
        We also have $u\perp_A w$ if and only if $L_A(u)$ vanishes on $w$.

        Suppose first that $Q$ is an $s$-compression space.
        By definition, there exist subspaces $U\subseteq \KK^n$ and $U'\subseteq (\KK^n)^*$ of dimensions $s+1+n-r$ and $s$, respectively such that $L_A(U)\subseteq U'$ for every $A\in Q$.
        Let $W=\ann(U')=\{w\in \KK^n:f(w)=0,\, \text{for every}\, f\in U'\}$.
        Then $\dim W=n-s$ and $U\perp_Q W$.
        Conversely, suppose $U$ and $W$ have dimensions $s+1+n-r$ and $n-s$ and $U\perp_Q W$.
        Take $U'=\ann(W)$, the subspace of $(\KK^n)^*$ vanishing on $W$.
        Then $L_A(U)\subseteq U'$ for every $A\in Q$.
        Because $\dim U'=s$, we deduce that $Q$ is an $s$-compression space.
    \end{proof}

    An immediate corollary is that when we consider $s$-compression spaces of symmetric or alternating matrices, we may assume $0\leq s\leq\lfloor\frac{r-1}{2}\rfloor$.

    \begin{corollary}
        \label{cor:s-comp}
        Every $s$-compression space $Q\subseteq \sym_n$ or $\alt_n$ is also an $(r-s-1)$-compression space.
    \end{corollary}
    \begin{proof}
        This follows from Proposition~\ref{prop:characterize-comp-space} by changing the roles of $U$ and $W$.
    \end{proof}

    We now describe compression spaces of $\sym_n$ and $\alt_n$ in terms of the standard ones.
    \begin{corollary}
        \label{cor:all-s-comp-spaces}
        Let $0\leq s\leq \lfloor\frac{r-1}{2} \rfloor$.
        Then every symmetric, respectively alternating, $s$-compression space is a subspace of a $\gl{n}$-translate of the standard symmetric, respectively alternating, $s'$-compression space for some $0\leq s'\leq s$.
    \end{corollary}
    \begin{proof}
        Let $Q\subseteq \sym_n$ or $\alt_n$ be an $s$-compression space.
        By Proposition~\ref{prop:characterize-comp-space}, there are subspaces $U,W\subseteq \KK^n$ of dimensions $s+1+n-r$ and $n-s$ respectively, such that $U\perp_Q W$.
        Define $s'=(\dim U\cap W) -(n-r+1)$.
        Note that $0\leq s'\leq s$.
        Since $\dim U\cap W=s'+1+n-r$, we have $\dim (U+W)=n-s'$ and $(U\cap W) \perp_Q (U+W)$.
        Choose a basis for $U\cap W$.
        Extend it to a basis for $U+W$ and further extend it to a basis $\mathcal{B}$ for $\KK^n$.
        With respect to a suitable ordering of the basis $\mathcal{B}$, every matrix in $Q$ has zero blocks similar to the standard $s'$-compression space in~\eqref{eq:standard-zero-pattern}.
    \end{proof}

    We distinguish the case where the orthogonal subspaces in Proposition~\ref{prop:characterize-comp-space} satisfy the additional property that one contains the other.
    \begin{definition}\label{def:nested-comp-space}
        Let $Q\subseteq \sym_n$ or $\alt_n$ be a linear subspace and let $0\leq s\leq \lfloor\frac{r-1}{2}\rfloor$.
        We say $Q$ is a \textit{nested $s$-compression space} if there exist subspaces $U\subseteq W\subseteq\KK^n$ of linear dimensions $s+1+n-r$ and $n-s$ respectively, such that $U\perp_Q W$.
    \end{definition}

    \begin{example}
        \label{ex:0-proper-comp-not-1-proper-comp}
        Nested $s$-compression spaces are exactly subspaces of $\gl{n}$-translates of the standard $s$-compression space.
        On the other hand, by Corollary~\ref{cor:all-s-comp-spaces}, an $s$-compression space is a nested $s'$-compression space for some $0\leq s'\leq s$ but $s'$ need not be the same as $s$.
        For example, let $n=r=3$.
        The subspace of $\sym_n$ given by
        \[
            Q=\begin{pmatrix}
                  z_0 & 0   & 0 \\
                  0   & z_1 & 0 \\
                  0   & 0   & 0
            \end{pmatrix},
        \]
        is both a $0$- and a $1$-compression space but only a nested $0$-compression space.
        Because if it were a nested $1$-compression space, then there would be a subspace $U\subseteq \KK^3$ of dimension 2 that is totally isotropic with respect to every matrix in $Q$.
        If $u=\alpha_1 e_1+\alpha_2 e_2+\alpha_3 e_3$ is an isotropic vector where $e_i$ are the standard basis vectors of $\KK^3$, then $\beta_Q(u,u)=\alpha_1 ^2 z_0 + \alpha_2^2 z_1=0$ for all $z_0,z_1\in \KK$.
        Hence $u\in \mathrm{span}\{e_3\}$, a contradiction to $\dim U=2$.
    \end{example}

    \begin{remark}\label{rem:comp-closed}
        For a fixed $k$, the locus $\compk{k}{s}$ of (nested) $s$-compression spaces of our Fano schemes is closed and irreducible.
        Fix $0\leq s\leq\left( \lfloor \frac{r-1}{2} \rfloor \right)$.
        Let $X_{\kappa(s)}$ be one of the schemes $\fano{\kappa(s)}{\SD[r]{n}}$ or $\fano{\kappa(s)}{\AD[r]{n}}$.
        We have a natural map
        \begin{equation}\label{eq:map-from-flag-variety}
        \flag{s+1+n-r,n-s}{\KK^n}\to X_{\kappa(s)},
        \end{equation}
        sending a flag $U\subseteq W\subseteq \KK^n$ with $\dim U=s+1+n-r$ and $\dim W=n-s$ to the subspace of symmetric or alternating bilinear forms with $U\perp_Q W$.
        The image of this projective map is the locus $\compk{\kappa(s)}{s}$ which is both closed and irreducible.

        Similarly, for  $0\leq s\leq r-1$, the $s$-compression space locus $\compk{\kappa(s)}{s}$ of $\fano{\kappa(s)}{D_{m,n}^r}$ is the image of the map
        \begin{equation}
            \label{eq:map-from-flag-variety-rectangular}
            \Gr(s+1+n-r,n)\times\Gr(n-s,n)\to\fano{\kappa(s)}{D_{m,n}^r},
        \end{equation}
        sending a pair $([U],[W])\in\Gr(s+1+n-r,n)\times\Gr(n-s,n)$ to the linear space of maps $\KK^m\to\KK^n$ that send $U$ into $W$.

        Now let $X_k$ be one of the schemes $\fano{k}{\SD[r]{n}}$, $\fano{k}{\AD[r]{n}}$, or $\fano{k}{D_{m,n}^r}$.
        Assume $k<\kappa(s)$ and let $\mathcal{U}(s)$ be the restriction of the universal bundle on $X_{\kappa(s)}$ to $\comp{s}$.
        The locus $\compk{k}{s}$ of $X_k$ is the image of the natural projective morphism
        \begin{equation}
            \label{eq:map-from-Grassmann-bundle}
            \Gr\left(k+1,\mathcal{U}(s)\right)\to X_k,
        \end{equation}
        that sends a point on the Grassmann bundle to the corresponding point on the Fano scheme.
        This proves the locus $\compk{k}{s}$ is both closed and irreducible.
    \end{remark}

    \begin{definition}\label{def:ck}
    When $X_k=\fano{k}{\SD[r]{n}}$ or $\fano{k}{\AD[r]{n}}$, we define $\compk{k}{s}$ to be the closed irreducible subscheme of $X_k$ consisting of nested $s$-compression spaces with the induced reduced subscheme structure.
    When $X_k=\fano{k}{D_{m,n}^r}$, the notation $\compk{k}{s}$ means the closed irreducible subscheme of $X_k$ consisting of $s$-compression spaces with the induced reduced subscheme structure, which is consistent with the notation of~\cite{fano-nathan-chan}.
    We equip $\ucompk{k}=\bigcup_s \compk{k}{s}$ with the induced reduced subscheme structure and call it the compression subscheme of the Fano scheme.
    \end{definition}

    \subsection{Fixed points}\label{subsec:fixed-points}
    Consider the \textit{Borel subgroup} $G_B(n)\subseteq\gl{n}$ consisting of all invertible upper triangular matrices and the unipotent subgroup $G_U(n)$ consisting of those matrices in $G_B(n)$ with ones on the diagonal.
    We say a point of $\fano{k}{\SD[r]{n}}$ or $\fano{k}{\AD[r]{n}}$ is a \textit{Borel fixed point}, respectively a \textit{unipotent fixed point}, if it is fixed under the action of $G_B(n)$, respectively $G_U(n)$.
    For the action of $\gl{m}\times\gl{n}$ on $\fano{k}{D_{m,n}^r}$, we define the Borel subgroup and the unipotent subgroup to be $G_B(m)\times G_B(n)$ and $G_U(m)\times G_U(n)$, respectively.
    We will frequently use fixed points to study the geometry of our Fano schemes.
    This section establishes several properties of the unipotent fixed points of the Fano schemes $\fano{k}{\SD[r]{n}}, \fano{k}{\AD[r]{n}}$ and $\fano{k}{D_{m,n}^r}$ and fully characterizes the Borel fixed points.

    \begin{notation}
        \label{not:matrix-basis-notation}
        Assume integers $m,n$ are fixed.
        We denote by $e_{(i,j)}$ the $m\times n$ matrix with 1 in the $(i,j)$ entry and zero elsewhere.
        Next assume $m=n$.
        Set $e_{\{i,j\}}=e_{(i,j)}+e_{(j,i)}$ when $i\neq j$ and $e_{\{i,j\}}=e_{(i,j)}$ when $i=j$.
        Also, for $i\neq j$, set $e_{\langle  i,j\rangle}=e_{(i,j)}-e_{(j,i)}$.
    \end{notation}

    \begin{lemma}
        \label{lem:properties-unipotent-fixed-point}
        Let $\mathcal{V}$ be one of the vector spaces $\sym_n, \alt_n$ or $\rect_{m,n}$.
        Suppose $[P]\in\Gr(k+1,\mathcal{V})$ is a unipotent fixed point, where $P=[p_{i,j}]$ is a $(k+1)$-dimensional subspace of $\mathcal{V}$ written as a symmetric, alternating or rectangular matrix whose entries $p_{i,j}$ are linear forms in $k+1$ free variables.
        \begin{enumerate}[label=(\roman*)]
            \item If $p_{i,j}=0$ for some $i,j$ (we also require $i\neq j$ if $\mathcal{V}=\alt_n$), then $p_{i',j'}=0$ for every $i',j'$ with $i'\geq i$ and $j'\geq j$.
            \label{unipotent-zero-entry}
            \item If $\mathcal{V}=\sym_n$ and $e_{\{i,j\}}\in P$ for some $i,j$, then $e_{\{i',j'\}}\in P$ for every $i',j'$ with $i'\leq i$ and $j'\leq j$.
            This statement holds in the alternating setting $\mathcal{V}=\alt_n$ by replacing $e_{\{i,j\}}$ and $e_{\{i',j'\}}$ with $e_{\langle i,j\rangle}$ and $e_{\langle i',j'\rangle}$ and imposing the condition $i'\neq j'$.
            It also holds in the rectangular case $\mathcal{V}=\rect_{m,n}$ by replacing $e_{\{i,j\}}$ and $e_{\{i',j'\}}$ with $e_{(i,j)}$ and $e_{(i',j')}$.
            \label{unipotent-basis-entry}
            \item Assume $p_{i,j}\neq 0$ for some $i,j$.
            Then for all $i',j'$ with $i'< i$, and $j'< j$ (and $i'\neq j'$ when $\mathcal{V}=\alt_n$), we have $e_{\{i',j'\}}\in P$ when $\mathcal{V}=\sym_n$, $e_{\langle  i',j'\rangle}\in P$ when $\mathcal{V}=\alt_n$ and $e_{(i',j')}\in P$ when $\mathcal{V}=\rect_{m,n}$ .
            \label{unipotent-nonzero-entry}
        \end{enumerate}
    \end{lemma}
    \begin{proof}
        See Example~\ref{ex:unipotent-and-borel-fixed} for what a unipotent fixed point looks like.
        We prove the lemma for $\mathcal{V}=\sym_n$.
        The other cases are proved similarly with minor adjustments.\\
        Part~\ref{unipotent-zero-entry}: It is enough to show the claim for $(i',j')=(i,j+1)$ and $(i',j')=(i+1,j)$.
        Assume to the contrary that for one of these two, $p_{i',j'}\neq 0$.
        We consider two cases.
        The first case is when $i\neq j$, i.e.\ $p_{i,j}$ is not on the diagonal.
        We may assume $i< j$.
        If $p_{i,j+1}\neq 0$, take $u=I_n+e_{(j,j+1)}\in G_U(n)$ where $I_n$ is the identity matrix of size $n$.
        The matrix $uPu^t$ is then obtained from $P$ by adding row/column $j+1$ to row/column $j$.
        We see that $[uPu^t]$ is not equal to $[P]$ because the $(i,j)$ entry in $uPu^t$ is nonzero, a contradiction.
        Similarly, if $p_{i+1,j}\neq 0$, we take $u=I_n+e_{(i,i+1)}$ and get $[uPu^t]\neq [P]$.

        The second case is when $i=j$, i.e.\ $P$ is as in Figure~\ref{fig:matrix-P}.
        Take $u=I_n+e_{(i,i+1)}$, then $uPu^t$ will be as in Figure~\ref{fig:matrix-uPut}.
        If $p_{i+1,i+1}+2p_{i,i+1}\neq 0$, then certainly $[uPu^t]\neq[P]$, a contradiction to $[P]$ being a unipotent fixed point.
        On the other hand, if $p_{i+1,i+1}+2p_{i,i+1}= 0$, then we take $u'=I_n+\frac{1}{2}e_{(i,i+1)}$ and then the $(i,i+1)$-entry of $u'P(u')^t$ is zero, while the same entry in $P$ is nonzero.

        \begin{figure}
            \centering
            \begin{subfigure}{.48\textwidth}
                \centering
                \hspace*{-.31\linewidth}
                $
                    \begin{pNiceArray}{ccccc}[first-row,first-col]
                        &&&\outmat{i}&\outmat{i+1}&\\
                        &&&&& \\
                        \outmat{i}&&&0&p_{i,i+1}& \\
                        \outmat{i+1}&&&p_{i,i+1}&p_{i+1,i+1}& \\
                        &&&&&\\
                    \end{pNiceArray}
                $
                \caption{The fixed point $P$}\label{fig:matrix-P}
            \end{subfigure}
            \begin{subfigure}{.48\textwidth}
                \centering
                \hspace*{-.25\linewidth}
                $
                    \begin{pNiceArray}{ccccc}[first-row,first-col]
                        &&&\outmat{i}&\outmat{i+1}&\\
                        &&&&& \\
                        \outmat{i}&&&p_{i+1,i+1}+2p_{i,i+1}&p_{i,i+1}+p_{i+1,i+1}& \\
                        \outmat{i+1}&&&p_{i,i+1}+p_{i+1,i+1}&p_{i+1,i+1}& \\
                        &&&&&\\
                    \end{pNiceArray}
                $
                \caption{$uPu^t$}\label{fig:matrix-uPut}
            \end{subfigure}
            \caption{The matrices in the proof of Lemma~\ref{lem:properties-unipotent-fixed-point} part~\ref{unipotent-zero-entry}}
        \end{figure}

        Part~\ref{unipotent-basis-entry}: It is enough to prove the statement for $(i',j')=(i-1,j)$ and $(i',j')=(i,j-1)$.
        We split the argument into two cases: $i\neq j$ and $i=j$.
        If $i\neq j$, we may assume $i<j$.
        Take $u=I_n+e_{(i-1,i)}$, then $e_{\{i-1,j\}}=ue_{\{i,j\}}u^t-e_{\{i,j\}}\in P$.
        A similar argument with $u=I_n+e_{(j-1,j)}$ proves $e_{\{i,j-1\}}\in P$.
        If $i=j$, take $u=I_n+e_{(i-1,i)}$ and define $v = ue_{\{i,i\}}u^t-e_{\{i,i\}}\in P$.
        Then $e_{\{i-1,i\}}=\frac{1}{2}(3v-uvu^t)\in P$.

        Part~\ref{unipotent-nonzero-entry}:
        Without loss of generality let $i\leq j$.
        By using part~\ref{unipotent-zero-entry} and increasing $i$ and $j$ (if necessary) we may assume $p_{i'',j''}=0$ for every $i'',j''$ such that $(i'',j'')\neq(i,j)$, $i''\geq i$, and $j''\geq j$.
        We show $e_{\{i-1,j-1\}}\in P$.
        The result then follows from part~\ref{unipotent-basis-entry}.
        By assumption $P$ is as in Figure~\ref{fig:fixed-point-P}.
        \begin{figure}
            \centering
            \begin{subfigure}{.49\textwidth}
                \centering
                $
                    \begin{pNiceArray}{ccccc}[first-row,first-col]
                        &&\outmat{j}&&&\\
                        &&&&&\\
                        \outmat{i}&&p_{i,j}&0&\cdots&0 \\
                        &&0&&&\\
                        &&\vdots&&\ddots&\\
                        &&0&&&0\\
                    \end{pNiceArray}
                $
                \caption{The fixed point $P$}\label{fig:fixed-point-P}
            \end{subfigure}
            \begin{subfigure}{.49\textwidth}
                \centering
                $
                \begin{pNiceArray}{ccccccc}[first-row,first-col]
                    &&&&&\outmat{j}&&\\
                    &&&&&0&&\\
                    &&&&&\vdots&& \\
                    &&&&&0&& \\
                    \outmat{i-1}&&&&&p_{i,j}&& \\
                    \outmat{i}&&&&&0&\ldots&0 \\
                    &&&&&\vdots&\ddots&\\
                    &&&&&0&&0\\
                \end{pNiceArray}
                $
                \caption{$P'=uPu^t-P$}\label{fig:matrix-P'}
            \end{subfigure}
            \caption{The matrices in the proof of Lemma~\ref{lem:properties-unipotent-fixed-point} part~\ref{unipotent-nonzero-entry}}
        \end{figure}
        Take $u=I_n+e_{(i-1,i)}$, then the subspace $P'= uPu^t-P$ of $P$ is as in Figure~\ref{fig:matrix-P'}.
        Defining $u'=I_n+e_{(j-1,j)}$, we get the subspace of $P$ with matrix
        $\alpha p_{i,j}e_{\{i-1,j-1\}}=u'P'(u')^t-P'$,
        where $\alpha=1$ if $i<j$ and $\alpha=2$ if $i=j$.
        Now $p_{i,j}\neq0$ implies $e_{\{i-1,j-1\}}\in P$.
    \end{proof}

    \begin{figure}
        \centering
        \begin{subfigure}[b]{.3\textwidth}
            \centering
            $
            \begin{pmatrix}
                *      &         &     & *      & P_1 \\
                *      &         & *   & P_2    & 0      \\
                *      &         & P_3 & 0      & 0      \\
                \vdots & \iddots &     & \vdots & \vdots \\
                P_t    & \cdots  & 0   & 0      & 0
            \end{pmatrix}
            $
            \caption{}\label{fig:general-blocks}
        \end{subfigure}
        \quad\quad\quad
        \begin{subfigure}[b]{.3\textwidth}
            \centering
            $
            \begin{pNiceMatrix}[r,left-margin=0.6em,right-margin=0.6em]
                &                     &   &   & \Block[draw]{2-2}{} & * \\
                &                     &   &   & *                   & * \\
                & \Block[draw]{3-3}{} &   & * & \mathbf{0}          & 0 \\
                &                     & * & * &     0                & 0 \\
                & *                   & * &   &     0                & 0 \\
                \Block[draw]{1-1}{}* & \mathbf{0}  & 0 & 0 & 0      & 0
            \end{pNiceMatrix}
            $
            \caption{}\label{fig:exmaple-blocks}
        \end{subfigure}
        \caption{Partitioning $P'$ in Proposition~\ref{prop:unipotent-fixed-is-compression}: (a) general case (b) the three blocks $P_1,P_2, P_3$ of sizes 2, 3 and 1 in a $6\times 6$ matrix determined by the marked zeros (in bold) on the second anti-diagonal}
    \end{figure}

    \begin{proposition}
        \label{prop:unipotent-fixed-is-compression}
        The following statements hold.
        \begin{enumerate}[label=(\roman*)]
            \item Every unipotent fixed point of $\fano{k}{\SD[r]{n}}$, respectively $\fano{k}{\AD[r]{n}}$, is a subspace of the standard symmetric, respectively alternating, $s$-compression space in~\eqref{eq:standard-zero-pattern} for some $s$.
            \label{unipotent-symmetric-alternating}
            \item Every unipotent fixed point of $\fano{k}{D_{m,n}^r}$ is a subspace of the standard $s$-compression space in~\eqref{eq:standard-rect-zero-pattern} for some $s$.
        \end{enumerate}
        In particular, every unipotent fixed point of $\fano{k}{\SD[r]{n}}$, $\fano{k}{\AD[r]{n}}$ or $\fano{k}{D_{m,n}^r}$ is an $s$-compression space for some $s$.
    \end{proposition}
    \begin{proof}
        Let $[P]$ be a unipotent fixed point of $\fano{k}{\SD[r]{n}}$, $\fano{k}{\AD[r]{n}}$ or $\fano{k}{D_{m,n}^r}$ where $P=[p_{i,j}]$ is a symmetric, alternating or rectangular matrix whose entries are linear forms in $k+1$ free variables.
        Let $P'$ be the upper left $r\times r$ submatrix of $P$.
        We show that at least one entry on the anti-diagonal of $P'$ is zero.
        By Lemma~\ref{lem:properties-unipotent-fixed-point} part~\ref{unipotent-zero-entry}, all entries to the right and bottom of this entry in $P$ will be zero and therefore for some $s$, $P$ will have zero blocks as in~\eqref{eq:standard-zero-pattern} when considering $\fano{k}{\SD[r]{n}}$ or $\fano{k}{\AD[r]{n}}$ and as in~\eqref{eq:standard-rect-zero-pattern} when considering $\fano{k}{D_{m,n}^r}$ (note that in the case of $\fano{k}{\AD[r]{n}}$, $r$ is even and the diagonal and anti-diagonal of $P'$ do not intersect).

        To show an entry on the anti-diagonal of $P'$ is zero, we use the action of the unipotent subgroup to produce from $P$ a matrix $Q$ filled with linear forms such that as vector spaces, $Q$ is a subspace of $P$ and as a matrix with linear forms, the upper left $r\times r$ minor of $Q$ is the product of the anti-diagonal entries of $P'$.
        Since this minor must be zero, then at least one entry on the anti-diagonal of $P'$ is zero.

        In a square matrix of size $n$, we call the entries $(i,j)$ with $i+j=n+2$, the second anti-diagonal.
        We claim we can partition the matrix $P'$ as in Figure~\ref{fig:general-blocks}, where $t\geq 1$ and each $P_i$ is a square matrix such that every entry on its second anti-diagonal is non-zero (in the alternating case, the middle $P_i$ block will have a zero on the intersection of its second anti-diagonal with the diagonal of $P'$, but that will not be a problem).
        Mark the zero entries on the second anti-diagonal of $P'$ (in the alternating case, the zero entry on the intersection of the diagonal and the second anti-diagonal is not marked).
        Assume there are $t-1$ marked entries, say $(1+\ell_i,n+1-\ell_i)$ for $1\leq i\leq t-1$ where $\ell_1,\ldots,\ell_{t-1}$ are in ascending order.
        Defining $\ell_0=0, \ell_t=n$, we set the size of $P_i$ to be $\ell_i-\ell_{i-1}$ for $1\leq i\leq t$.
        See Figure~\ref{fig:exmaple-blocks} for an example.

        By construction every entry on the second anti-diagonal of each $P_i$ is non-zero (except the entry on the intersection of the diagonal and the second anti-diagonal of $P'$ in the alternating case).
        By Lemma~\ref{lem:properties-unipotent-fixed-point} part~\ref{unipotent-zero-entry}, the marked zeros force every entry below and to the right of the $P_i$ blocks to be zero.
        We want to change the entries in $P$ so that in each $P_i$ block the entries above the anti-diagonals are zero and such that the changed matrix is still a subspace of $P$.
        Then the determinant in each  $P_i$ block will be the product of the entries on the anti-diagonal.

        By applying Lemma~\ref{lem:properties-unipotent-fixed-point} part~\ref{unipotent-nonzero-entry} to every entry on the second anti-diagonal of every $P_i$, we can subtract from $P$ suitable multiples of those $e_{\{i,j\}}$, $e_{\langle i,j \rangle}$ or $e_{(i,j)}$  in $P$ (depending on whether we are working with symmetric, alternating, or rectangular matrices respectively) in such a way that in each $P_i$ block, every entry above the anti-diagonal is zero.
        Let $Q$ be the resulting matrix with $Q'$ its upper left $r\times r$ submatrix and denote by $Q_i$ the submatrix of $Q$ corresponding to $P_i$.
        Then as subspaces, $Q$ is contained in $P$ and as a matrix with linear forms we have
        \[
            0=\det Q'=\pm\prod_{i=1}^t \det Q_{i}=\pm\prod_{i=1}^t p_{i,n-i+1}.
        \]
        This shows an entry on the anti-diagonal of $P'$ is zero, completing the proof.
    \end{proof}

    Next we characterize Borel fixed points of the Grassmannian and our Fano schemes.

    \begin{proposition}
        \label{prop:characterize-borel-fixed-points-grassmannian}
        Assume $\mathcal{V}$ is one of the vector spaces $\sym_n$, $\alt_n$ or $\rect_{m,n}$ and let $Q\subseteq \mathcal{V}$ be a linear subspace of dimension $k+1$.
        Then $[Q]\in\Gr\left(k+1,\mathcal{V}\right)$ is a Borel fixed point if and only if $Q$ can be represented as a matrix in linear forms satisfying the following conditions:
        \begin{enumerate}
            \item For every zero entry $(i,j)$ in $Q$ (we also require $i\neq j$ if $\mathcal{V}=\alt_n$), all entries $(i',j')$ with $i'\geq i$ and $j'\geq j$ are also zero,
            \label{borel-condition-propagate-zeros}
            \item Every non-zero entry of $Q$ is a free variable when $\mathcal{V}=\rect_{m,n}$, and every non-zero entry on and above the main diagonal is a free variable when $\mathcal{V}=\sym_n$ or $\alt_n$.
            \label{borel-condition-free-variable}
        \end{enumerate}
    \end{proposition}
    \begin{proof}
        See Example~\ref{ex:unipotent-and-borel-fixed} for a matrix satisfying these conditions.
        Suppose $Q=[q_{i,j}]$ satisfies the two conditions of the statement.
        To prove the claim for $\mathcal{V}=\rect_{m,n}$, let $A=[a_{i,j}]\in G_B(m)$ and $A'=[a'_{i,j}]\in G_B(n)$ be two upper-triangular matrices.
        If $q_{i,j}=0$ for some $(i,j)$, then the $(i,j)$ entry of $AQ(A')^t$ is also zero
        \[
            \sum_{\substack{1\leq \ell\leq m\\ 1\leq \ell'\leq n}} a_{i,\ell} q_{\ell,\ell'} a'_{j,\ell'} =
            \sum_{\substack{i\leq \ell\leq m\\ j\leq \ell'\leq n}} a_{i,\ell} q_{\ell,\ell'} a'_{j,\ell'} = 0,
        \]
        where the first equality is due to $A,A'$ being upper triangular and the second equality is due to condition~\ref{borel-condition-propagate-zeros}.

        Since non-zero entries of $Q$ vary freely, $Q$ has $k+1$ non-zero entries.
        This means the number of non-zero entries of $AQ(A')^t$ is at most $k+1$.
        But $AQ(A')^t$ is also a subspace of $\mathcal{V}$ of dimension $k+1$.
        Therefore, the number of non-zero entries of $AQ(A')^t$ is exactly $k+1$ and they are linearly independent and are positioned exactly the same as non-zero entries of $Q$.
        Thus $Q$ and $AQ(A')^t$ are equal as subspaces of $\mathcal{V}$.
        The proof for $\mathcal{V}=\sym_n$ or $\alt_n$ is similar except that we take $A'=A$ and consider only non-zero entries on and above the diagonals of $Q$ and $AQA^t$.

        Conversely, assume $Q$ is fixed under the action of the Borel subgroup.
        Since every Borel fixed point is a unipotent fixed point, by Lemma~\ref{lem:properties-unipotent-fixed-point} part~\ref{unipotent-zero-entry}, the matrix $Q$ satisfies condition~\ref{borel-condition-propagate-zeros}.
        To show $Q$ satisfies~\ref{borel-condition-free-variable}, it is enough to prove if $q_{i,j}\neq 0$ for some $i,j$, then the subspace $Q$ contains $e_{\{i,j\}}$, $e_{\langle i,j \rangle}$ or $e_{(i,j)}$ depending on whether $\mathcal{V}=\sym_n$, $\alt_n$ or $\rect_{m,n}$ respectively.

        Assume $\mathcal{V}=\sym_n$ or $\alt_n$.
        First we consider the case $i\neq j$.
        For a matrix $A\in G_B(n)$ and an $n\times n$ matrix $M$ in linear forms, denote by $A\cdot M$ the matrix $AMA^t-M$.
        Obviously if $M$ gives a subspace of $Q$, then so does $A\cdot M$.
        Also for $\alpha\in\KK^*$ and $1\leq \ell\leq n$, define $E_{n,\ell,\alpha}$ to be the diagonal $n\times n$ matrix with entry $(\ell,\ell)$ equal to $\alpha$ and ones elsewhere on the diagonal.
        When $\mathcal{V}=\sym_n$, for every $\alpha,\beta\in\KK^*$ we have
        \[
            (\beta-1)(\alpha-1)q_{i,j}e_{\{i,j\}} =
            E_{n,j,\beta}\cdot\left(E_{n,i,\alpha}\cdot Q\right),
        \]
        representing a subspace of $Q$, which gives $e_{\{i,j\}}\in Q$.
        When $\mathcal{V}=\alt_n$, a similar equation with $e_{\{i,j\}}$ replaced by $e_{\langle i,j \rangle}$ holds showing $e_{\langle i,j \rangle}\in Q$.

        If $i=j$, then necessarily $\mathcal{V}=\sym_n$ and the matrix  $E_{n,i,\alpha}\cdot Q$ represents a subspace of $Q$ and has nonzero entries only in the $i$-th row (and column) and in particular the $(i,i)$ entry equals $(\alpha^2-1)q_{i,i}$.
        Using the proof of the claim for entries not on the diagonal, we can zero out all the entries in the $i$-th row (and column) except the entry on the diagonal.
        This shows $e_{\{i,i\}}\in Q$.

        Finally, for $\mathcal{V}=\rect_{m,n}$ we use the notation $(A,B)\cdot M$ to mean $AMB^t-M$ and a similar computation on the matrix $Q$ to get
        \[
            (\beta-1)(\alpha-1)q_{i,j}e_{(i,j)} =
            (I_m,E_{n,j,\beta})\cdot\left( (E_{m,i,\alpha},I_n)\cdot Q \right),
        \]
        representing a subspace of $Q$, hence $e_{(i,j)}\in Q$.
    \end{proof}
    \begin{corollary}
        \label{cor:characterize-borel-fixed-points-fano-scheme}
        Let $X_k$ be one of the schemes $\fano{k}{\SD[r]{n}}$, $\fano{k}{\AD[r]{n}}$ or $\fano{k}{D_{m,n}^r}$.
        A Borel fixed point $[Q]$ of the Grassmannian lies on $X_k$ if and only if for some $s$, $Q$ is a subspace of the standard $s$-compression space in~\eqref{eq:standard-rect-zero-pattern} (in the symmetric and alternating settings) or in~\eqref{eq:standard-zero-pattern} (in the rectangular setting).
        In particular, Borel fixed points of $X_k$ are compression spaces.
    \end{corollary}
    \begin{proof}
        If for some $s$, $Q$ is a subspace of the standard $s$-compression space, then clearly $[Q]\in X_k$.
        Conversely, if $[Q]$ is a Borel fixed point of $X_k$, then it is unipotent fixed and by Proposition~\ref{prop:unipotent-fixed-is-compression}, $Q$ has the required form.
    \end{proof}

    \begin{example}
        \label{ex:unipotent-and-borel-fixed}
        The following are two examples of unipotent fixed points of $\Gr(5,\sym_4)$.
        The point on the right is Borel fixed and also lies on $\fano{4}{\SD[4]{4}}$.
        \[
            \begin{pmatrix}
                z_0 & z_1 & z_2 & z_3 \\
                z_1 & z_4 & z_3 & 0   \\
                z_2 & z_3 & 0   & 0   \\
                z_3 & 0   & 0   & 0
            \end{pmatrix},\qquad
            \begin{pmatrix}
                z_0 & z_1 & z_2 & z_3 \\
                z_1 & z_4 & 0   & 0   \\
                z_2 & 0   & 0   & 0   \\
                z_3 & 0   & 0   & 0
            \end{pmatrix}.
        \]
    \end{example}

    Using Borel fixed points, we can characterize when the schemes $\fano{k}{\SD[r]{n}}$, $\fano{k}{\AD[r]{n}}$ and $\fano{k}{D_{m,n}^r}$ are non-empty.
    Flanders~\cite[Theorem 1]{flanders} proved this for rectangular matrices (see also~\cite[Theorem 1]{dieudonne} for the case $m=n=r$) and Meshulam~\cite[Theorem 1]{meshulam-max-dim} gave a proof for the symmetric and alternating cases.
    We give a unified geometric proof for symmetric, alternating, and rectangular matrices.
    Our argument is similar to the proof of Ilten and Chan~\cite[Proposition 2.6]{fano-nathan-chan} for $\fano{k}{D_{m,n}^r}$ that uses fixed points of a torus action.

    \begin{proposition}
        \label{prop:max-dim-non-empty}
        Let $X_k$ be one of the Fano schemes $\fano{k}{\SD[r]{n}}$, $\fano{k}{\AD[r]{n}}$ or $\fano{k}{D_{m,n}^r}$.
        Then $X_k$ is non-empty if and only if
        \[
            k\leq\max\left\{\kappa(0),\kappa\left( s_{\max} \right)\right\},
        \]
        where $\kappa(s)$ is as in~\eqref{eq:kappa}, $s_{\max}=\left\lfloor\frac{r-1}{2}\right\rfloor$ when $X=\fano{k}{\SD[r]{n}}$ or $\fano{k}{\AD[r]{n}}$ and $s_{\max}=r-1$ when $X=\fano{k}{D_{m,n}^r}$.
    \end{proposition}
    \begin{proof}
        If $k$ satisfies the inequality, then the standard $s$-compression space lies on $X_k$ for $s=0$ or $s=s_{\max}$ and $X_k\neq \emptyset$.
        Conversely, if the projective scheme $X_k$ is non-empty, by Borel's fixed point theorem~\cite[Theorem 10.4]{borel}, it contains a Borel fixed point which by Corollary~\ref{cor:characterize-borel-fixed-points-fano-scheme} is a subspace of the standard $s$-compression space for some $0\leq s\leq s_{\max}$.
        This implies $k\leq\kappa(s)$ and since $\kappa(s)$ is convex in $s$, the integer $k$ will satisfy the inequality of the statement.
    \end{proof}

    \section{Connectedness}\label{sec:connectedness}
    In this section we characterize when the schemes $\fano{k}{\SD[r]{n}}$, $\fano{k}{\AD[r]{n}}$ and $\fano{k}{D_{m,n}^r}$ are connected.
    First we prove Theorem~\ref{thm:connectedness} by using the unipotent fixed points.

    \begin{proof}
    [Proof of Theorem~\ref{thm:connectedness}]
        We saw in Remark~\ref{rem:comp-closed} that $\ucompk{k}$ is a closed subscheme of $X_k$.
        Let $Y_k$ be the unipotent fixed locus of $X_k$.
        By Proposition~\ref{prop:unipotent-fixed-is-compression}, we have the inclusions $Y_k\rightarrow\ucompk{k}\rightarrow X_k$, inducing the maps
        \begin{equation}
            \label{eq:maps-connected-components}
            \pi_0(Y_k)\rightarrow\pi_0(\ucompk{k})\rightarrow\pi_0\left( X_k \right).
        \end{equation}
        The closed subscheme $\ucompk{k}$ is Borel-stable.
        By Borel's fixed point theorem~\cite[Theorem 10.4]{borel}, every connected component of $\ucompk{k}$ contains a Borel fixed point which is unipotent fixed.
        The first map in~\eqref{eq:maps-connected-components} is thus surjective.
        By Horrocks' theorem~\cite[Theorem 6.2]{horrocks}, when an additive group, such as the unipotent subgroups $G_U(n)$ or $G_U(m)\times G_U(n)$, acts on a proper scheme, say $X_k$, then the set of connected components of the unipotent fixed locus is in bijection with the set of connected components of the scheme.
        This implies that the composition in~\eqref{eq:maps-connected-components} is bijective and therefore the second map is bijective.
    \end{proof}

    Our goal is to obtain a combinatorial description of connectedness for the schemes $\fano{k}{\SD[r]{n}}$, $\fano{k}{\AD[r]{n}}$ and $\fano{k}{D_{m,n}^r}$ and prove Corollary~\ref{cor:connectedness}.
    Equivalently, by Theorem~\ref{thm:connectedness}, we may obtain a combinatorial description of connectedness for the compression subscheme $\ucompk{k}=\bigcup_s \compk{k}{s}$.
    We will achieve this by associating each of our Fano schemes with a graph.
    This graph will have the important property that its connected components are in bijection with the connected components of $\ucompk{k}$.

    \begin{definition}
        \label{def:graph}
        Let $X_k$ be one of the schemes $\fano{k}{\SD[r]{n}}$, $\fano{k}{\AD[r]{n}}$ or $\fano{k}{D_{m,n}^r}$ where $r\leq m\leq n$ and $k$ are integers and set
        \[
            s_{\max}=
            \begin{cases}
                \left\lfloor \frac{r-1}{2} \right\rfloor & \text{if } X=\fano{k}{\mathrm{SD}_n^r}\, \text{ or }\, \fano{k}{\mathrm{Pf}_n^r},\\
                r-1 & \text{if } X=\fano{k}{D_{m,n}^r}.
            \end{cases}
        \]
        Define $\mathcal{\tilde{G}}_{X_k}$ to be the complete simple graph $(\mathbf{V},\mathbf{E})$ on the vertex set $\mathbf{V}=\{s\in\ZZ:0\leq s\leq s_{\max}\}$ equipped with a function $g:\mathbf{V}\cup \mathbf{E}\to \ZZ_{>0}$ that labels the vertices and edges of $\mathcal{\tilde{G}}_{X_k}$ as follows.
        For a vertex $s\in \mathbf{V}$ we define $g(s)=\kappa(s)$ and for the edge $\{s,s'\}\in \mathbf{E}$ with $s<s'$ we set
        \begin{equation}
            \label{eq:edge-label}
            g(\{s,s'\})=\kappa(s')-(s+1+n-r)(s'-s).
        \end{equation}
        Next, define the graph $\mathcal{G}_{X_k}$ associated to $X_k$ to be the subgraph of $\mathcal{\tilde{G}}_{X_k}$ obtained by deleting those vertices and edges whose labels are less than $k$.
    \end{definition}
    \begin{remark}
        The graph $\mathcal{G}_{X_k}$ depends on the parameters $m$, $n$, $r$, and $k$ and whether we are in the symmetric, alternating, or rectangular setting while $\mathcal{\tilde{G}}_{X_k}$ depends on all these parameters except $k$.
    \end{remark}
    \begin{example}
        Let $X=\fano{k}{\SD[6]{6}}$.
        Figure~\ref{fig:graph-n=6} shows the graph $\mathcal{G}_{X_k}$ and its labels for three values of $k$.
        The graph in Figure~\ref{fig:complete-graph} is also the complete graph $\mathcal{\tilde{G}}_{X_k}$ for all $k$.
    \end{example}
    \begin{figure}
        \centering
        \begin{subfigure}[b]{.3\textwidth}
            \centering
            \scalebox{0.9}{
            \begin{tikzpicture}[node distance={14mm}, thick, main/.style = {draw, circle}]
                \node[main] (0) {$0$};
                \node[above] at (0.north) {14};
                \node[main] (1) [below right of=0] {$1$};
                \node[below] at (1.south) {11};
                \node[main] (2) [above right of=1] {$2$};
                \node[above] at (2.north) {11};
                \draw (0) -- node[midway, below left, pos=.5] {10} (1);
                \draw (0) -- node[midway, above,  pos=.5] {9} (2);
                \draw (1) -- node[midway, below right, pos=.5] {9} (2);
            \end{tikzpicture}}
            \caption{$k=9$}\label{fig:n=6-k=9}
            \label{fig:complete-graph}
        \end{subfigure}
        \quad
        \begin{subfigure}[b]{.3\textwidth}
            \centering
            \scalebox{0.9}{
            \begin{tikzpicture}[node distance={14mm}, thick, main/.style = {draw, circle}]
                \node[main] (0) {$0$};
                \node[above] at (0.north) {14};
                \node[main] (1) [below right of=0] {$1$};
                \node[below] at (1.south) {11};
                \node[main] (2) [above right of=1] {$2$};
                \node[above] at (2.north) {11};
                \draw (0) -- node[midway, below left, pos=.5] {10} (1);
            \end{tikzpicture}}
            \caption{$k=10$}\label{fig:n=6-k=10}
        \end{subfigure}
        \quad
        \begin{subfigure}[b]{.3\textwidth}
            \centering
            \scalebox{0.9}{
            \begin{tikzpicture}[node distance={14mm}, thick, main/.style = {draw, circle}]
                \node[main] (0) {0};
                \node[above] at (0.north) {14};
                \node (1) [below of=0] {};
            \end{tikzpicture}}
            \caption{$k=12$}\label{fig:n=6-k=12}
        \end{subfigure}
        \caption{The graph $\mathcal{G}_{X_k}$ for $X_k=\fano{k}{\SD[6]{6}}$ and three $k$ values}\label{fig:graph-n=6}
    \end{figure}
    \begin{remark}
        \label{rem:interpret-labels}
        When $X_k=\fano{k}{\SD[r]{n}}$ or $\fano{k}{\AD[r]{n}}$, in Definition~\ref{def:graph}, the label $g(s)$ of the vertex $s\in \mathbf{V}$ denotes the highest integer $k$ such that the scheme $X_k$ contains a Borel fixed point that is a nested $s$-compression space, or equivalently $\compk{k}{s}\neq\emptyset$.
        Similarly, by the next lemma, the label $g\left(\{s,s'\}\right)$ of the edge $\{s,s'\}$ denotes the highest integer $k$ such that $X_k$ contains a Borel fixed point that is both a nested $s$- and a nested $s'$-compression space, or equivalently $\compk{k}{s}\cap\compk{k}{s'}\neq\emptyset$.
        A similar remark holds when $X=\fano{k}{D_{m,n}^r}$ with the expression \textit{nested compression space} substituted by \textit{compression space}.
    \end{remark}

    \begin{lemma}
        \label{lem:borel-compression-is-sub-standard}
        If $X_k=\fano{k}{\SD[r]{n}}$ or $\fano{k}{\AD[r]{n}}$ and $[Q]\in X_k$ is a Borel fixed point that is a nested $s$-compression space, then $Q$ is a subspace of the standard $s$-compression space.
        The statement also holds for $X_k=\fano{k}{D_{m,n}^r}$ if $[Q]$ is a Borel fixed $s$-compression space.
    \end{lemma}
    \begin{proof}
        Let $X_k$ be $\fano{k}{\SD[r]{n}}$ or $\fano{k}{\AD[r]{n}}$.
        By assumption there are subspaces $U\subseteq W\subseteq \KK^n$ of dimensions $s+1+n-r$ and $n-s$ with $U\perp_Q W$.
        Moreover, by Proposition~\ref{prop:characterize-borel-fixed-points-grassmannian}, we can write $Q$ as a matrix in linear forms such that   the $(i,j)$ entry with $i\leq j$ is either zero, or a free variable, say $z_{i,j}$.

        Let $\mathcal{B}_U$ be a $(s+1+n-r)\times n$ matrix with entries in $\KK$ whose rows form a basis for $U$ with respect to the standard basis of $\KK^n$.
        Define a similar $(n-s)\times n$ matrix $\mathcal{B}_W$ for $W$.
        We may assume that $\mathcal{B}_U$ and $\mathcal{B}_W$ are in row-echelon form.
        Let $I_U\subseteq\{1,\ldots,n\}$ be the indices of columns containing leading 1's in $\mathcal{B}_U$ and similarly define $I_W$ for $\mathcal{B}_W$.
        The sets $I_U$ and $I_W$ have sizes $s+1+n-r$ and $n-s$ respectively.
        Assume the elements of $I_U$ and $I_W$ in ascending order are
        \[
                I_U = \{i_1,\ldots,i_{s+1+n-r}\},\qquad I_W = \{i'_1,\ldots,i'_{n-s}\}.
        \]

        We claim the $(i,i')$ entry in $Q$ is zero for every $i\in I_U$ and $i'\in I_W$.
        By construction, there are elements $u=e_i+\sum_{j>i}\alpha_j e_j\in U$ and $w=e_{i'}+\sum_{j>i'}\beta_j e_j\in W$ for suitable scalars $\alpha_j,\beta_j$.
        If the $(i,i')$ entry in $Q$ is not zero, the term $z_{i,i'}$ appears in $uQw^t$, contrary to $u\perp_Q w$.

        Next, we show that for every $i\geq r-s$ and $i'\geq s+1$, the $(i,i')$ entry in $Q$ is zero which puts $Q$ in the form~\eqref{eq:standard-zero-pattern}.
        When $X=\fano{k}{\SD[r]{n}}$, apply Proposition~\ref{prop:characterize-borel-fixed-points-grassmannian} part~\ref{borel-condition-propagate-zeros} to the zero entry $(i_1,i'_1)$ in $Q$.
        Because the dimensions of $U$ and $W$ imply $i_1\leq r-s$ and $i'_1\leq s+1$, we have the required zero blocks in $Q$.

        When $X=\fano{k}{\AD[r]{n}}$, if $i_1\neq i'_1$, then the same argument as in the symmetric case applies.
        Suppose now that $i_1=i'_1\leq s+1$.
        If $i'_2\leq s+1$, we can apply Proposition~\ref{prop:characterize-borel-fixed-points-grassmannian} to the zero entry $(i_1,i'_2)$ because $i_1\neq i'_2$ and we are done.
        If $i'_2>s+1$, then $\dim W=n-s$ gives $i'_2=s+2$.
        On the other hand, $r$ is even in the alternating case and $\dim U\leq \dim W$ yields $i'_2=s+2\leq r-s$.
        Another application of Proposition~\ref{prop:characterize-borel-fixed-points-grassmannian} to the zero entry $(i_1,i'_2)$ in $Q$ gives that for every $i\geq s+1$ and $i'\geq r-s$, the $(i,i')$ entry in $Q$ is zero, putting this matrix in the desired form.
        The proof is the same for $X=\fano{k}{D_{m,n}^r}$ by viewing an $m\times n$ matrix $A\in Q$ as a bilinear map $\KK^m\times\KK^n\to\KK$ via $(u,w)\mapsto uAw^t$ and characterizing compression spaces in terms of fixed orthogonal subspaces.
    \end{proof}

    \begin{lemma}
        \label{lem:connectedness-scheme-graph}
        The set of connected components of a scheme $X_k=\fano{k}{\SD[r]{n}}$, $\fano{k}{\AD[r]{n}}$ or $\fano{k}{D_{m,n}^r}$ is in bijection with the set of connected components of its associated graph $\mathcal{G}_{X_k}$.
        In particular, $X_k$ is connected if and only if $\mathcal{G}_{X_k}$ is.
    \end{lemma}

    \begin{proof}
        Consider the subscheme $\ucompk{k}=\bigcup_s \compk{k}{s}$ of $X$ consisting of all compression spaces.
        By definition $\mathcal{G}_{X_k}$ has a node for each $s$ such that $\compk{k}{s}\neq\emptyset$ and by Remark~\ref{rem:interpret-labels}, two nodes corresponding to $s,s'$ are adjacent if and only if $\compk{k}{s}\cap\compk{k}{s'}\neq \emptyset$.
        When $\compk{k}{s}\cap\compk{k}{s'}\neq \emptyset$, we have that $\compk{k}{s}\cup\compk{k}{s'}$ is connected since each $\compk{k}{s}$ is irreducible (see Remark~\ref{rem:comp-closed}).
        This means the connected components of the graph $\mathcal{G}_{X_k}$ are in bijection with the connected components of $\ucompk{k}$.
        Now use Theorem~\ref{thm:connectedness}.
    \end{proof}

    We now prove Corollary~\ref{cor:connectedness} characterizing when $\mathcal{G}_{X_k}$, and as a result $X_k$, is connected.

    \begin{proof}
    [Proof of Corollary~\ref{cor:connectedness}]
        By Lemma~\ref{lem:connectedness-scheme-graph}, the scheme $X_k$ is disconnected if and only if the graph $\mathcal{G}_{X_k}$ is.
        The set $\mathbf{V}_k$ is the vertex set of $\mathcal{G}_{X_k}$ and the two conditions of the statement are equivalent to the fact that the graph $\mathcal{G}_{X_k}$ has at least two vertices and the cycle $s_1-s_2-\cdots-s_{\ell}-s_1$ is missing at least two edges.
        We need to show $\mathcal{G}_{X_k}$ is disconnected if and only if these two conditions on $\mathcal{G}_{X_k}$ hold.

        One direction is clear.
        Assume that $\mathcal{G}_{X_k}$ has at least two vertices and the cycle is missing at least two edges.
        If $\ell=2$, $\mathcal{G}_{X_k}$ is disconnected, so assume $\ell\geq 3$.
        Consider the graph $\tilde{\mathcal{G}}_{X_k}$ in Definition~\ref{def:graph}.
        Build the subgraph $\mathcal{G}'_{X_k}$ of $\tilde{\mathcal{G}}_{X_k}$ by removing those edges whose labels are less than $k$.
        This is similar to building $\mathcal{G}_{X_k}$ from $\tilde{\mathcal{G}}_{X_k}$ except that we are not deleting any vertices.
        This makes $\mathcal{G}_{X_k}$ a subgraph of $\mathcal{G}'_{X_k}$.

        We claim the following statements hold for $\mathcal{G}'_{X_k}$ (we use $\sim$ for the relation of two vertices being adjacent):
        \begin{enumerate}[wide, labelwidth=!, labelindent=0pt]
            \item If $i\not\sim i+1$ then $s\not\sim i+1$ for $0\leq s\leq i$,
            \label{statement-not-adjacent-before}
            \item If $i\not\sim i+1$ then $i\not\sim s$ for $i+1\leq s\leq s_{\max}$.
            \label{statement-not-adjacent-after}
        \end{enumerate}

        \begin{figure}[t]
            \centering
            \begin{subfigure}[b]{0.49\textwidth}
                \centering
                \scalebox{0.9}{
                \begin{tikzpicture}
                    \coordinate (a) at (0, 0);
                    \coordinate (b) at (.75,.45);
                    \coordinate (c) at (.7,2.1);
                    \coordinate (d) at (-.7,2.1);
                    \coordinate (e) at (-.75,.45);
                    \coordinate (f) at (-1.3,1.1);
                    \coordinate (g) at (1.3,1.1);
                    \draw[fill=black] %
                    (a) circle (2pt) node[below] {$i+1$}
                    (b) circle (2pt) node[below right] {$i+2$}
                    (c) circle (2pt) node[above right] {$s_{\max}$}
                    (d) circle (2pt) node[above left] {$1$}
                    (e) circle (2pt) node[below left] {$i$}
                    (f) circle (2pt) node[left] {$i-1$}
                    (g) circle (2pt) node[right] {$i+3$};

                    \draw[dashed] (g) .. controls +(up:3mm) and +(right:3mm) .. (c);
                    \draw[dashed] (f) .. controls +(up:3mm) and +(left:3mm) .. (d);
                    \draw (c) -- (d);
                    \draw (e) --(f);
                    \draw (b) --(g);
                    \draw (d) -- (e);
                \end{tikzpicture}}
                \caption{The isolated vertex $i+1$}
                \label{fig:single-vertex-isolated}
            \end{subfigure}%
            \begin{subfigure}[b]{0.49\textwidth}
                \centering
                \scalebox{0.9}{
                \begin{tikzpicture}
                    \coordinate (a) at (0, 0);
                    \coordinate (b) at (.75,.45);
                    \coordinate (c) at (.7,2.4);
                    \coordinate (d) at (-.7,2.4);
                    \coordinate (e) at (-.75,.45);
                    \coordinate (f) at (-1.3,1.1);
                    \coordinate (g) at (1.3,1.1);
                    \node at (a) {$\cdots$};
                    \draw[fill=black] %
                    (b) circle (2pt) node[below right] {$j$}
                    (c) circle (2pt) node[above right] {$s_{\max}$}
                    (d) circle (2pt) node[above left] {$1$}
                    (e) circle (2pt) node[below left] {$i+1$}
                    (f) circle (2pt) node[left] {$i$}
                    (g) circle (2pt) node[right] {$j+1$};

                    \draw[dashed] (g) .. controls +(up:3mm) and +(right:3mm) .. (c);
                    \draw[dashed] (f) .. controls +(up:3mm) and +(left:3mm) .. (d);
                    \draw (c) -- (d);
                    \draw (d) -- (g);
                \end{tikzpicture}}
                \caption{Isolated vertices $i+1,\ldots,j$}
                \label{fig:range-vertices-isolated}
            \end{subfigure}\\[1.5ex]
            \begin{subfigure}{\linewidth}
                \centering
                \scalebox{0.9}{
                \begin{tikzpicture}
                    \coordinate (a) at (0,0);
                    \coordinate (b) at (-1,-.8);
                    \coordinate (bc) at (-1.5,-1.5);
                    \coordinate (c) at (-1,-2.2);
                    \coordinate (d) at (0,-3);
                    \coordinate (e) at (2,0);
                    \coordinate (f) at (3,-.8);
                    \coordinate (fg) at (3.5,-1.5);
                    \coordinate (g) at (3,-2.2);
                    \coordinate (h) at (2,-3);
                    \draw[fill=black] %
                    (a) circle (2pt) node[above] {$0$}
                    (b) circle (2pt) node[above left] {$1$}
                    (c) circle (2pt) node[below left] {$i-1$}
                    (d) circle (2pt) node[below] {$i$}
                    (e) circle (2pt) node[above] {$s_{\max}$}
                    (f) circle (2pt) node[above right] {$s_{\max}-1$}
                    (g) circle (2pt) node[below right] {$i+2$}
                    (h) circle (2pt) node[below] {$i+1$};
                    \draw (a) -- (b);
                    \draw[dashed] (b) .. controls (bc) and (bc) .. (c);
                    \draw (c) -- (d);
                    \draw (e) -- (f);
                    \draw[dashed] (f) .. controls (fg) and (fg) .. (g);
                    \draw (g) -- (h);
                    \draw (a) -- (c);
                    \draw (h) -- (f) -- (g);

                    \draw[dashed, rounded corners] (-2, .4) rectangle (.9, -3.5) {};
                    \draw[dashed, rounded corners] (4.6, .4) rectangle (1.1, -3.5) {};
                    \node at (-.5,.7) {$\mathbf{V}'$};
                    \node at (2.5,.7) {$\mathbf{V}''$};
                \end{tikzpicture}}
                \caption{If $i\not\sim i+1$ and $s_{\max}\not\sim 0$ then there are no edges between $\mathbf{V}'$ and $\mathbf{V}''$.}
                \label{fig:disconnected-sets}
            \end{subfigure}
            \caption{The graph $\mathcal{G}'_{X_k}$ in the proof of Corollary~\ref{cor:connectedness}}
            \label{fig:isolated-vertices}
        \end{figure}

        The following are immediate consequences of the above two statements:
        \begin{enumerate}[label=(\alph*),wide, labelwidth=!, labelindent=2pt]
            \item If $i\not\sim i+1\not\sim i+2$, then $i+1$ is an isolated vertex, i.e.\ not adjacent to any other vertex (see Figure~\ref{fig:single-vertex-isolated}),
            \label{claim-one-isolated-vertex}
            \item If $i\not\sim i+1$ and $j\not\sim j+1$ for $0\leq i<j<s_{\max}$, then every vertex $s$ with $i+1\leq s \leq j$ is an isolated
            vertex (see Figure~\ref{fig:range-vertices-isolated}): The function $g(\{s,s+1\})$ is convex in $s$ and by assumption $g(\{i,i+1\}),g(\{j,j+1\})<k$.
            This shows $g(\{s,s+1\})<k$ or equivalently $s\not\sim s+1$ for $i\leq s \leq j$.
            Now use statement~\ref{claim-one-isolated-vertex}.
            \label{claim-multiple-isolated-vertices}
            \item If $i\not\sim i+1$ for some $0\leq i<s_{\max}$ and $s_{\max}\not\sim 0$, then there are no edges between $\mathbf{V}'=\{0,\ldots,i\}$ and $\mathbf{V}''=\{i+1,\ldots,s_{\max}\}$ (see Figure~\ref{fig:disconnected-sets}): By statements~\ref{statement-not-adjacent-before} and~\ref{statement-not-adjacent-after}, $0\not\sim i+1$ and $i\not\sim s_{\max}$, i.e.\ $k$ is larger than all of
            \begin{equation}\label{eq:four-points}
                g(\{0,i+1\}), g(\{i,i+1\}), g(\{i,s_{\max}\}), g(\{0,s_{\max}\}).
            \end{equation}
            We claim $g(\{s,s'\})<k$ for all $0\leq s\leq i$ and $i+1\leq s'\leq s_{\max}$, which is equivalent to no edges between $\mathbf{V}'$ and $\mathbf{V}''$.
            Momentarily think of the quadratic polynomial defining $g(\{s,s'\})$ as a function that accepts a pair of real numbers $(s,s')$ with domain $s<s'$.
            This function is  convex since its Hessian is positive semidefinite.
            In particular when four points in the domain are the vertices of a convex quadrilateral, such as those in~\eqref{eq:four-points}, then the value of the function at every convex combination of the these four points is at most the maximum of the values at the vertices of the quadrilateral.
            This means $g(\{s,s'\})<k$ for $s,s'\in\ZZ_{\geq0}$, $0\leq s\leq i$, and $i+1\leq s'\leq s_{\max}$.
            \label{claim-disconnected-sets}
        \end{enumerate}

        Assuming these claims hold, we prove that $\mathcal{G}_{X_k}$ is disconnected.
        By assumption, at least two edges are missing from $\{s_1,s_2\},\{s_2,s_3\},\ldots,\{s_{\ell-1},s_{\ell}\},\{s_{\ell},s_1\}$.
        There are two cases: If $\{s_i,s_{i+1}\}$ and $\{s_j,s_{j+1}\}$ are missing in $\mathcal{G}_{X_k}$ for some $1\leq i<j<\ell$, then $\{s_i,s_i +1\}$ and $\{s_j,s_j +1\}$ are also missing in $\mathcal{G}'_{X_k}$.
        This is because either $s_{i+1}=s_i+1$ or $s_{i+1}>s_i+1$.
        If the latter holds, then the vertex $s_i+1$ is not present in $\mathcal{G}_{X_k}$, meaning that $k>g(\{s_i+1,s\})$ for all $s\neq s_i+1$, i.e.\ $s_i+1$ is isolated in $\mathcal{G}'_{X_k}$.
        A similar argument shows $s_{j}+1$ is isolated in $\mathcal{G}'_{X_k}$.
        By statement~\ref{claim-multiple-isolated-vertices}, $s_j$ is isolated in $\mathcal{G}'_{X_k}$ and hence in $\mathcal{G}_{X_k}$.

        The second case is when $\{s_i,s_{i+1}\}$ and $\{s_{\ell},s_1\}$ are missing in $\mathcal{G}_{X_k}$ for some $1\leq i<\ell$.
        Then $\{s_i,s_{i}+1\}$ and $\{s_{\max},0\}$ are missing in $\mathcal{G}'_{X_k}$ and by statement~\ref{claim-disconnected-sets}, there are no edges between $\{0,\ldots,s_i\}$ and $\{s_i +1,\ldots,s_{\max}\}$ in $\mathcal{G}'_{X_k}$ and therefore no edges between $\{s_1,\ldots,s_i\}$ and $\{s_{i+1},\ldots,s_{\ell}\}$ in $\mathcal{G}_{X_k}$, completing the proof.

        Proof of statement~\ref{statement-not-adjacent-before}: By assumption $g(\{i,i+1\})<k$.
        Because
        \[
            g(\{i,i+1\})-g(\{0,i+1\})= (n-r)i\geq 0,
        \]
        we also have $g(\{0,i+1\})<k$.
        Now $g(\{s,i+1\})$ is convex in $s$ when $s<i+1$, thus $g(\{s,i+1\})<k$ for $0\leq s\leq i$, i.e.\ $s\not\sim i+1$.

        Proof of statement~\ref{statement-not-adjacent-after}: If $i+1=s_{\max}$, there is nothing to prove, so let $i+2\leq s_{\max}$.
        In the symmetric case $s_{\max}\leq \frac{r-1}{2}$, and we have $i+2\leq\frac{r-1}{2}$.
        Therefore,
        \[
            \Delta=g(\{i,i+1\})-g\left(\{i,\frac{r-1}{2}\}\right)=\frac{1}{8}\left( r-2i-5 \right)\left( r-2i-3 \right)\geq 0,
        \]
        which together with $g(\{i,i+1\})<k$ and the convexity of the edge label function $g$ in each variable yields $g(\{i,s\})<k$ for $i+1\leq s\leq\lfloor\frac{r-1}{2}\rfloor$, i.e.\ $i\not\sim s$.
        In the alternating and rectangular cases, $s_{\max}=\frac{r-2}{2}$ and $s_{\max}=r-1$ respectively, and a similar argument using $\Delta=g(\{i,i+1\})-g(\{i,s_{\max}\})$ goes through.
    \end{proof}

    \section{On the components}\label{sec:on-the-components}
    The rest of this paper focuses solely on the Fano schemes $\fano{k}{\SD[r]{n}}$.
    In this section, we prove in particular cases that the closed subscheme $\compk{k}{s}$ is equal to an irreducible component of $\fano{k}{\SD[r]{n}}$ taken with the reduced subscheme structure.
    We expect this statement to hold for every $k$ (see Conjecture~\ref{conj:symmetric-ck-component}).
    Table~\ref{tab:summary-symmetric-components} lists all cases for which we prove the statement in this paper.

    \begin{table}
        \centering
        \begin{tabular}[t]{ c c l }
            \toprule
            $s$ & $k$ & Statement\\
            \midrule
            $0$ & $\kappa(0)-(r-2)<k\leq \kappa(0)$ & Corollary~\ref{cor:component-s=0}\\
            \addlinespace
            $1$ & $k=\kappa(1)$  & Corollary~\ref{cor:component-s>=1} \\
            \addlinespace
            $1$ & $k=\kappa(1)-1$ when $n>6$, $r=6$ & Proposition~\ref{prop:s=1-special-case-component}\\
            \addlinespace
            $\geq 2$ & $\kappa(s)-1\leq k\leq \kappa(s)$ & Corollary~\ref{cor:component-s>=1} \\
            \addlinespace
            $\frac{r-1}{2}$ & arbitrary $k$ & Proposition~\ref{prop:smooth-points}\\
            \addlinespace
            arbitrary $s$ & 1 & Theorem~\ref{thm:fano-lines-components-complete}\\
            \bottomrule
        \end{tabular}
        \caption{Cases where we show $\compk{k}{s}$ is an irreducible component of $\fano{k}{\SD[r]{n}}$ taken with the reduced subscheme structure}
        \label{tab:summary-symmetric-components}
    \end{table}

    Our general strategy to prove $\compk{k}{s}$ gives a component is the following.
    We choose a suitable Borel fixed nested $s$-compression space $[Q]\in\compk{k}{s}$.
    By Proposition~\ref{prop:characterize-borel-fixed-points-grassmannian}, we may assume every entry on and above the diagonal of the matrix $Q$ is a free variable or zero.
    There is an affine chart of the Grassmannian $\Gr\left( k+1,\sym_n \right)$, called the \textit{distinguished} affine chart,  corresponding to $[Q]$ consisting of all points $[P]$ where the matrix $P$ is obtained from $Q$ by replacing its zero entries with linear forms in its free variables.
    For every point $[P]$ of this neighborhood lying on $\fano{k}{\SD[r]{n}}$, we partition the matrix $P$ into blocks as
    \begin{equation}\label{eq:partition-typical-point-affine-chart}
    \begin{pNiceArray}{cccc}[first-row,first-col,margin]
        & \outmat{s}       & \outmat{r-2s-1} & \outmat{n-r+1} & \outmat{s}\\
        \outmat{s}      &  Z_1    &  Z_2   &  Z_3  & Z_4\\
        &         &        &       &    \\
        \outmat{r-2s-1} &   *  &  Z_5   &   C_1 & C_2\\
        &         &        &       &    \\
        \outmat{n-r+1}  &   *  &   * &   C_3 & C_4\\
        &         &        &       &    \\
        \outmat{s}      &   *  &   * &  * & C_5
    \end{pNiceArray},
    \end{equation}
    and perform a particular sequence of \rowcol{} operation, called \textit{simplifying} $P$, to reduce $C_i$ blocks to zero.
    This will show $[P]$ is a nested $s$-compression space, i.e.\ $[P]\in\compk{k}{s}$ implying that $\compk{k}{s}$ is an irreducible component of $\fano{k}{\SD[r]{n}}$ taken with the reduced subscheme structure.

    \begin{definition}[simplifying matrices]\label{def:symmetric-simplification}
    Let $[Q]\in\fano{k}{\SD[r]{n}}$ be a Borel-fixed nested $s$-compression space, and assume $[P]$ is a point in the distinguished affine chart corresponding to $[Q]$.
    Partition $[P]$ as in~\eqref{eq:partition-typical-point-affine-chart}, where exactly $k+1$ entries in the $Z_i$ blocks are independent variables $z_{i,j}$ and the rest of the entries in the $Z_i$ blocks along with all the entries in the $C_i$ blocks are linear forms in the $z_{i,j}$.
    We call the entries on the anti-diagonal of $Z_4$ and $Z_5$ \textit{marked entries} (see Figure~\ref{fig:marked-entries}).

    When the marked entries of $Q$ (and hence every entry above and to the left of them) are independent variables $z_{i,j}$, we call the following sequence of \rowcol{} operations \textit{simplifying} the matrix $P$:
    Label the marked entries as $1,2,\ldots,r-s-1$ starting from the bottom left and going up and to the right (the first one is the $(r-s-1,s+1)$ entry, see Figure~\ref{fig:marked-entries}).
    For each $i=1,\ldots,r-s-1$, perform row-column operations to remove marked entry with label $i$ from the $c_{i,j}$ entries below it in that column.
    After each \rowcol{} operation, via a change of variables and relabeling, we may assume these operations do not change any of the $z_{i,j}$ entries, and we do not disturb the elimination we did for the previous marked entries.
    We call the resulting matrix \textit{simplified}.
    \end{definition}

    \begin{figure}
        \centering
        \begin{subfigure}[b]{.4\textwidth}
            \centering
            $
            \begin{pmatrix}
                z_{1,1}        & \cdots & \bm{z_{1,n-1}} & c_{1,n}   \\
                \vdots         & \bm{\iddots} & \vdots         & \vdots    \\
                \bm{z_{1,n-1}} & \cdots & z_{n-1,n-1}    & c_{n-1,n} \\
                c_{1,n}        & \cdots & c_{n-1,n}      & c_{n,n}
            \end{pmatrix}
            $
            \caption{$s=0$, $k=\kappa(0)$}\label{fig:marked-entries-s=0}
        \end{subfigure}
        \quad\quad
        \begin{subfigure}[b]{.4\textwidth}
            \centering
            $
            \begin{pmatrix}
                z_{1,1} & z_{1,2} & z_{1,3} & z_{1,4} &\bm{z_{1,5}} \\
                z_{1,2} & z_{2,2} & \bm{z_{2,3}} & c_{2,4}        &  c_{2,5}       \\
                z_{1,3} & \bm{z_{2,3}} & z_{3,3} &  c_{3,4}       &  c_{3,5}       \\
                z_{1,4} & c_{2,4}        &  c_{3,4}       &   c_{4,4}      &  c_{4,5}       \\
                z_{1,5}   &  c_{2,5}       &   c_{3,5}      &  c_{4,5}       &c_{5,5}
            \end{pmatrix}
            $
            \caption{$n=5$, $s=1$, $k=\kappa(1)=7$}\label{fig:marked-entries-s=1}
        \end{subfigure}
        \caption{Marked entries in Definition~\ref{def:symmetric-simplification} are given in bold (when $r=n$); the $c_{i,j}$ are linear forms in the free variables $z_{i,j}$.}
        \label{fig:marked-entries}
    \end{figure}

    \begin{lemma}\label{lem:from-n-to-r}
    Fix $0\leq s\leq \lfloor\frac{r-1}{2}\rfloor$.
    Suppose $[P]\in\Gr\left( k+1,\sym_n \right)$ is given by a symmetric matrix $P$ of the block form~\eqref{eq:partition-typical-point-affine-chart} with entries in $\KK[z_0,\ldots,z_k]_1$.
    Assume $P$ satisfies the following conditions
    \begin{enumerate}
        \item $\det Z_4\neq 0$ when $s\neq 0$ and $\det Z_5\neq0$ when $r-2s-1\neq 0$, \label{part:detz4-z5-nonzero}
        \item the $r\times r$ minors of $P$ vanish, \label{part:every-rbyr-vanishes}
        \item the vanishing of a principal $r \times r$ minor of $P$ involving the first $r-s-1$ and the last $s$ rows and columns implies those entries of the minor belonging to the $C_i$ blocks are zero. \label{part:every-principal-gives-zero-cij}
    \end{enumerate}
    Then $C_i=0$ for $1\leq i\leq 5$.
    \end{lemma}
    \begin{proof}
        For each $r-s\leq \ell \leq n-s$, consider the principal $r\times r $ minor of $P$ involving rows and columns $\{1,\ldots,r-s-1,\ell,n-s+1,\ldots,n\}$.
        By condition~\ref{part:every-rbyr-vanishes}, this minor vanishes and by condition~\ref{part:every-principal-gives-zero-cij} the $c_{i,j}$ entries in this minor are zero.
        This means all entries in $C_{1}, C_2,C_4,C_5$ and the diagonal entries in $C_3$ are zero.

        To show $C_3=0$, let $c_{i,j}$ be the $(i,j)$-entry of $P$ located in the block $C_3$.
        Then consider the $r\times r$ minor given by the first $r-s-1$ rows and columns, the last $s$ rows and columns, and the $i$-th row and $j$-th column.
        This minor is $\pm c_{i,j}(\det Z_4)^2\det Z_5$ which is zero by condition~\ref{part:every-rbyr-vanishes}.
        Using condition~\ref{part:detz4-z5-nonzero} gives $c_{i,j}=0$ as desired.
    \end{proof}

    The next two lemmas pave the way to proving the main result of this section.

    \begin{lemma}\label{lem:simplified-s=0}
    Let $s=0$, and $\kappa(0)-(r-2)< k\leq \kappa(0)$.
    Assume $[Q_n]\in\fano{k}{\SD[r]{n}}$ is the Borel fixed point obtained by taking the standard 0-compression space in Figure~\ref{fig:standard-0-comp-space} and setting $z_{r-m,r-1}, \ldots, z_{r-1,r-1}$ equal to zero where $m=\kappa(0)-k$.
    Suppose $P_n$ is a symmetric matrix obtained from $Q_n$ by replacing the zero entries of $Q_n$ with linear forms $c_{i,j}$ in the variables $z_{i,j}$.
    If $P_n$ is partitioned with $C_i$ blocks as in~\eqref{eq:partition-typical-point-affine-chart} and the $r\times r$ minors of $P_n$ vanish, then simplifying $P_n$ reduces all $C_i$ blocks to 0, i.e.\ $P_n$ is a nested $0$-compression space.
    \end{lemma}
    \begin{proof}
        \begin{figure}
            \centering
            \begin{subfigure}[b]{.4\textwidth}
                \centering
                $
                \left(\begin{array}{ccc|ccc}
                          z_{1,1}        & \cdots & \bm{z_{1,r-1}} &&  & \\
                          \vdots         & \bm{\iddots} & \vdots         &&  \Large{0}&   \\
                          \bm{z_{1,r-1}} & \cdots & z_{r-1,r-1}    &&&  \\
                          \hline
                          &\Large{0} & &    & \Large{0}  &
                \end{array}\right)_{n\times n}
                $
                \caption{Standard 0-compression space}\label{fig:standard-0-comp-space}
            \end{subfigure}
            \quad\quad
            \begin{subfigure}[b]{.4\textwidth}
                \centering
                $
                \left(\begin{array}{ccc|c}
                          z_{1,1}        & \cdots & \bm{z_{1,n-1}} &c_{1,n} \\
                          \vdots       & & \vdots         & \vdots \\
                          \vdots&&\makecell{z_{n-m-1,n-1}\\c_{n-m,n-1}}&\vdots\\
                          \vdots &&\vdots&\vdots\\
                          \bm{z_{1,n-1}} & \cdots & c_{n-1,n-1}    &c_{n-1,n}  \\
                          \hline
                          c_{1,n}&\cdots &c_{n-1,n}  & c_{n,n}
                \end{array}\right)_{n\times n}
                $
                \caption{The matrix $P_n$ when $r=n$}\label{fig:typical-point-s=0}
            \end{subfigure}
            \caption{Marked entries in Lemma~\ref{lem:simplified-s=0} shown in bold}
            \label{fig:matrices-in-prop-for-s=0}
        \end{figure}

        We consider the case $r=n$ and induct on $n$.
        See Figure~\ref{fig:typical-point-s=0} for the matrix $P_n$.
        To be able to use the induction hypothesis we change the claim slightly and assume that $c_{i,j}$ are linear forms in the variables $z_{i,j}$ and potentially other independent variables.
        It is also convenient to take $n=2$ as the base case,
        \[
            P_2=\begin{pmatrix}
                    \bm{z_{1,1}} & c_{1,2} \\
                    c_{1,2} & c_{2,2}
            \end{pmatrix}.
        \]
        After simplifying $P_2$, i.e.\ removing the marked entry $z_{1,1}$ from $c_{1,2}$ by \rowcol{} operation, $\det P_2=0$ implies $c_{1,2}=c_{2,2}=0$, establishing the base of induction.

        Assume $n>2$ and the claim holds for simplified matrices of size smaller than $n$.
        Modulo $z_{1,1},z_{1,2},\ldots,z_{1,n-1}$, we have $0=\det P_{n}=\pm c_{1,n}^2\det R_1$, where $R_1$ is the submatrix obtained by deleting rows and columns $\{1,n\}$ in $P_{n}$.
        We claim $\det R_1\neq 0$ which then gives $c_{1,n}=0$.

        To prove $\det R_1\neq 0$, we use $\kappa(0)-(r-2)< k$ to get $m=\kappa(0)-k< r-2=n-2$.
        This means the first two entries in column $n-1$ of $P_n$ are independent variables $z_{1,n-1}$ and $z_{2,n-1}$.
        Therefore, the entries on and above the anti-diagonal of $R_1$ are distinct independent variables and after specializing these variables to 1 and the other independent variables to 0, we get $\det R_1\neq 0$.
        We conclude that $c_{1,n}=0$ modulo $z_{1,1},z_{1,2},\ldots,z_{1,n-1}$, i.e.\ $c_{1,n}$ involves only the $z_{i,j}$ in the first row except $z_{1,n-1}$ (in simplifying $P_{n}$ we removed $z_{1,n-1}$ from $c_{1,n}$).

        Modulo $z_{1,1}, z_{1,2},\ldots,z_{1,n-2}$, we have $c_{1,n}=0$ and
        \begin{equation}\label{eq:det-base-case}
        0=\det P_{n}=
        \pm z_{1,n-1}^2\det
        \begin{pmatrix}
            z_{2,2}   & \cdots & z_{2,n-2}   & c_{2,n}   \\
            \vdots    & \ddots & \vdots      & \vdots    \\
            z_{2,n-2} & \cdots & z_{n-2,n-2} & c_{n-2,n} \\
            c_{2,n}   & \cdots & c_{n-2,n}   & c_{n,n}
        \end{pmatrix},
        \end{equation}
        where the matrix on the right is simplified.
        Using the induction hypothesis gives that $c_{2,n}, \ldots, c_{n-2,n}, c_{n,n}$ are zero (modulo $z_{1,1},\ldots,z_{1,n-2}$) and thus do not involve $z_{1,n-1}$.
        Because $P_n$ is simplified, we conclude that $z_{1,n-1}$ appears only in entries $(1,n-1)$ and $(n-1,1)$.
        Therefore in the expansion of $\det P_{n}$, the coefficient of $z_{1,n}^2$ is exactly the determinant on the right-hand side of~\eqref{eq:det-base-case} which must vanish (this time not modulo any variables).
        Another application of the induction hypothesis gives $c_{2,n},\ldots,c_{n-2,n},c_{n,n}$ are indeed zero.

        It remains to show $c_{1,n}=c_{n-1,n}=0$.
        If a term in the expansion of $\det P_{n}$ involves the entry $(1,n-1)$, then due to the zeros in the last column and row, this term has to involve $c_{n-1,n}$ in the last column and $c_{1,n}$ in the last row.
        Combining all such terms gives $\pm z_{1,n-1}c_{n-1,n}c_{1,n}\det R_2$ where $R_2$ is obtained from $P_{n}$ by deleting rows and columns $\{1,n-1,n\}$.
        The coefficient of $z_{1,n-1}$ in $\det P_{n}$ is thus $\pm 2z_{1,n-1}c_{n-1,n}c_{1,n}\det R_2$.
        Because $\det R_2\neq 0$, we get $c_{n-1,n}c_{1,n}=0$.
        Since one of $c_{n-1,n}$ or $c_{1,n}$ is zero, by expanding $\det P_{n}$ along the last column, we get the other one is also zero.
        This completes the proof for $r=n$.

        When $r<n$, $P_n$ satisfies the conditions in Lemma~\ref{lem:from-n-to-r}.
        condition~\ref{part:every-principal-gives-zero-cij} of the lemma holds for $P_n$ by the case $r=n$ proved above and because a principal $r\times r$ submatrix of $P_n$ involving the first $r-1$ rows and columns is simplified.
        By the lemma, all $c_{i,j}$ not in the upper-left $(r-1)\times(r-1)$ block of $P_n$ are zero, as required.
    \end{proof}

    \begin{lemma}\label{lem:simplified-s>=1}
    Let $1\leq s\leq \lfloor\frac{r-1}{2}\rfloor$.
    Assume $[Q_{n,s}]\in\fano{\kappa(s)}{\SD[r]{n}}$ is the $n\times n$ standard $s$-compression space in~\eqref{eq:standard-zero-pattern} where the star blocks are filled with $\kappa(s)+1$ independent variables $z_{i,j}$
    Suppose the symmetric matrix $P_{n,s}$ is obtained from $Q_{n,s}$ by replacing the zero entries of $Q_{n,s}$ with linear forms $c_{i,j}$ in the variables $z_{i,j}$.
    If $P_{n,s}$ is partitioned with $C_i$ blocks as in~\eqref{eq:partition-typical-point-affine-chart} and the $r\times r$ minors of $P_{n,s}$ vanish, then simplifying $P_{n,s}$ reduces all $C_i$ blocks to 0, i.e.\ $P_{n,s}$ is a nested $s$-compression space.
    Moreover, the statement holds when $s\geq 2$, and we set $z_{s,n}=0$ in $Q_{n,s}$.
    \end{lemma}

    \begin{proof}
        We prove the claim first for $r=n$ via inducting on $n$.
        Call the marked entries $(1,n)$ and $(s,n-s+1)$ in $P_{n,s}$ \textit{pivots}  (see Figure~\ref{fig:pivots-s}).
        Consider the pivot $z_{1,n}$.
        Modulo the variables in the first row of $P_{n,s}$ except the pivot we have
        \[
            0=\det P_{n,s}=\pm z_{1,n}^2\det P_{n-2,s-1},
        \]
        where $P_{n-2,s-1}$ is the matrix obtained from $P_{n,s}$ by deleting rows and columns $\{1,n\}$.
        Because $P_{n-2,s-1}$ is simplified, if $s-1\geq 1$, by the induction hypothesis, its $c_{i,j}$ entries are zero modulo the variables in the first row of $P_{n,s}$ except the pivot.
        If $s-1=0$, Lemma~\ref{lem:simplified-s=0} gives that all $c_{i,j}$ in $P_{n-2,s-1}$ are zero.

        This shows that the $c_{i,j}$ entries in $P_{n-2,s-1}$ do not involve the pivot $z_{1,n}$ and since $P_{n,s}$ is simplified, the only entries involving $z_{1,n}$ are the entries $(1,n)$ and $(n,1)$.
        Hence, the coefficient of $z_{1,n}^2$ in $\det P_{n,s}$ is $\det P_{n-2,s-1}$ which must be zero (this time not modulo any variables) and another application of the induction hypothesis (or Lemma~\ref{lem:simplified-s=0} when $s-1=0$) gives that the $c_{i,j}$ entries in $P_{n-2,s-1}$ are indeed zero (see Figure~\ref{fig:ell=2} for an example).

        \begin{figure}
            \centering
            \begin{subfigure}[b]{.45\textwidth}
                \centering
                $
                \begin{pmatrix}
                    \cdots & z_{1,n-s}   & z_{1,n-s+1}      &              & \bm{z_{1,n}} \\
                    &    \vdots         &       \vdots           & \iddots &     \vdots         \\
                    \cdots  & z_{s,n-s}   & \bm{z_{s,n-s+1}} &              &     z_{s,n}         \\
                    \cdots & c_{s+1,n-s} & c_{s+1,n-s+1}      &              & c_{s+1,n}      \\
                    &  \vdots     &    \vdots              &              & \vdots
                \end{pmatrix}
                $
                \caption{The two pivots in $P_{n,s}$ (in bold)}\label{fig:pivots-s}
            \end{subfigure}
            \quad\quad
            \begin{subfigure}[b]{.45\textwidth}
                \centering
                $
                \begin{pmatrix}
                    \cdots & z_{1,n-2} & z_{1,n-1}      & \bm{z_{1,n}} \\
                    \cdots & z_{2,n-2} & \bm{z_{2,n-1}} & z_{2,n}      \\
                    & 0         & 0              & c_{3,n}      \\
                    & 0         & 0              & \vdots       \\
                    & 0         & 0              & c_{n-1,n}    \\
                    \cdots & c_{n-2,n} & c_{n-1,n}      & c_{n,n}
                \end{pmatrix}
                $
                \caption{Using pivot $(1,n)$ in $P_{n,2}$}\label{fig:ell=2}
            \end{subfigure}
            \caption{Pivots in the proof of Lemma~\ref{lem:simplified-s>=1}}
            \label{fig:making-zeros}
        \end{figure}

        The above argument shows that those $c_{i,j}$ in $P_{n,s}$ not in the last column and not in the last row are zero.
        To show $c_{i,n}=0$ for $i>s$, we consider two cases:

        The first case is $s\geq 2$.
        We then have a second pivot $(s,n-s+1)$ as in Figure~\ref{fig:making-zeros}.
        We show $c_{n-s+1,n}$ depends only on the $z_{i,j}$ variables in column $n-s+1$, except the pivot.
        Modulo the $z_{i,j}$ variables in column $n-s+1$, we have $0=\det P_{n,s}=\pm c_{n-s+1,n}^2\det R_3$, where $R_3$ is the matrix obtained by deleting rows and columns $\{n-s+1,n\}$ in $P_{n,s}$.
        Because $\det R_3\neq 0$, we have $c_{n-s+1,n}=0$, proving the claim for $c_{n-s+1,n}$.
        Recall that $c_{n-s+1,n}$ does not involve the pivot $z_{s,n-s+1}$ due to $P_{n,s}$ being simplified.

        Modulo the $z_{i,j}$ variables in column $n-s+1$ except the pivot we have $0=\det P_{n,s}=\pm z_{s,n-s+1}^2 \det R_4$, where $R_4$ is the simplified matrix of the form $P_{n-2,s-1}$ obtained by deleting rows and columns $\{s,n-s+1\}$ in $P_{n,s}$.
        Applying the induction hypothesis (or Lemma~\ref{lem:simplified-s=0} when $s-1=0$) to $\det R_4=0$, we get all the $c_{i,j}$ in $R_4$ are zero (modulo the mentioned variables), and hence these $c_{i,j}$ do not depend on the pivot $z_{s,n-s+1}$.

        We conclude that the only entries in $P_{n,s}$ involving $z_{s,n-s+1}$ are $(s,n-s+1)$ and $(n-s+1,s)$.
        Therefore the coefficient of $z_{s,n-s+1}^2$ in $\det P_{n,s}$ is $\pm \det R_4$ giving $\det R_4=0$ (this time not modulo any variables).
        Applying the induction hypothesis (or Lemma~\ref{lem:simplified-s=0} when $s-1=0$) to $R_4$ gives all $c_{i,j}$ in $R_4$ are zero.

        It remains to prove $c_{n-s+1,n}=0$.
        We consider those terms of $\det P_{n,s}$ involving $z_{1,n}$.
        Since $P_{n,s}$ is simplified, the only entries involving $z_{1,n}$ are $(1,n)$ and $(n,1)$.
        A term of $\det P_{n,s}$ is formed by picking one entry from each column and row.
        Assume a term in the expansion of $\det P_{n,s}$ involves entry $(1,n)$.
        Due to the zero entries of $P_{n,s}$, that term has to pick the entry $(n,n-s+1)$ in column $n-s+1$.
        Otherwise, if it picks an entry $(\ell',n-s+1)$ for some $2\leq \ell'\leq s$, then after deleting rows $\{1,\ell'\}$ and columns $\{n-s+1,n\}$ we are left with a square matrix of size $n-2$ that is a compression space, i.e.\ the term is zero.
        Thus, combining the terms involving $z_{1,n}$ with exponent 1 gives $\pm 2z_{1,n}c_{n-s+1,n}\det R_5$ where $R_5$ is the submatrix of $P_{n,s}$ obtained by deleting rows $\{1,n\}$ and columns $\{n-s+1,n\}$.
        Because $\det R_5\neq0$, we have $c_{n-s+1,n}=0$.

        The second case is when $s=1$.
        Then we have $z_{1,n}$ as the only pivot.
        We claim $c_{n-1,n}=0$.
        Since none of the $c_{i,j}$ involve $z_{1,n}$, the terms in $\det P_{n,s}$ involving $z_{1,n}$ with exponent 1 are $\pm 2z_{1,n}c_{n-1,n}\det R_6$, where $R_6$ is obtained from $P_{n,s}$ by deleting rows $\{1,n\}$ and columns $\{n-1, n\}$.
        Using the fact that $\det R_6\neq 0$, we get $c_{n-1,n}=0$.
        Expanding $\det P_{n,s}$ along column $n-1$ gives
        \[
            0=\det P_{n,s}=\pm z_{1,n-1}^2 \det
            \begin{pmatrix}
                z_{2,2}   & \cdots & z_{2,n-2}   & c_{2,n}   \\
                \vdots    & \ddots & \vdots      & \vdots    \\
                z_{2,n-2} & \cdots & z_{n-2,n-2} & c_{n-2,n} \\
                c_{2,n}   & \cdots & c_{n-2,n}   & c_{n,n}
            \end{pmatrix}.
        \]
        Since the above matrix is simplified, Lemma~\ref{lem:simplified-s=0} gives $c_{i,n}=0$.
        This completes the argument for $r=n$.

        When $r<n$, the matrix $P_{n,s}$ satisfies the conditions in Lemma~\ref{lem:from-n-to-r} (Condition~\ref{part:every-principal-gives-zero-cij} is simply what we proved above in the case $r=n$).
        The lemma now gives that $C_i=0$ for every $i$.
    \end{proof}

    \begin{corollary}\label{cor:component-s=0}\label{cor:component-s>=1}
    For arbitrary $n,r$, the subscheme $\compk{k}{s}$ of $\fano{k}{\SD[r]{n}}$ is an irreducible component taken with the reduced subscheme structure in each of the following cases.
    \begin{enumerate}[label=(\roman*)]
        \item $s=0$, and $\kappa(0)-(r-2)< k\leq \kappa(0)$,\label{part:component-part-s=0}
        \item $s=1$, and $k=\kappa(s)$,\label{part:component-part-s=1}
        \item $2\leq s\leq \lfloor\frac{r-1}{2}\rfloor$ and $\kappa(s)-1\leq k\leq \kappa(s)$.\label{part:component-part-s>=2}
    \end{enumerate}
    \end{corollary}
    \begin{proof}
        Consider part~\ref{part:component-part-s=0}.
        By Lemma~\ref{lem:simplified-s=0} there is a Borel fixed point $[Q_n]\in\compk{k}{0}$ such that the distinguished affine chart of $\fano{k}{\SD[r]{n}}$ corresponding to $[Q_n]$ is fully contained in $\compk{k}{0}$.
        This implies $\compk{k}{0}$ is an irreducible component.
        The proofs for parts~\ref{part:component-part-s=1} and~\ref{part:component-part-s>=2} are similar by using Lemma~\ref{lem:simplified-s>=1} instead.
    \end{proof}

    When $k=\kappa(s)$, the components $\compk{\kappa(s)}{s}$ are familiar objects as we now describe.
    \begin{corollary}
        \label{cor:kappa-s-component}
        For $0\leq s\leq\lfloor\frac{r-1}{2}\rfloor$, the closed subscheme $\comp{s}$ is an irreducible component of $\fano{\kappa(s)}{\SD[r]{n}}$ taken with the reduced subscheme structure, and is isomorphic to the flag variety $\flag{s+1+n-r,n-s}{\KK^n}$.
        In particular, $\comp{s}$ is the $\gl{n}$ orbit of the standard $s$-compression space and does not intersect other components (if any) of the scheme.
    \end{corollary}
    \begin{proof}

        By Corollary~\ref{cor:component-s=0}, $\comp{s}$ is an irreducible component of $\fano{\kappa(s)}{\SD[r]{n}}$ taken with the reduced subscheme structure.
        The  subscheme $\compk{\kappa(s)}{s}$ is the image of the morphism in~\eqref{eq:map-from-flag-variety}.
        Denote by $\phi$ the composition of $\fano{\kappa(s)}{\SD[r]{n}}\to\Gr(\kappa(s),\sym_n)$ with the map in~\eqref{eq:map-from-flag-variety}.
        Then $\phi$ is $\gl{n}$-equivariant and bijective onto its image.

        A standard argument using affine charts of the flag variety and the Grassmannian shows the differential $d\phi$ is injective at the flag that maps to the standard $s$-compression space.
        Because $\gl{n}$ acts transitively on $\flag{s+1+n-r,n-s}{\KK^n}$, $d\phi$ is injective everywhere, hence $\phi$ is a closed embedding.

        This implies $\comp{s}$ is isomorphic to $\flag{s+1+n-r,n-s}{\KK^n}$.
        In particular, it is the $\gl{n}$-orbit of any single point such as the standard $s$-compression space.
        Finally, if $\compk{\kappa(s)}{s}$ intersects another irreducible component, then the intersection would be $\gl{n}$-stable and $\compk{\kappa(s)}{s}$ will have more than one orbit, a contradiction.
    \end{proof}

    As a corollary, we recover Loewy and Radwan's theorem~\cite[Theorem 6.1]{loewy-radwan}.

    \begin{corollary}
        \label{cor:loewy-radwan}
        Let $Q\subseteq \sym_n$ be a subspace whose elements have rank less than $r$ and assume $Q$ is of maximal dimension $\max\left\{\kappa(0),\kappa\left(\lfloor\frac{r-1}{2}\rfloor\right)\right\}+1$.
        Then $Q$ is a $\gl{n}$-translate of the standard $s$-compression space for $s=0$ or $\lfloor\frac{r-1}{2}\rfloor$.
    \end{corollary}
    \begin{proof}
        Let $k=\max \{\kappa(0),\kappa(s_{\max})\}$ with $s_{\max}=\lfloor\frac{r-1}{2}\rfloor$, and consider the Fano scheme $\fano{k}{\SD[r]{n}}$.
        If $k=\kappa(0)$ or $k=\kappa(s_{\max})$, then Corollary~\ref{cor:kappa-s-component} gives that $\compk{k}{0}$, respectively $\compk{k}{s_{\max}}$, is an irreducible component of $\fano{k}{\SD[r]{n}}$ (taken with the reduced subscheme structure) not intersecting other components and is the $\gl{n}$-orbit of the standard $0$-compression space, respectively $\left(s_{\max}\right)$-compression space.

        Since $\kappa(s)$ is convex in $s$, for $0< s <s_{\max}$ we have $k>\kappa(s)$ and thus every Borel fixed point is a nested $0$- or $\left(s_{\max}\right)$-compression space.
        By Borel's theorem~\cite[Theorem 10.4]{borel}, every irreducible component contains a Borel fixed point.
        Therefore, the Fano scheme $\fano{k}{\SD[r]{n}}$ is the disjoint union of $\compk{k}{0}$ and $\compk{k}{s_{\max}}$.
    \end{proof}

    We recall a result of de Seguins Pazzis that generalizes Corollary~\ref{cor:loewy-radwan}.
    \begin{theorem}[{\cite[Theorem 1.3]{pazzis-classification}}]\label{thm:pazzis}
        Let $k>\max\{\kappa(1),\kappa(s_{\max}-1)\}$ where $s_{\max}=\lfloor\frac{r-1}{2}\rfloor$.
        Then every closed point of $\fano{k}{\SD[r]{n}}$ is a nested $s$-compression space for $s=0$ or $s_{\max}$.
    \end{theorem}

    In Lemma~\ref{lem:simplified-s>=1} and Corollary~\ref{cor:component-s>=1}, we excluded the case $s=1$ and $k=\kappa(s)-1$ because when $n=r$, it often happens that $\compk{k}{1}\cap\compk{k}{0}\neq \emptyset$.
    The affine chart corresponding to a Borel fixed point on this intersection contains both nested $0$- and $1$-compression spaces and our technique fails.
    The next proposition deals with a case when this intersection is empty.
    We will use this result in §\ref{sec:smoothness} when characterizing the smoothness of our Fano schemes $\fano{k}{\SD[r]{n}}$.

    \begin{proposition}\label{prop:s=1-special-case-component}
    Let $n>6$, $r=6$, $k=\kappa(1)-1$.
    Then $\compk{k}{1}$ is an irreducible component of $\fano{k}{\SD[r]{n}}$ taken with the reduced subscheme structure.
    \end{proposition}
    \begin{proof}
        Let $s=1$ and assume $Q\subseteq \sym_n$ is the Borel fixed point obtained by taking the $n\times n$ standard  $1$-compression space in~\eqref{eq:standard-zero-pattern} where the star blocks are filled with $\kappa(1)+1$ independent variables $z_{i,j}$ and setting $z_{1,n}=0$.
        Take a point $[P]$ in the distinguished affine chart of $\fano{k}{\SD[r]{n}}$ corresponding to $[Q]$.
        We show $P$ is a nested 1-compression space proving $\compk{k}{1}$ gives a component.
        The matrix $P$ is of the form
        \vspace{3mm}
        \[
            P=\begin{pNiceMatrix}[left-margin=0.6em,right-margin=0.6em]
                  \Block[draw]{3-5}{} &              &&& & z_{1,5} & \cdots     & z_{1,n-1} & c_{1,n} \\
                  &\bm{z_{i,j}}  &&& &           &            &           &         \\
                  &              &&& &           &            &           &         \\
                  z_{1,5}           &              &&& &           &            &           &         \\
                  \vdots              &              &&& &           &\bm{c_{i,j}}&           &         \\
                  z_{1,n-1}           &              &&& &           &            &           &         \\
                  c_{1,n}             &              &&& &           &            &           &
                  \CodeAfter
                  \OverBrace[shorten,yshift=3pt]{1-1}{2-5}{\outmat{r-s-1=4}}
            \end{pNiceMatrix},
        \]
        where the $c_{i,j}$ are linear forms in the $z_{i,j}$ and the $r\times r$ minors of $P$ vanish.

        Let $P_1$ be the upper left $(n-1)\times (n-1)$ submatrix of $P$.
        Then $[P_1]\in\fano{k'}{\SD[r]{{n-1}}}$, where $k'=\kappa(1,n-1,r)$.
        By Lemma~\ref{lem:simplified-s>=1}, we can reduce the $c_{i,j}$ entries in $P_1$ to zero by \rowcol{} operations on $P$.
        It remains to reduce the linear forms in the last row and column of $P$, except possibly $c_{1,n}$, to zero.

        Eliminate the variable in each of the entries $(4,2),(3,3),(2,4),(1,n-1)$ from the $c_{i,j}$ below it in the last row of $P$.
        We claim this reduces $c_{2,n},\ldots,c_{n,n}$ to zero.
        To see this, consider the $6\times 6$ submatrix $P_2$ given by rows and columns $\{1,\ldots,4,n-1,n\}$ of $P$ (see Figure~\ref{subfig:P2}).

        Modulo $z_{1,n-1}$, we have $0= \det P_2=\pm c_{n-1,n}^2\det P_2'$ where $P_2'$ is obtained by deleting rows and columns $\{5,6\}$ in $P_2$.
        Because $\det P_2'\neq 0$, we have $c_{n-1,n}=0$, i.e.\ $c_{n-1,n}$ depends only on $z_{1,n-1}$.
        But we eliminated $z_{1,n-1}$ from $c_{n-1,n}$, which gives $c_{n-1,n}=0$ (not modulo any variables).

        Therefore $0=\det P_2=\pm z_{1,n-1}^2\det P_2''$, where $P_2''$ is obtained by deleting rows and columns $\{1,5\}$ in $P_2$.
        Since $\det P_2''= 0$ and $P_2''$ is simplified, by Lemma~\ref{lem:simplified-s=0}, $c_{2,n}=c_{3,n}=c_{4,n}=c_{n,n}=0$.

        Next, fix $5\leq \ell < n-1$, and consider the $6\times 6$ submatrix $P_3$ given by rows and columns $\{1,\ldots,4,\ell,n\}$ of $P$ (see Figure~\ref{subfig:P3}).
        A similar argument shows that $c_{\ell,n}$ depends only on $z_{1,\ell}$, say $c_{\ell,n}=\alpha z_{1,\ell}$.
        Expand $\det P_3$ along the last column to get
        \[
            0=\det P_3= 2\alpha z_{1,\ell}^2 c_{1,n}\det
            \begin{pmatrix}
                z_{2,2} & z_{2,3} & z_{2,4}\\
                z_{2,3} & z_{3,3} & z_{3,4}\\
                z_{2,4} & z_{3,4} & z_{4,4}
            \end{pmatrix}
            -\alpha^2 z_{1,\ell}^2\det
            \begin{pmatrix}
                z_{1,1} & \cdots & z_{1,4}\\
                \vdots & \ddots & \vdots \\
                z_{1,4} & \cdots & z_{4,4}
            \end{pmatrix}.
        \]
        Because the determinant of a generic symmetric $4\times 4$ matrix is irreducible (see~\cite[Theorem 1]{kutz74}), we conclude $\alpha=0$, i.e.\ $c_{\ell,n}=0$.
        Since $\ell$ was arbitrary, we are done.

        \begin{figure}
            \centering
            \begin{subfigure}{.49\textwidth}
                \centering
                \hspace*{-.05\linewidth}
                $
                \begin{pmatrix}
                    z_{1,1} & z_{1,2} & z_{1,3} & z_{1,4} & z_{1,n-1} & c_{1,n}\\
                    z_{1,2} & z_{2,2} & z_{2,3} & z_{2,4} & 0 & c_{2,n}\\
                    z_{1,3} & z_{2,3} & z_{3,3} & z_{3,4} & 0 & c_{3,n}\\
                    z_{1,4} & z_{2,4} & z_{3,4} & z_{4,4} & 0 & c_{4,n}\\
                    z_{1,n-1} & 0 & 0 & 0 & 0 & c_{n-1,n}\\
                    c_{1,n} & c_{2,n} & c_{3,n} & c_{4,n} & c_{n-1,n} & c_{n,n}
                \end{pmatrix}
                $
                \caption{The matrix $P_2$}\label{subfig:P2}
            \end{subfigure}
            \begin{subfigure}{.49\textwidth}
                \centering
                \hspace*{.05\linewidth}
                $
                \begin{pmatrix}
                    z_{1,1} & z_{1,2} & z_{1,3} & z_{1,4} & z_{1,\ell} & c_{1,n}\\
                    z_{1,2} & z_{2,2} & z_{2,3} & z_{2,4} & 0 & 0\\
                    z_{1,3} & z_{2,3} & z_{3,3} & z_{3,4} & 0 & 0\\
                    z_{1,4} & z_{2,4} & z_{3,4} & z_{4,4} & 0 & 0\\
                    z_{1,\ell} & 0 & 0 & 0 & 0 & c_{\ell,n}\\
                    c_{1,n} & 0 & 0 & 0 & c_{\ell,n} & 0
                \end{pmatrix}
                $
                \caption{The matrix $P_3$}\label{subfig:P3}
            \end{subfigure}
            \caption{The matrices in the proof of Proposition~\ref{prop:s=1-special-case-component}}
        \end{figure}
    \end{proof}

    \section{Fano scheme of lines}\label{sec:fano-scheme-of-lines}
    In this section we give a complete description of the irreducible components of $\fano{1}{\SD[r]{n}}$.
    Using this description, we will characterize when $\fano{k}{\SD[r]{n}}$ is irreducible.
    To prove the main theorem of this section, we need several preliminary results.

    \begin{lemma}
        \label{lem:generic-matrix}
        Fix integers $n\geq2$ and $0<s\leq\frac{n}{2}$.
        For a general 2-dimensional linear subspace $Q\subseteq \sym_n$, there are no linear subspaces $U,W\subseteq\KK^n$ of dimensions $s$ and $n-s$ such that $U\subseteq W$ and $U\perp_Q W$.
    \end{lemma}
    \begin{proof}
        Define the incidence variety $Y\subseteq\flag{s,n-s}{\KK^n}\times\Gr(2,\sym_n)$ by
        \[
            Y=\left\{ \big( [U,W], [P] \big)\in\flag{s,n-s}{\KK^n}\times\Gr(2,\sym_n) \,|\, U\perp_{P} W\right\}.
        \]
        Consider the diagram
        \[
            \begin{tikzcd}
                Y \arrow[r, "\pi_1"] \arrow[d, "\pi_2"'] & \flag{s,n-s}{\KK^n} \\
                \Gr(2,\sym_n)                                      &
            \end{tikzcd}
            ,\]
        where $\pi_1$ and $\pi_2$ are the projections of $\flag{s,n-s}{\KK^n}\times\Gr(2,\sym_n)$ onto its factors.

        We show that the restriction of $\pi_2$ to $Y$ is not surjective by comparing the dimensions of $Y$ and $\Gr(2,\sym_n)$.
        Fix a point $[U,W]\in\flag{s,n-s}{\KK^n}$, where $U\subseteq W$ are subspaces of $\KK^n$ of dimensions $s$ and $n-s$ respectively.
        The $\gl{n}$ action on $\KK^n$ induces a natural action of $\gl{n}$ on $\flag{s,n-s}{\KK^n}$
        and on $Y$.
        The restriction of $\pi_1$ to $Y$ is $\gl{n}$-equivariant.
        We may thus assume
        \[
            U = \spn\{e_{n-s+1},\ldots,e_n\},\qquad W = \spn\{e_{s+1}, \ldots,e_n\},
        \]
        where $e_i$ denotes the $i$-th standard basis vector in $\KK^n$.
        The fiber of $\pi_1$ over $[U,W]$ consists of all 2-dimensional subspaces of $\tilde{P}\subseteq \sym_n$ where $\tilde{P}$ is the subspace of all symmetric matrices of the form
        \[
            \begin{pNiceArray}{cccccc}[first-row,first-col]
                &&\outmat{s} &&\outmat{n-2s}&\outmat{s}&\\
                \outmat{s}&&*&&*&*& \\
                & &  &&&& \\
                \outmat{n-2s}&&*&&*&0&\\
                & &&  &&& \\
                \outmat{s}&&*&&0&0&\\
            \end{pNiceArray}.
        \]
        Therefore,
        \[
                \dim\pi_1^{-1}([U,W]) =\dim\Gr(2, \tilde{P})=2\left( \binom{n+1}{2}-\left( \binom{s+1}{2}+s(n-2s) \right) -2 \right).
        \]

        The dimension of $Y$ is then
        \[
                \dim Y=\dim\flag{s,n-s}{\KK^n}+\dim\pi_1^{-1}([U,W])=n^2+n-s-4.
        \]
        We now see that
        $
            \dim \Gr(2,\sym_n)-\dim Y = s>0.
        $
        Since $\pi_2$ is not surjective and its image is closed, the claim is proved.
    \end{proof}

    \begin{lemma} \label{lem:middle-compression-space}
        Let $r=2s+1$ and assume $[Q]\in\fano{k}{\SD[r]{n}}$ is a nested $t$-compression space only for $t=s$.
        \begin{enumerate}[label=(\roman*)]
            \item There exists a unique flag $U\subseteq W\subseteq\KK^n$ with $\dim U=s+1+n-r$, $\dim W=n-s$, and $U\perp_Q W$. \label{lem-part:unique-flag}
            \item If a flag $U'\subseteq W'\subseteq \KK^n$ satisfies $\dim U'+\dim W'\geq 2n-r+1$ and $U'\perp_Q W'$, then $\dim U'+\dim W'= 2n-r+1$ and $U'\subseteq W'$ is the unique flag in part~\ref{lem-part:unique-flag}. \label{lem-part:inequality-forces-equality}
        \end{enumerate}
    \end{lemma}
    \begin{proof}
        Part~\ref{lem-part:unique-flag}: The existence of the flag $U\subseteq W$ follows from the assumption on $[Q]$.
        To show uniqueness, assume $U_1\subseteq W_1$ is another flag satisfying the conditions.
        We have $\dim U\cap U_1\geq \dim U+\dim U_1 - n = 1+n-r$.
        Let $p\geq 0$ be such that $\dim U\cap U_1=p+1+n-r$.
        It follows $\dim (U+U_1)=n-p$.
        Since $U\cap U_1\perp_Q U+U_1$, $[Q]$ is a nested $p$-compression space, forcing $p=s$ and therefore $U_1=U$.

        Next, we show $W_1=W$.
        If $s=0$, then $W=W_1=\KK^n$.
        If $s>0$ and $W\neq W_1$, then $\dim W+W_1> n-s$.
        Choose subspaces $U_2\subseteq U$ and $W_2\subseteq W+W_1$ with $U_2\subseteq W_2$, $\dim U_2=(s-1)+1+n-r$ and $\dim W_2=n-(s-1)$.
        Then $U_2\perp_Q W_2$ turns $[Q]$ into a nested $(s-1)$-compression space, a contradiction.

        Part~\ref{lem-part:inequality-forces-equality}: The inequality implies $\dim U'\geq 1+n-r$.
        Let $p\geq0$ be such that $\dim U'=p+1+n-r$.
        Then $\dim W'\geq n-p$, and $U'\perp_Q W'$ makes $[Q]$ a nested $p$-compression space forcing $p=s$.
        If $\dim W'>n-s$, then similar to the proof of part~\ref{lem-part:unique-flag}, $[Q]$ becomes a nested $(s-1)$-compression space.
        We conclude $\dim W'=n-s$.
        Now use part~\ref{lem-part:unique-flag}.
    \end{proof}

    \begin{lemma}\label{lem:direct-sum-of-lines}
        Let $0<n'< n$ and assume $[Q]\in\Gr(k+1,\sym_n)$ is given by an $n\times n$ matrix of the form
        $Q=\begin{pmatrix}
               P & 0\\
               0 & P'
        \end{pmatrix}$
         where $P$ and $P'$ are square blocks of size $n'$ and $n-n'$ respectively, with linear forms as entries.
        Let $U=\spn\{e_1,\ldots,e_{n'}\}$ and $W=\spn\{e_{n'+1},\ldots,e_n\}$ and define $\pi_U$ and $\pi_W$ to be projections of $\KK^n$ onto $U$ and $W$.
        For a flag $U'\subseteq W'\subseteq\KK^n$ satisfying $U'\perp_Q W'$,
        \begin{enumerate}[label=(\roman*)]
            \item $U'\cap U \perp_P \pi_U(W')$ and $\pi_W(U')\perp_{P'} W'\cap W$, \label{lem-part:orthogonal}
            \item if $\dim U'+\dim W'\geq n+\ell+1$ for some $0\leq \ell<n$, then $\dim U'\cap U+ \dim\pi_U(W')\geq n'+\ell+1$ or $\dim \pi_W(U')+ \dim W'\cap W\geq n-n'$, \label{lem-part:dims}
            \item if $\dim U'+\dim W'= n+\ell+1$ for some $0\leq \ell<n$, we have $\dim U'\cap U+ \dim\pi_U(W')= n'+\ell+1$ if and only if $\dim \pi_W(U')+ \dim W'\cap W= n-n'$. \label{lem-part:two-equality}
        \end{enumerate}
    \end{lemma}
    \begin{proof}
        Part~\ref{lem-part:orthogonal}: Let $u\in U'\cap U$ and $v\in\pi_U(W')$.
        Then $v=w'-w''$ for some $w'\in W'$ and $w''\in \pi_W(W')$.
        Since $U'\perp_Q W'$ and $U\perp_Q W$ we have $u\perp_Q w'$ and $u\perp_Q w''$, hence $u\perp_Q v$.
        This shows $U'\cap U \perp_Q \pi_U(W')$ and since $U'\cap U,\pi_U(W')\subseteq U$, they are also orthogonal with respect to $P$.
        The proof of $\pi_W(U')\perp_{P'} W'\cap W$ is similar.

        Part~\ref{lem-part:dims}: The exact sequence $0\to U'\cap U\to U' \xrightarrow{\pi_W} W$ gives $\dim U'=\dim U'\cap U+ \dim \pi_W(U')$ and similarly $\dim W'=\dim\pi_U(W')+ \dim W'\cap W$.
        Assume to the contrary that $\dim U'\cap U+ \dim\pi_U(W')< n'+\ell+1$ and $\dim \pi_W(U')+ \dim W'\cap W< n-n'$.
        Adding these inequalities yields
        \[
            \left(\dim U'\cap U+ \dim \pi_W(U')\right)+ \left(\dim\pi_U(W')+ \dim W'\cap W\right) < n+\ell+1.
        \]
        But the left-hand side is $\dim U'+\dim W'\geq n+\ell+1$, a contradiction.

        Part~\ref{lem-part:two-equality} follows from $\left( \dim U'\cap U +\dim \pi_U(W') \right)+\left(\dim \pi_W(U')+\dim W'\cap W\right)=\dim U' +\dim W'$.
    \end{proof}

    \begin{lemma}\label{lem:existence-only-nested-s-compression-space}
    Let $n$, $r$, and $s$ be integers with $r=2s+1$.
    Consider the point $[P]\in\fano{1}{\SD[r]{n}}$ defined by
    \[
        P=\begin{pmatrix}
              0   & B \\
              B^t & 0
        \end{pmatrix},\quad\text{with}\quad
        B=\left(
              \begin{array}{cccc|c}
                  z_0 & z_1    &        &     &                 \\
                  & \ddots & \ddots &     & \mathlarger{\mathlarger{0}}  \\
                  &        & z_0    & z_1 &
              \end{array}
        \right)_{s\times (s+1+n-r)},
    \]
    where the entries not shown in $B$ are all zero.
    Then $[P]$ is a nested $t$-compression space only for $t=s$.
    \end{lemma}
    \begin{proof}
        In~\cite[Proof of Proposition 4.5]{fano-nathan-chan}, the authors prove $[P]$ is not a $t$-compression space for $t<s$ and by Corollary~\ref{cor:s-comp}, it is also not a $t$-compression space for $t>s$ since $s$ is the middle integer in $0,\ldots,r-1$.
        Therefore, $[P]$ is a nested $t$-compression space only for $t=s$.
    \end{proof}

    We now prove Theorem~\ref{thm:fano-lines-components} as part of the next two results.

    \begin{theorem}
        \label{thm:fano-lines-components-complete}
        The following hold for the Fano scheme of lines $\fano{1}{\SD[r]{n}}$.
        \begin{enumerate}[label=(\roman*)]
            \item The subschemes $\compk{1}{s}$ for $s=0, 1,\ldots,\lfloor\frac{r-1}{2}\rfloor$ are distinct, and are exactly the irreducible components of $\fano{1}{\SD[r]{n}}$.
            Moreover, $\bigcap_s \compk{1}{s}\neq\emptyset$.
            \label{fano-lines-components-intersect}
            \item Let $\,\mathcal{U}(s)$ be the restriction of the universal bundle on $\fano{\kappa(s)}{\SD[r]{n}}$ to $\comp{s}$.
            The fiber of the natural map $\Gr\left( 2,\mathcal{U}(s) \right)\to \fano{1}{\SD[r]{n}}$ over a general point of $\compk{1}{s}$ is finite.
            \label{fano-lines-components-finite-fiber}
        \end{enumerate}
    \end{theorem}
    \begin{proof}
        Part~\ref{fano-lines-components-intersect}: By the Kronecker canonical form of a pencil of singular matrices~\cite[XII.4]{gantmacher}, every line on $\SD[r]{n}$ is a compression space, hence a nested one by Corollary~\ref{cor:all-s-comp-spaces}.
        Therefore, the subschemes $\compk{1}{s}$ cover $\fano{1}{\SD[r]{n}}$.
        For every $0\leq s\leq \lfloor\frac{r-1}{2}\rfloor$, we show $\compk{1}{s}\backslash\bigcup_{t\neq s}\compk{1}{t}\neq \emptyset$ which proves that $\compk{1}{s}$ is an irreducible component.

        Set $r'=2s+1$ and $n'=r'+n-r$.
        By Lemma~\ref{lem:existence-only-nested-s-compression-space}, there exists a pencil $[P]\in\fano{1}{\mathrm{SD}^{r'}_{n'}}$ with $P$ an $n'\times n'$ matrix over $\KK[z_0,z_1]_1$ that is a nested $t$-compression space only for $t=s$.
        Also, take $P'$ to be a general symmetric matrix of size $n-n'$ over $\KK[z_0,z_1]_1$ and define the $n\times n$ matrix
        \[
            Q=\begin{pmatrix}
                  P & 0  \\
                  0 & P'
            \end{pmatrix}.
        \]
        Define $U = \spn\{e_1,\ldots,e_{n'}\}$ and $W = \spn\{e_{n'+1},\ldots,e_n\}$, so that $\KK^n=U\oplus W$ and the matrices $P$ and $P'$ are pencils of bilinear forms on $U$ and $W$ respectively.
        We claim $[Q]\in\fano{1}{\SD[r]{n}}$ and is a nested $t$-compression space if and only if $t=s$.

        To show $Q$ is a nested $s$-compression space, let $U_1\subseteq W_1$ be a flag in $U$ with $U_1\perp_P W_1$ making $P$ into a nested $s$-compression space.
        By definition $\dim U_1=s+1+n'-r'$ and $\dim W_1=n'-s$.
        Then the flag $U_1\subseteq W_1+W$ in $\KK^n$ satisfies $U_1\perp_Q W_1+W$ and makes $Q$ into a nested $s$-compression space because $\dim U_1=s+1+n-r$ and $\dim W_1+W=n-s$.
        Here we used $n-r=n'-r'$.

        Conversely, assume $Q$ is a nested $t$-compression space and $U'\subseteq W'\subseteq \KK^n$ is a flag satisfying $\dim U'=t+1+n-r$, $\dim W'=n-t$ and $U'\perp_Q W'$.
        Our goal is to show $t=s$.
        By Lemma~\ref{lem:direct-sum-of-lines} part~\ref{lem-part:orthogonal}, we have $U'\cap U\perp_P \pi_U(W')$ and $\pi_W (U')\perp_{P'} W'\cap W$.
        Since $P'$ is not singular, $\dim \pi_W (U')+ \dim W'\cap W\leq \dim W=n-n'$.
        By Lemma~\ref{lem:direct-sum-of-lines} parts~\ref{lem-part:dims} and~\ref{lem-part:two-equality} with $\ell=n-r=n'-r'$, the orthogonal flag $U'\cap U\subseteq \pi_U(W')$ in $U$ satisfies $\dim U'\cap U+\dim\pi_{U}(W')\geq 2n'-r'+1$.
        Applying Lemma~\ref{lem:middle-compression-space} part~\ref{lem-part:inequality-forces-equality}, we have that this last inequality is an equality and that $U'\cap U\subseteq \pi_{U}(W')$ is the unique flag making $P$ a nested $s$-compression space.
        In particular $\dim U'\cap U=s+1+n'-r'$ and $\dim \pi_U(W')=n'-s$.
        Lemma~\ref{lem:direct-sum-of-lines} part~\ref{lem-part:two-equality} now gives $\dim \pi_W(U')+\dim W'\cap W=n-n'$.

        We compare some of the above dimensions and use $n-r=n'-r'$ and $r'=2s+1$.
        Comparing the dimensions of $U'\cap U$ and $U'$ gives $s\leq t$.
        Next, note that $U'\cap U$ is contained in $\pi_U(W')$ and they have the same dimension which means $\pi_U(W')=U'\cap U$.
        We conclude that $\pi_U(U')=U'\cap U$ implying $\pi_W(U')=U'\cap W$.
        Therefore $\pi_W(U')\subseteq W'\cap W$ is a flag in $W$ and by Lemma~\ref{lem:direct-sum-of-lines} part~\ref{lem-part:orthogonal}, $\pi_W(U')\perp_{P'} W'\cap W$.

        The exact sequences $0\to U'\cap U\to U' \xrightarrow{\pi_W} W$ and $0\to W'\cap W\to W' \xrightarrow{\pi_U} U$  give the dimensions of the flag $\pi_W(U')\subseteq W'\cap W$ as
        \begin{equation}\label{eq:difference-in-dims}
            \begin{split}
                \dim \pi_W(U')&=\dim U' - \dim U'\cap U = t-s,\\
                \dim W'\cap W &= \dim W' - \dim \pi_U(W') = (n-n')-(t-s).
            \end{split}
        \end{equation}
        Since $P'$ is general and $\pi_W(U')\perp_{P'} W'\cap W$, by Lemma~\ref{lem:generic-matrix} we must have $t-s=0$.

        The last statement in~\ref{fano-lines-components-intersect} holds because $\bigcap_s \compk{1}{s}$ contains the point
        \[
            \begin{pmatrix}
                z_0 & z_1 &                 & \\
                z_1 &     &                 & \\
                &     & \text{\Large 0} & \\
                &     &                 &
            \end{pmatrix}.
        \]

        Part~\ref{fano-lines-components-finite-fiber}: The map is proper with image $\compk{1}{s}$.
        For a proper map, the dimension of fibers is upper semi-continuous on the target, therefore it is enough to show the fiber over $[Q]$ is a single point.
        Since $\Gr(k+1,\mathcal{U}(s))$ is a bundle over $\comp{s}$ and by Corollary~\ref{cor:kappa-s-component}$, \comp{s}$ is isomorphic in its reduced structure to $\flag{s+1+n-r,n-s}{\KK^n}$, it is enough show a flag $[U',W']\in\flag{s+1+n-r,n-s}{\KK^n}$ satisfying $U'\perp_Q W'$ is uniquely determined by $[Q]$.

        Using~\eqref{eq:difference-in-dims} and $t-s=0$, we have that $U'=U'\cap U$ and we saw in the proof of part~\ref{fano-lines-components-intersect} that $U'\cap U$ is uniquely determined by $[P]$ (the uniqueness is due to Lemma~\ref{lem:middle-compression-space} and does not depend on the representative matrix $P$).
        We claim $W'=U'\oplus W$, which shows that $W'$ is also uniquely determined by $[P]$.
        By~\eqref{eq:difference-in-dims}, $\dim W'\cap W=n-n'=\dim W$, i.e.\ $W\subseteq W'$.
        We have subspaces $U',W\subseteq W'$ with $U'\cap W=0$ (since $U'\subseteq U$ and $U\cap W=0$) and
        \[
            \dim U'+\dim W=(s+1+n'-r')+ (n-n')=n-s=\dim W'.
        \]
        Therefore, $W'=U'\oplus W$.
    \end{proof}

    In §\ref{sec:smoothness} we will prove Proposition~\ref{prop:when-gen-non-reduced} which gives that all irreducible components of $\fano{1}{\SD[r]{n}}$  (except $\compk{1}{\frac{r-1}{2}}$ when $r$ is odd) are generically non-reduced.

        \begin{corollary}
            \label{cor:dim-compk}
            The closed subscheme $\compk{k}{s}$ of $\fano{k}{\SD[r]{n}}$ has dimension
            \[
                \dim \compk{k}{s} = (s+1+n-r)(r-2s-1)+s(n-s)+(\kappa(s)-k)(k+1).
            \]
            In particular for $k=1$, each irreducible component $\compk{1}{s}$ of $\fano{1}{\SD[r]{n}}$ has dimension $nr+(s-1)(n-r)-5$ and when $r=n$, the components are equidimensional and of the expected dimension.
        \end{corollary}
        \begin{proof}
            The second statement follows from the first.
            By Theorem~\ref{thm:fano-lines-components-complete} part~\ref{fano-lines-components-finite-fiber}, there exists a line $[L]\in\compk{1}{s}$ for which there are only finitely many flags $U\subseteq W$ in $\flag{s+1+n-r,n-s}{\KK^n}$ such that $U\perp_L W$.
            After acting by $\gl{n}$, we may assume the standard $s$-compression space $[Q]\in\compk{\kappa(s)}{s}$ contains this line.
            Take a $k$-plane $[P]\in\compk{k}{s}$ with $L\subseteq P\subseteq Q$.
            Then there are only finitely many flags $U\subseteq W$ with $U\perp_P W$.
            The map in~\eqref{eq:map-from-Grassmann-bundle} is proper with image $\compk{k}{s}$ and the fiber over $[P]$ is finite.
            By upper semi-continuity of fiber dimensions, the fiber over a general point of $\compk{k}{s}$ is finite.
            Therefore,
            \begin{align*}
                \dim \compk{k}{s} &= \dim \Gr\left(k+1,\mathcal{U}(s)\right)\\
                &=\dim \comp{s}+\dim \Gr\left(k+1,\kappa(s)+1\right)\\
                &=\dim \flag{s+1+n-r,n-s}{\KK^n}+(\kappa(s)-k)(k+1)\qquad\text{(by Corollary~\ref{cor:kappa-s-component})}\\
                &=(s+1+n-r)(r-2s-1)+s(n-s)+(\kappa(s)-k)(k+1).
            \end{align*}
        \end{proof}

    \begin{corollary}\label{cor:s-s'-on-distinct-comp}
    Let $0\leq s,s'\leq\lfloor\frac{r-1}{2}\rfloor$ with $s\neq s'$.
    Assume $k\leq\min\{\kappa(s),\kappa(s')\}$ (equivalently the subschemes $\compk{k}{s}$ and $\compk{k}{s'}$ of $\fano{k}{\SD[r]{n}}$ are non-empty).
    Then
    \begin{enumerate}[label=(\roman*)]
        \item $\compk{k}{s}$ and $\compk{k}{s'}$ are distinct and one does not contain the other,
        \label{part:s-s'-distinct}
        \item No irreducible components of $\fano{k}{\SD[r]{n}}$ contains both $\compk{k}{s}$ and $\compk{k}{s'}$.
        \label{part:s-s'-different-irred-comp}
    \end{enumerate}
    \end{corollary}
    \begin{proof}
        Let $Q$ be the standard $s$-compression space.
        By Theorem~\ref{thm:fano-lines-components-complete} part~\ref{fano-lines-components-intersect}, there is a line $L$ inside $Q$ such that $L$ is not a nested $s'$-compression space.
        Any $k$-plane inside $Q$ containing $L$ is also not a nested $s'$-compression space, i.e.\ $\compk{k}{s}$ is not contained in $\compk{k}{s'}$.
        This proves part~\ref{part:s-s'-distinct} since the argument is symmetric in $s,s'$.

        For part~\ref{part:s-s'-different-irred-comp}, assume an irreducible component $Z$ of $\fano{k}{\SD[r]{n}}$ contains both $\compk{k}{s},\compk{k}{s'}$.
        Let $\mathcal{U}$ be the restriction to $Z$ of the universal bundle on $\fano{k}{\SD[r]{n}}$.
        The image of the natural map from the Grassmann bundle $\Gr(2,\mathcal{U})$ to $\fano{1}{\SD[r]{n}}$ is irreducible and contains both $\compk{1}{s},\compk{1}{s'}$ contradicting Theorem~\ref{thm:fano-lines-components-complete} part~\ref{fano-lines-components-intersect}.
    \end{proof}

    Next, we prove Theorem~\ref{thm:irreducibility} characterizing the irreducibility of $\fano{k}{\SD[r]{n}}$.

    \begin{proof}[Proof of Theorem~\ref{thm:irreducibility}]
        If $\fano{k}{\SD[r]{n}}$ is irreducible, then there is only one $0\leq s\leq \left\lfloor\frac{r-1}{2}\right\rfloor$ such that $k\leq \kappa(s)$.
        This is because if $s\neq s'$ and $k\leq \min\{\kappa(s),\kappa(s')\}$, then $\compk{k}{s},\compk{k}{s'}$ would be non-empty and contained in the irreducible component $\fano{k}{\SD[r]{n}}$, which contradicts Corollary~\ref{cor:s-s'-on-distinct-comp} part~\ref{part:s-s'-different-irred-comp}.
        Using the convexity of $\kappa(s)$ in $s$, we see that $k\leq \kappa(s)$ for only one integer $s$ if and only if the conditions on $k$ hold.

        Conversely, assume $k$ satisfies the conditions.
        By Theorem~\ref{thm:pazzis}, every point of $\fano{k}{\SD[r]{n}}$ is a compression space, i.e.\ $\compk{k}{0}$ or $\compk{k}{s_{\max}}$ covers $\fano{k}{\SD[r]{n}}$, depending on $k$ satisfying the first or second condition.
        The Fano scheme is thus irreducible.
    \end{proof}

    \section{Smoothness}\label{sec:smoothness}

    \subsection{Dimensions of tangent spaces}\label{subsec:dimensions-of-tangent-spaces}

    In this section, we use deformation theory to compute the dimensions of Zariski tangent spaces at select points of $\fano{k}{\SD[r]{n}}$ and investigate the smoothness of the scheme.
    The general setup is as follows.

    Let $S=\KK[x_{1,1},x_{1,2},\ldots,x_{n,n}]$ be the homogeneous coordinate ring of the projective space $\PP^{\binom{n+1}{2}-1}$ of all symmetric matrices and let $I$ be the ideal of $\SD[r]{n}$, generated by the $r\times r$ minors of the generic symmetric matrix in~\eqref{eq:gen-sym-matrix}.
    Let $J$ be the homogeneous ideal of a $k$-plane $Q$ inside $\SD[r]{n}$ and choose an isomorphism $S/J=\KK[z_0,\ldots,z_k]$.
    Then we can write $Q$ as an $n\times n$ matrix whose entries are linear forms in $z_0,\ldots, z_k$.
    The tangent space to the Fano scheme $\fano{k}{\SD[r]{n}}$ at $[Q]$
    is isomorphic to the space of first-order deformations of $Q$ inside $\SD[r]{n}$ which is known to be isomorphic to $\Hom_{\mathcal{O}_{\SD[r]{n}}}(\mathcal{J},\mathcal{O}_Q)$, where $\mathcal{J}$ is the ideal sheaf of $Q$ in $\SD[r]{n}$ (see~\cite[Theorem 2.4]{hartshorne-deformation}).

    We recall the equivalence of the categories of quasi-coherent sheaves on $\SD[r]{n}$ and saturated $S/I$-modules (see~\cite[§15.7]{vakil-rising-sea}).
    The ring $S/I$ is a saturated $S/I$-module since it is integrally closed~\cite[Theorem in §1]{goto}.
    This implies that $J/I$ is a saturated $S/I$-module as it is the largest homogenous ideal defining the $k$-plane $Q$ in  $S/I$.
    Moreover $S/J$ is saturated as it is isomorphic to a polynomial ring.
    By the equivalence of the above two categories, we have the $\KK$-vector space isomorphism
    $
        \Hom_{\mathcal{O}_{\SD[r]{n}}}(\mathcal{J},\mathcal{O}_Q)\simeq \Hom_{S/I}(J/I,S/J)_0,
    $
    where the right-hand side is the space of degree-preserving maps of graded $S/I$-modules.

    We will compute dimensions of Zariski tangent spaces at a general point $[Q]$ of the subscheme $\compk{k}{s}$ of $\fano{k}{\SD[r]{n}}$ by computing $\dim \Hom_{S/I}(J/I,S/J)_0$.
    We say a minor of a matrix is $\ell$-anchored if it involves the first $\ell$ rows and columns.

    \begin{lemma}
        \label{lem:account-of-tangent-dimension}
        Fix integers $3\leq r\leq n$, $k\geq 1$, and $0\leq s\leq \lfloor \frac{r-1}{2} \rfloor$.
        Let $[Q]\in\fano{k}{\SD[r]{n}}$ be a nested $s$-compression space of the form
        \begin{equation}
            \label{eq:s-comp-space-block-form-r}
            \begin{pNiceArray}{ccccc}[first-row,first-col]
                &&\outmat{s}&&\outmat{r-2s-1}&\outmat{s+1+n-r}\\
                \outmat{s}&&B&&C&D \\
                & &  &&& \\
                \outmat{r-2s-1}&&C^t&&E&0\\
                & &&  && \\
                \outmat{s+1+n-r}&&D^t&&0&0\\
            \end{pNiceArray},
        \end{equation}
        where the blocks $B,\ldots,E$ have entries in $\KK[z_0,\ldots,z_k]_1$.
        If $r-2s-1>0$ and $\det E=0$, then
        \[
            \dim T_{[Q]}\fano{k}{\SD[r]{n}}=\left( \binom{n+1}{2}-k-1 \right)(k+1).
        \]
        On the other hand if $r-2s-1=0$ or $\det E\neq 0$, then
        \[
            \dim T_{[Q]}\fano{k}{\SD[r]{n}}=a_{\det} + (s+1+n-r)(r-2s-1)(k+1)+(\kappa(s)-k)(k+1),
        \]
        where $a_{\det}$ is the dimension of the $\KK$-vector space of all symmetric matrices $A$ of size $s+1+n-r$ with entries in $\KK[z_0,\ldots,z_k]_1$ such that the $s$-anchored $(2s+1)\times(2s+1)$ minors of the following matrix vanish
        \begin{equation}
            \label{eq:permissible-A}
            \begin{pNiceArray}{cccc}[first-row,first-col]
                &&\outmat{s}&\outmat{s+1+n-r}\\
                \outmat{s}&&0&D \\
                &&& \\
                \outmat{s+1+n-r}&&D^t&A\\
            \end{pNiceArray}.
        \end{equation}
    \end{lemma}
    \begin{proof}
        As explained above, $\dim T_{[Q]}\fano{k}{\SD[r]{n}}=\dim\Hom_{S/I}(J/I,S/J)_0$.
        A set of generators for $J/I$ consists of $x_{i,j}$ corresponding to the entries in the zero blocks in~\eqref{eq:s-comp-space-block-form-r} together with $\kappa(s)-k$ linearly independent linear forms $f_i$.
        An element $\phi\in\Hom_{S/I}(J/I,S/J)_0$ is determined by specifying the images $\phi(x_{i,j})$ and $\phi(f_i)$ in $\KK[z_0,\ldots,z_k]_1$.
        These images must be chosen in a way that every relation satisfied by the generators $x_{i,j}, f_i$ is also satisfied by their images.

        The Koszul relations $a_1(a_2)-a_2(a_1)$ for $a_i\in J$ are obviously mapped to zero.
        The only other relations come from writing the $r\times r$ minors of the generic symmetric matrix (generators of $I$) as $S/I$-linear combinations of the generators of $J$.
        A minor that is not $(r-s-1)$-anchored has at least $s+2$ of the last $s+1+n-r$ columns and is automatically sent to zero in $S/J$ since every term in the minor involves two $x_{i,j}$ corresponding to positions in the zero blocks, hence two generators of $J/I$.

        We inspect each term of an $(r-s-1)$-anchored minor up to Koszul relations.
        Such a minor involves $s+1$ of the last $s+1+n-r$ columns in~\eqref{eq:s-comp-space-block-form-r}.
        If a term involves a position in block $B$ or $C$, then it  involves fewer than $s$ positions in block $D$, hence involves two positions in the zero blocks and is mapped to zero.
        Similarly, if a term has a position in the $(r-2s-1)\times (s+1+n-r)$ zero block, then it will have fewer than $r-2s-1$ positions in block $E$, thus involves a position in the $(s+1+n-r) \times (r-2s-1)$ block, i.e.\ two positions in the zero blocks.
        To remove the terms that are mapped to zero from consideration, we may replace blocks $B, C, C^t$, and the $(r-2s-1)\times (s+1+n-r)$ and $(s+1+n-r) \times (r-2s-1)$ blocks with zeros in the generic symmetric matrix.

        Let $A$ be the matrix filled with $\phi(x_{i,j})$ for those $(i,j)$ in the lower right block of size $s+1+n-r$ in~\eqref{eq:s-comp-space-block-form-r}.
        The entries of $A$ must be such that the $(r-s-1)$-anchored $r\times r$ minors of
        \begin{equation}
            \label{eq:s-comp-blocks-zero-substituted}
            \begin{pNiceArray}{ccccc}[first-row,first-col]
                &&\outmat{s}&&\outmat{r-2s-1}&\outmat{s+1+n-r}\\
                \outmat{s}&&0&&0&D \\
                & &  &&& \\
                \outmat{r-2s-1}&&0&&E&0\\
                & &&  && \\
                \outmat{s+1+n-r}&&D^t&&0&A\\
            \end{pNiceArray},
        \end{equation}
        vanish, since these minors are the images in $S/J$ of the $r\times r$ minors of the generic symmetric matrix under the map $\phi$.
        Note that each anchored minor involves exactly one position in $A$.
        Therefore, the fact that $\phi$ sends generators of $I$ to zero imposes conditions only on $\phi(x_{i,j})$ in the $A$ block.
        This means there are no restrictions on the images of $x_{i,j}$ in the $(r-2s-1)\times (s+1+n-r)$ block and the images of $f_i$s.
        Let $a_{\det}$ be the dimension of the space of all symmetric matrices $A$ such that the $(r-s-1)$-anchored $r\times r$ minors of~\eqref{eq:s-comp-blocks-zero-substituted} vanish.

        An $(r-s-1)$-anchored $r\times r$ minor of the matrix in~\eqref{eq:s-comp-blocks-zero-substituted} is up to a sign the product of $\det E$ and an $s$-anchored $(2s+1)\times(2s+1)$ minor of~\eqref{eq:permissible-A}.
        If $\det E=0$, then there are no restrictions on what the entries of $A$ are.
        Therefore we can freely map generators of $J/I$ to elements in $\KK[z_0,\ldots,z_k]_1$ and we get
        \[
            \dim \Hom_{S/I}(J/I,S/J)_0=\left( \binom{n+1}{2}-k-1 \right)(k+1).
        \]

        If $r-2s-1=0$ or $\det E\neq 0$, then $a_{\det}$ is also the dimension of the space of all symmetric matrices such that the $s$-anchored $(2s+1)\times (2s+1)$ minors in~\eqref{eq:permissible-A} vanish.
        Since there are no restrictions on the images of $x_{i,j}$ in the $(r-2s-1)\times (s+1+n-r)$ block and $f_i$s we have
        \[
            \dim \left(\Hom_{S/I}(J/I,S/J)_0\right) = a_{\det} + (s+1+n-r)(r-2s-1)(k+1)+(\kappa(s)-k)(k+1).\qedhere
        \]
    \end{proof}

    To compute the dimension of the tangent space at a general point of the subscheme $\compk{k}{s}$ of $\fano{k}{\SD[r]{n}}$, we first consider the case $r=n$ with $n$ odd, and $s=\frac{r-1}{2}$.
    But before that, we need two preliminary lemmas.

    \begin{lemma}[Jensen's formula~{\cite[§3]{jensen}}]
        \label{lem:jensen}
        Let $L\geq 0$ be an integer and assume $\alpha,\beta,\gamma$ are complex numbers.
        Then
        \[
            \sum_{\ell=0}^L \binom{\alpha+\beta \ell}{\ell} \binom{\gamma-\beta \ell}{L-\ell} =\sum_{\ell=0}^L \binom{\alpha+\gamma-\ell}{L-\ell} \beta^{\ell}.
        \]
        In particular the sum $\sum_{\ell=0}^L \binom{\alpha-\ell}{\ell} \binom{\gamma+\ell}{L-\ell}$ depends only on $\alpha+\gamma$ and $L$.
    \end{lemma}

    \begin{lemma}
        \label{lem:expand-det}
        Let $n\geq3$ be odd and set $s=\frac{n-1}{2}$.
        Suppose $A=[a_{i,j}]$ and $D$ are $(s+1)\times(s+1)$ and $s\times (s+1)$ matrices respectively.
        Then
        \[
            \det\left(
            \begin{array}{c|ccc}
                0   & & D & \\
                \hline
                & &   & \\
                D^t & & A & \\
                & &   & \\
            \end{array}
            \right)_{n\times n}=-\sum_{1\leq i,j\leq s+1}(-1)^{i+j}a_{i,j}p_i p_j,
        \]
        where $p_i$ is the minor of $D$ after deleting its $i$-th column.
    \end{lemma}

    \begin{proof}
        Suppose $a_{i,j}$ are independent forms.
        Note that we are not assuming $A$ is symmetric.
        In the expansion of the above determinant, due to the upper left zero block, every term involves exactly one entry in $A$.
        To get the coefficient of $a_{i,j}$, we expand the determinant along row $i$ and obtain $-(-1)^{i+j}p_i p_j$ as the coefficient.
    \end{proof}
    \begin{notation}
        For an integer $s\geq 1$, define the matrix $D_s(z_0, z_1,z_2)$ as
        \begin{equation}\label{eq:block-D}
            D_s(z_0,z_1,z_2)=
            \begin{pmatrix}
                z_0 & z_1    & z_2    &        &     \\
                & \ddots & \ddots & \ddots &     \\
                &        & z_0    & z_1    & z_2 \\
                &        &        & z_0    & z_1
            \end{pmatrix}_{s\times (s+1)},
        \end{equation}
        where the entries not shown are all zero.
        When $s=1$, the variable $z_2$ is not present.
    \end{notation}
    \begin{proposition}
        \label{prop:tangent-dim-smooth}
        Let $r=n$ be odd and suppose $1\leq k\leq \kappa(s)$ where $s=\frac{r-1}{2}$.
        There exists a point $[Q]$ on the subscheme $\compk{k}{s}$ of $\fano{k}{\SD[r]{n}}$ such that
        \[
            \dim T_{[Q]}\fano{k}{\SD[r]{n}}=s(s+1)+\left(\kappa(s)-k\right)(k+1).
        \]
    \end{proposition}
    \begin{proof}
        Consider a nested $s$-compression space $[Q]$ given by the symmetric matrix
        \begin{equation}
            \label{eq:middle-s-compression}
            Q=\left(
            \begin{array}{c|ccc}
                B   & & D & \\
                \hline
                & &   & \\
                D^t & & 0 & \\
                & &   & \\
            \end{array}
            \right),
        \end{equation}
        where $B$ and $D$ are  $s\times s$ and $s\times (s+1)$ matrices respectively, with entries in $\KK[z_0,\ldots,z_k]_1$ such that $Q$, when viewed as a space of matrices, has $\KK$-dimension $k+1$.
        We claim that $[Q]$ is the desired point for suitable $D$ and arbitrary $B$.

        By Lemma~\ref{lem:account-of-tangent-dimension}, $\dim T_{[Q]}\fano{k}{\SD[r]{n}}=a_{\det}+\left(\kappa(s)-k\right)(k+1)$,
        where $a_{\det}$ is the $\KK$-dimension of the space of $(s+1)\times(s+1)$ matrices $A$ over $\KK[z_0\ldots,z_k]_1$ with
        \begin{equation}
            \label{eq:det-0}
            \det\left(
            \begin{array}{c|ccc}
                0   & & D & \\
                \hline
                & &   & \\
                D^t & & A & \\
                & &   & \\
            \end{array}
            \right)=0.
        \end{equation}

        It remains to exhibit for $1\leq k\leq \kappa(s)$ a matrix $D$ for which $a_{\det}=s(s+1)$.
        We will choose a matrix $D$ such that its rows are $\KK$-linearly independent and satisfy the following condition: The space $\mathcal{D}_{\mathrm{row}}$ of $(s+1)\times(s+1)$ matrices whose rows are $\KK$-linear combinations of the rows of $D$ has trivial intersection with the space $\mathcal{D}_{\mathrm{col}}$ of $(s+1)\times(s+1)$ matrices whose columns are $\KK$-linear combinations of the columns of $D^t$.
        Note that $\mathcal{D}_{\mathrm{col}}=\{\tilde{D}^t\,\big|\, \tilde{D}\in \mathcal{D}_{\mathrm{row}}\}$.

        The conditions on $D$ ensure that $\tilde{D}+\tilde{D}^t=0$ implies $\tilde{D}=0$ when $\tilde{D}\in \mathcal{D}_{\mathrm{row}}$, and the space of $(s+1)\times(s+1)$ matrices $\mathcal{D}=\{ \tilde{D}+\tilde{D}^t\,\big|\, \tilde{D}\in\mathcal{D}_{\mathrm{row}} \}$
        has linear dimension $s(s+1)$.

        The next step will be to show that $\mathcal{D}$ is the space of those matrices $A$ that satisfy~\eqref{eq:det-0}.
        Obviously, a matrix $A\in\mathcal{D}$ satisfies~\eqref{eq:det-0}, because by construction the block $A$ can be reduced to zero by applying \rowcol{} operations on the matrix in~\eqref{eq:det-0}.
        We just need to prove the converse, that if an $(s+1)\times(s+1)$ matrix $A$ satisfies~\eqref{eq:det-0}, then $A\in \mathcal{D}$.

        \emph{Case $k=1$:}
        Take $D=D_s(z_0,z_1,0)$  (see~\eqref{eq:block-D}).
        Let $A=[a_{i,j}]$ be a symmetric matrix satisfying~\eqref{eq:det-0} with $a_{i,j}\in\KK[z_0,z_1]_1$.
        We show $A\in \mathcal{D}$ by proving that we can zero out the entries of $A$ in~\eqref{eq:det-0} through suitable \rowcol{} operations using the blocks $D, D^t$.
        We may assume (after \rowcol{} operations) that no entries of $A$ except possibly $a_{s+1,s+1}$ involves $z_0$.
        But $a_{s+1,s+1}$ cannot involve $z_0$ because otherwise by Lemma~\ref{lem:expand-det}, the determinant in~\eqref{eq:det-0} would have a $z_0^{2s+1}$ term.
        Therefore, we may write $a_{i,j}=\alpha_{i,j}z_1$ for some $\alpha_{i,j}\in\KK$.

        Furthermore, by using successive \rowcol{} operations, we can push $z_1$-terms in the matrix $A$ to the boundary (i.e.\ the outermost entries) of this matrix without producing $z_0$ in the entries of $A$ (see Figure~\ref{fig:pushing-z1}).
        We may thus assume $a_{i,j}=0$ for $1<i,j<s+1$.
        By Lemma~\ref{lem:expand-det}, expanding the determinant in~\eqref{eq:det-0} gives
        $
            \sum_{i,j}(-1)^{i+j}a_{i,j}p_i p_j=0,
        $
        where $p_i=z_0^{i-1} z_1^{s+1-i}$, $1\leq i\leq s+1$.
        Hence
        \[
            \sum_{i,j}(-1)^{i+j}\alpha_{i,j} z_0^{i+j-2} z_1^{2s+3-(i+j)}=0.
        \]
        Inspecting the $(z_0,z_1)$-degree, we see that for a fixed $\ell$, the sum $\sum_{i+j=\ell}(-1)^{i+j}\alpha_{i,j}$ vanishes.
        Since $A$ has entries only at its boundary, we get $\alpha_{i,j}=0$ for all $i,j$.

        \begin{figure}
            \centering
            \begin{subfigure}{.48\textwidth}
                \centering
                $
                \left(
                \begin{array}{cc|ccc}
                    & & z_0 & z_1 & \\
                    & & & z_0 & z_1 \\
                    \hline
                    z_0 & & & & \\
                    z_1 & z_0 & & \alpha z_1 & \\
                    & z_1 & & &
                \end{array}
                \right)
                $
                \caption{Original matrix}
            \end{subfigure}
            \begin{subfigure}{.48\textwidth}
                \centering
                $
                \left(
                \begin{array}{cc|ccc}
                    & & z_0 & z_1 & \\
                    & & & z_0 & z_1 \\
                    \hline
                    z_0 & & &  &\frac{\alpha}{2}z_1 \\
                    z_1 & z_0 &  & 0 & \\
                    & z_1 & \frac{\alpha}{2}z_1& &
                \end{array}
                \right)
                $
                \caption{After two operations}
            \end{subfigure}
            \caption{Pushing $z_1$ to the boundary by row-column operations}\label{fig:pushing-z1}
        \end{figure}

        \emph{Case $k=2$:}
        Take $D=D_s(z_0,z_1,z_2)$.
        Let $A=[a_{i,j}]$ be a symmetric matrix satisfying~\eqref{eq:det-0} where $a_{i,j}\in\KK[z_0,z_1,z_2]_1$.
        Similar to the previous case we show that we can reduce the block $A$ to zero through \rowcol{} operations by using the blocks $D, D^t$.
        Using what we proved for $k=1$, we may assume that after suitable \rowcol{} operations on the matrix in~\eqref{eq:det-0}, the entries of $A$ are scalar multiples of $z_2$.
        We claim this reduces $A$ to zero.
        In Lemma~\ref{lem:lin-indep-minors}, we will prove that the polynomials $p_i p_j$ for $1\leq i\leq j\leq s+1$ are $\KK$-linearly independent where $p_i$ is the minor of $D$ after deleting the $i$-th column.
        The determinant~\eqref{eq:det-0} is
        $
            \sum_{i,j} (-1)^{i+j} a_{i,j} p_i p_j =0,
        $
        which by linear independence of $p_i p_j$ gives $a_{i,j}=0$.

        \emph{Case $k>2$:}
        Take $D=D_s(z_0,z_1,z_2)+D'$ where $D'$ is any $s\times (s+1)$ matrix with entries in $\KK[z_3,\ldots,z_k]_1$ such that $Q$, as a space of matrices, has $\KK$-dimension $k+1$.
        Assume $A$ has entries in $\KK[z_0,\ldots,z_k]_1$ and satisfies~\eqref{eq:det-0}.
        Using the case $k=2$, we can assume that after \rowcol{} operations, $A$ does not involve $z_0,z_1,z_2$.
        We prove this reduces $A$ to zero.
        To see that for a fixed $2<\ell\leq k$, the entries of $A$ do not involve $z_{\ell}$, set all the $z_i$ except $z_0,z_1,z_2,z_{\ell}$ equal to zero.
        The entries of $A$ are now scalar multiples of $z_{\ell}$.
        By the same argument as in the case $k=2$, all entries of $A$ are zero.
        Since $\ell$ was arbitrary, we are done.

    \end{proof}

    We used the following lemma in the proof of the previous proposition.

    \begin{lemma}
        \label{lem:lin-indep-minors}
        Let $s\geq 1$ and consider $D_s(z_0,z_1,z_2)$ in~\eqref{eq:block-D}.
        For $1\leq i \leq s+1$, let $p_i$ be the minor of $D_s(z_0,z_1,z_2)$ after deleting the $i$-th column.
        Then
        \begin{equation}
            \label{eq:expr-p_i}
            p_i=\sum_{\ell=0}^{\frac{s-i+1}{2}}(-1)^{\ell} \binom{s-i-\ell+1}{\ell} z_0^{i+\ell-1} z_1^{s-i+1-2\ell} z_2^{\ell}.
        \end{equation}
        Furthermore,
        \begin{enumerate}
        [label=(\roman*)]
            \item for $1\leq i,j,i',j'\leq s+1$ with $i+j\neq i'+j'$, the monomials appearing in $p_i p_j$ are distinct from those in $p_{i'} p_{j'}$ ,
            \label{expr-p_i-part1}
            \item for a fixed integer $2\leq d \leq 2(s+1)$, the polynomials $p_i p_j$ for $1\leq i\leq j\leq s+1$ and $i+j=d$ are $\KK$-linearly independent.
            \label{expr-p_i-part2}
        \end{enumerate}
        In particular, the polynomials $p_i p_j$ for $1\leq i\leq j\leq s+1$ are $\KK$-linearly independent.
    \end{lemma}
    \begin{proof}
        Let us denote by $p_{i,s}$ the determinant of the matrix obtained by deleting column $i$ in $D_s(z_0,z_1,z_2)$.
        We have
        \begin{equation}
            \label{eq:initial-conditions}
            p_{s+1,s}=z_0^s,\quad p_{s,s}=z_0^{s-1} z_1,\quad p_{s-1,s}=z_0^{s-2}z_1^2-z_0^{s-1}z_2.
        \end{equation}
        Moreover, by expanding the determinant along the last row we obtain the recursion
        \begin{equation}
            \label{eq:recursion}
            p_{i,s}=z_1 p_{i,s-1}-z_0 z_2 p_{i,s-2}\qquad(1\leq i\leq s-2).
        \end{equation}
        The equations in~\eqref{eq:initial-conditions} and~\eqref{eq:recursion} uniquely determine $p_{i,s}$ for all $s$ and $1\leq i\leq s+1$.
        To prove the expression for $p_i$ in~\eqref{eq:expr-p_i}, it is enough to show that this expression satisfies both~\eqref{eq:initial-conditions} and~\eqref{eq:recursion}.
        Let $q_{i,s}$ be the right-hand side of~\eqref{eq:expr-p_i}.
        It is easy to verify $q_{i,s}$ satisfies~\eqref{eq:initial-conditions}.
        To show it satisfies the recursion, we compute
        \begin{equation*}
            \begin{split}
                z_1 & q_{i,s-1}-z_0 z_2 q_{i,s-2}\\
                & = z_1\sum_{\ell=0}^{\frac{s-i}{2}}(-1)^{\ell}\binom{s-i-\ell}{\ell}z_0^{i+\ell-1} z_1^{s-i-2\ell} z_2^{\ell} - z_0 z_2\sum_{\ell=0}^{\frac{s-i-1}{2}}(-1)^{\ell}\binom{s-i-1-\ell}{\ell} z_0^{i+\ell-1} z_1^{s-i-1-2\ell} z_2^{\ell}\\
                & = \sum_{\ell=0}^{\frac{s-i}{2}}(-1)^{\ell}\binom{s-i-\ell}{\ell}z_0^{i+\ell-1} z_1^{s-i+1-2\ell} z_2^{\ell} + \sum_{\ell=1}^{\frac{s-i+1}{2}} (-1)^{\ell}\binom{s-i-\ell}{\ell-1} z_0^{i+\ell-1} z_1^{s-i+1-2\ell} z_2^{\ell}\\
                & = \sum_{\ell=0}^{\frac{s-i+1}{2}}(-1)^{\ell}\binom{s-i-\ell+1}{\ell} z_0^{i+\ell-1} z_1^{s-i+1-2\ell} z_2^{\ell}\\
                & = q_{i,s}.
            \end{split}
        \end{equation*}

        Proof of part~\ref{expr-p_i-part1}: Using the expression for $p_i$, we have
        \[
            p_i p_j =\sum_{\ell=0}^{\frac{s-i+1}{2}} \sum_{\ell'=0}^{\frac{s-j+1}{2}} (-1)^{\ell+\ell'} \binom{s-i-\ell+1}{\ell}\binom{s-j-\ell'+1}{\ell'}z_0^{i+j+(\ell+\ell')-2} z_1^{2s-(i+j)+2-2(\ell+\ell')} z_2^{\ell+\ell'},
        \]
        and a similar expression for $p_{i'}p_{j'}$.
        If $p_i p_j$ and $p_{i'}p_{j'}$ share a monomial, then by inspecting the $(z_0,z_2)$-degree, we see that $i+j=i'+j'$.

        Proof of part~\ref{expr-p_i-part2}: Let $\tilde{p}_{i,j}=(p_i p_j)|_{z_0=z_1=1}$.
        It is enough to show that the polynomials $\tilde{p}_{i,j}$ for $1\leq i\leq j\leq s+1$, $i+j=d$ are $\KK$-linearly independent.
        Sort these polynomials in ascending order for the index $i$ as $\tilde{p}_{i_0,j_0},\tilde{p}_{i_0+1,j_0-1},\ldots,\tilde{p}_{i_1-1,j_1+1},\tilde{p}_{i_1,j_1}$, where $i_0$, respectively $i_1$, is the smallest, respectively largest, integer $i$ for which there exists an integer $j=j_0$, respectively $j=j_1$, satisfying $1\leq i\leq j\leq s+1$, $i+j=d$.
        Note that we may instead prove that the polynomials of successive differences
        $\tilde{p}_{i_0,j_0},\tilde{p}_{i_0,j_0}-\tilde{p}_{i_0+1,j_0-1},\ldots,\tilde{p}_{i_1-1,j_1+1}-\tilde{p}_{i_1,j_1}$
        are $\KK$-linearly independent.
        We prove this by showing that if $i\leq j-2$, then the term with the smallest degree in $\tilde{p}_{i,j}-\tilde{p}_{i+1,j-1}$ has degree $s-j+2$ (note that $\tilde{p}_{i_0,j_0}$ has a degree zero term).
        This would complete the proof of part~\ref{expr-p_i-part2}.

        To show the smallest degree in $\tilde{p}_{i,j}-\tilde{p}_{i+1,j-1}$ is $s-j+2$, in the expression
        \[
            \tilde{p}_{i,j}=\sum_{\ell=0}^{\frac{s-i+1}{2}} \sum_{\ell'=0}^{\frac{s-j+1}{2}} (-1)^{\ell+\ell'} \binom{s-i-\ell+1}{\ell}\binom{s-j-\ell'+1}{\ell'} z_2^{\ell+\ell'},
        \]
        we make the change of variable $L=\ell+\ell'$ to get
        \begin{alignat}{2}
            \tilde{p}_{i,j} &= && \sum_{L=0}^{s-\frac{i+j}{2}+1}\sum_{\substack{0\leq\ell\leq L \\ 0\leq L-\ell\leq s-j+1}}(-1)^L \binom{s-i-\ell+1}{\ell}\binom{s-j-(L-\ell)+1}{L-\ell}z_2^L \nonumber\\
            &= && \sum_{L=0}^{s-j+2}\sum_{\substack{0\leq\ell\leq L \\ 0\leq L-\ell\leq s-j+1}}(-1)^L \binom{s-i-\ell+1}{\ell}\binom{s-j-(L-\ell)+1}{L-\ell}z_2^L \label{eq:sum-p_ij}\\
            & &&+\,\text{(terms with degree higher than $s-j+2$)}\nonumber\\
            & = && \sum_{L=0}^{s-j+1}\sum_{\ell=0}^{L}(-1)^L \binom{s-i-\ell+1}{\ell}\binom{s-j-(L-\ell)+1}{L-\ell}z_2^L\nonumber\\
            & &&+\sum_{\ell=1}^{s-j+2}(-1)^{s-j+2}\binom{s-i-\ell+1}{\ell}\binom{\ell-1}{s-j+2-\ell} z_2^{s-j+2}\nonumber\\
            & &&+\,\text{(terms with degree higher than $s-j+2$)}\nonumber\\
            & = && \sum_{L=0}^{s-j+2}\sum_{\ell=0}^{L}(-1)^L \binom{s-i-\ell+1}{\ell}\binom{s-j-(L-\ell)+1}{L-\ell}z_2^L\label{eq:ptilde1}\\
            & && -(-1)^{s-j+2}\binom{s-i+1}{0}\binom{-1}{s-j+2} z_2^{s-j+2}\nonumber\\
            & && +\,\text{(terms with degree higher than $s-j+2$)}\nonumber.
        \end{alignat}

        Similar to~\eqref{eq:sum-p_ij},
        \begin{equation}
            \label{eq:ptilde2}
            \begin{split}
                \tilde{p}_{i+1,j-1} = & \sum_{L=0}^{s-j+2}\sum_{\ell=0}^{L}(-1)^L \binom{s-i-\ell}{\ell}\binom{s-j-(L-\ell)+2}{L-\ell}z_2^L\\
                &+\text{(terms with degree higher than $s-j+2$)}.
            \end{split}
        \end{equation}
        By Jensen's formula (Lemma~\ref{lem:jensen}), the double summations in~\eqref{eq:ptilde1} and~\eqref{eq:ptilde2} are equal, hence
        \[
            \tilde{p}_{i,j}-\tilde{p}_{i+1,j-1}=-z_2^{s-j+2} +\, \text{(terms with higher degrees)}.
        \]
        This completes the proof of part~\ref{expr-p_i-part2}.
        The last statement of the lemma follows from~\ref{expr-p_i-part1} and~\ref{expr-p_i-part2}.
    \end{proof}

    We need the next lemma in computing the dimension of tangent spaces to $\fano{k}{\SD[r]{n}}$ at a general point of $\compk{k}{s}$.

    \begin{corollary}
        \label{lem:generic-block-row-span}
        Let $s\geq 1$ and assume $D$ is a general $s\times (s+1)$ matrix with entries in $\KK[z_0,\ldots,z_k]_1$.
        Then a symmetric matrix $A$ of size $s+1$ over $\KK[z_0,\ldots,z_k]_1$ satisfies $A=\tilde{D}+\tilde{D}^t$ with $\tilde{D}$ an $(s+1)\times (s+1)$ matrix whose rows are in the row span of $D$ if and only if
        \begin{equation}
            \label{eq:general-D}
            \det
            \begin{pmatrix}
                0   & D \\
                D^t & A
            \end{pmatrix}
            =0.
        \end{equation}
    \end{corollary}
    \begin{proof}
        If $A=\tilde{D}+\tilde{D}^t$, use \rowcol{} operations to reduce $A$ to zero, then use Lemma~\ref{lem:expand-det}.
        Conversely, assume~\eqref{eq:general-D} holds.
        Let $n=r=2s+1$ and consider the Fano scheme $\fano{k}{\SD[r]{n}}$.
        We may assume $\compk{k}{s}\neq \emptyset$ or equivalently $k\leq \kappa(s)$:
        By a straightforward argument using Lemma~\ref{lem:expand-det} and the fact that in the lemma, the polynomials $p_i p_j$ are $\KK$-linearly independent for $1\leq i\leq j\leq s+1$ and general $D$, one can show that $A$ can involve only those variables present in $D$.
        Since $\kappa(s)+1$ is larger than the number of entries of $D$, we may assume $k\leq \kappa(s)$.

        Up to $\gl{n}$ action, a general point $[Q]$ on the subscheme $\compk{k}{s}$ of $\fano{k}{\SD[r]{n}}$ is of the form~\eqref{eq:middle-s-compression} where $B$ and $D$ are matrices filled with general linear forms.

        By Corollary~\ref{cor:dim-compk},  $\dim \compk{k}{s}=s(s+1)+\left(\kappa(s)-k\right)(k+1)$ and by Proposition~\ref{prop:tangent-dim-smooth}, the Fano scheme contains a point in $\compk{k}{s}$ with tangent space dimension equal to $\dim\compk{k}{s}$.
        This means the dimension of the tangent space to $\fano{k}{\SD[n]{n}}$ at $[Q]$ equals $\dim\compk{k}{s}$.
        This dimension together with Lemma~\ref{lem:account-of-tangent-dimension} give that for a general $D$, the dimension of all matrices $A$ satisfying~\eqref{eq:general-D} is $a_{\det}=s(s+1)=s(n-s)$.
        Since $D$ is general, the space of all matrices $A$ satisfying~\eqref{eq:general-D} contains the subspace of those matrices of the form $\tilde{D}+\tilde{D}^t$ where $\tilde{D}$ is an $(s+1)\times (s+1)$ matrix in the row span of $D$ and this subspace has dimension $s(n-s)$ too.
        Therefore every matrix $A$ satisfying~\eqref{eq:general-D} is of the required form.
    \end{proof}

    Next, we compute the dimension of the tangent space at a general point of $\compk{k}{s}$.

    \begin{proposition}
        \label{prop:dim-tangent-space-generic-point-r}
        Let $[Q]$ be a general point of the subscheme $\compk{k}{s}$ of $\fano{k}{\SD[r]{n}}$.
        Then
        \[
            \dim T_{[Q]}\left( \fano{k}{\SD[r]{n}} \right)=s(s+1+n-r)+(s+1+n-r)(r-2s-1)(k+1)+(\kappa(s)-k)(k+1).
        \]
    \end{proposition}
    \begin{proof}
        Up to $\gl{n}$ action, $Q$ is a matrix of the form~\eqref{eq:s-comp-space-block-form-r} where the entries in blocks $B,\ldots,E$ are general linear forms in $z_0,\ldots,z_k$.
        Since $E$ is general, $\det E\neq 0$ and by Lemma~\ref{lem:account-of-tangent-dimension},
        \[
            \dim T_{[Q]}\left( \fano{k}{\SD[r]{n}} \right)=a_{\det}+(s+1+n-r)(r-2s-1)(k+1)+(\kappa(s)-k)(k+1).
        \]
        Our task is to show $a_{\det}=s(s+1+n-r)$.
        Call a matrix $A$ of size $s+1+n-r$ over $\KK[z_0,\ldots,z_k]_1$ permissible if all $s$-anchored $(2s+1)\times(2s+1)$ minors of~\eqref{eq:permissible-A} vanish.
        The space of all permissible $A$ contains a subspace of dimension $s(s+1+n-r)$, namely the set of all matrices $\tilde{D}+\tilde{D}^t$ where $\tilde{D}$ is a square matrix of size $s+1+n-r$ whose rows are in the row span of the general matrix $D$.

        For a subset $T\subseteq\{1,\ldots s+1+n-r\}$, denote by $A_T$ the principal submatrix of $A$ given by rows and columns in $T$ (e.g.\ $A_{\{1,2\}}$ is the upper left $2\times 2$ submatrix of A).
        Also, denote by $D_T$ the submatrix of $D$ given by the columns in $T$.
        Consider the principal $s$-anchored $(2s+1)\times(2s+1)$ minors of~\eqref{eq:permissible-A}.
        By assumption each such minor vanishes and is of the form
        \[
            \det\begin{pNiceArray}{cc}[first-row]
                    \outmat{s} & \outmat{s+1}\\
                    0 & D_T \\
                    D_T^t & A_T\\
            \end{pNiceArray},
        \]
        for some $T$ of size $s+1$.
        By Lemma~\ref{lem:generic-block-row-span}, the vanishing of the above determinant implies $A_T=\tilde{D}_T+\tilde{D}_T^t$ where $\tilde{D}_T$ is a matrix whose rows are in the row span of $D_T$.
        We claim since this holds for every $T$ of size $s+1$, we have $A=\tilde{D}+\tilde{D}^t$ where $\tilde{D}$ is a matrix with rows in the row span of $D$.

        The fact that for a $T$ of size $s+1$, we have $A_T=\tilde{D}_T+\tilde{D}_T^t$, means that by \rowcol{} operations on
        \[
            \begin{pNiceArray}{ccc}[first-row]
                    &\outmat{s} & \outmat{s+1+n-r}\\
                    &0 & D \\
                    &D^t & 0\\
            \end{pNiceArray},
        \]
        we can produce the submatrix $A_T$ as a submatrix (given by rows and columns in $T$) of the lower right zero block.
        We show that for two subsets $T,T'\subseteq\{1,\ldots,s+1+n-r\}$ of size $s+1$, the \rowcol{} operations used to build $A_{T}$ is \textit{compatible} with the \rowcol{} operations used to build $A_{T'}$.
        Being compatible means if $i\in T\cap T'$, then the linear combination of the rows of $D$ (and columns of $D^t$) used to build the row of $A_{T}$ corresponding to row $i$ in $A$ is the same as the linear combination used to build that row in $A_{T'}$ (see Figure~\ref{fig:two-submatrices}).

        \begin{figure}
            \centering
            $
            \begin{pNiceMatrix}[first-col,left-margin=0.6em,right-margin=0.6em]
                &       \Large{0}     &   & D  &  \\
                &                     & \Block[draw,rounded-corners]{2-2}{} * & *  & * \\
                \text{row $i$ of $A\quad$} &           D^t       & * & \Block[draw,rounded-corners]{2-2}{} *  & * \\
                &                     & * & *  & *
                \CodeAfter  
                \tikz \draw  (1-|2) ++(-3pt,0) -- ++(0,-2cm);
                \tikz \draw  (2-|1) ++(0,2pt) -- ++(2.7cm,0);
            \end{pNiceMatrix}
            $
            \caption{Submatrices $A_T$ and $A_{T'}$ in the proof of Proposition~\ref{prop:dim-tangent-space-generic-point-r} when $s=1$, $n-r=1$, $T=\{1,2\}$, $T'=\{2,3\}$, and $i=2\in T\cap T'$}\label{fig:two-submatrices}
        \end{figure}

        To show compatibility for $T$ and $T'$, each of size $s+1$, assume first that they differ in one element.
        Set $T_0=T\cap T'$.
        Then $A_{T}$ and $A_{T'}$ share the $s\times s$ principal submatrix $A_{T_0}$ of $A$.
        A linear algebra argument based on dimensions of vector spaces shows that $A_{T_0}$ can be built in a unique way by \rowcol{} operations using rows of $D_{T_0}$ (and columns of $D_{T_0}^t$), proving that $T$ and $T'$ are compatible. 
        If $T$ and $T'$ differ in more than one element, there is a sequence $T=T_1, \ldots, T_{\ell}=T'$ of subsets of $\{1,\ldots,s+1+n-r\}$,
        each of size $s+1$, such that $T_{j}$ and $T_{j+1}$ differ only in one element for $1\leq j<\ell$.
        Since $T_{j}$ and $T_{j+1}$ are compatible for all $j$, we deduce $T$ and $T'$ are compatible.
        This completes the proof.
    \end{proof}

    \subsection{Smoothness}\label{subsec:smoothness}
    In this section we use the dimensions of tangent spaces to detect components and to characterize the smoothness of $\fano{k}{\SD[r]{n}}$.

    \begin{proposition}
        \label{prop:smooth-points}
        Let $r$ be odd, $s=\frac{r-1}{2}$ and $k\leq \kappa(s)$.
        The subscheme $\compk{k}{s}$ of $\fano{k}{\SD[r]{n}}$ is an irreducible component that contains a smooth point of $\fano{k}{\SD[r]{n}}$.
        Moreover, when $k=\kappa(s)$, every point of $\compk{k}{s}$ is a smooth point of $\fano{k}{\SD[r]{n}}$.
    \end{proposition}
    \begin{proof}
        By Proposition~\ref{prop:dim-tangent-space-generic-point-r} and Corollary~\ref{cor:dim-compk}, when $r-2s-1=0$,
        a general point of $\compk{k}{s}$ has tangent space dimension equal to $\dim\compk{k}{s}$.
        Therefore $\compk{k}{s}$ is a component containing a smooth point of $\fano{k}{\SD[r]{n}}$.
        When $k=\kappa(s)$, by Corollary~\ref{cor:kappa-s-component}, $\compk{\kappa(s)}{s}$ is a single $\gl{n}$ orbit, proving the last statement.
    \end{proof}

    Next, we prove Proposition~\ref{prop:when-gen-non-reduced} which states that our Fano schemes can have generically non-reduced components.
    \begin{proof}[Proof of Proposition~\ref{prop:when-gen-non-reduced}]
        The irreducible component $\compk{k}{s}$ is generically reduced if and only if $\dim T_{[Q]}\left( \fano{k}{\SD[r]{n}} \right)=\dim \compk{k}{s}$ for a general point $[Q]$ of $\compk{k}{s}$.
        By Proposition~\ref{prop:dim-tangent-space-generic-point-r} and Corollary~\ref{cor:dim-compk}
        \[
            \dim T_{[Q]}\left( \fano{k}{\SD[r]{n}} \right)-\dim \compk{k}{s} = (r-2s-1)\left((s+1+n-r)k-s\right).
        \]
        We now have the desired result since $\left((s+1+n-r)k-s\right)>0$.
    \end{proof}

    Finally, we prove the main result of this section, namely Theorem~\ref{thm:characterize-smoothness}, that characterizes when $\fano{k}{\SD[r]{n}}$ is smooth.

    \begin{proof}[Proof of Theorem~\ref{thm:characterize-smoothness}]
        Assume $r$ is odd and $k$ satisfies~\eqref{eq:smoothness-criterion}.
        By Theorem~\ref{thm:pazzis}, $\compk{k}{\frac{r-1}{2}}$ equals $\fano{k}{\SD[r]{n}}$ in its reduced structure.
        In Lemma~\ref{lem:smoothness-of-middle-s}, we will prove that a Borel fixed point of $\compk{k}{\frac{r-1}{2}}$ that is a nested $s$-compression space only for $s=\frac{r-1}{2}$ is a smooth point of $\fano{k}{\SD[r]{n}}$.
        Here every Borel fixed of the scheme is a nested $s$-compression space only for $s=\frac{r-1}{2}$ and hence is smooth.
        Because the non-smooth locus of $\fano{k}{\SD[r]{n}}$ is Borel-stable, if this locus were non-empty it would contain a Borel fixed point.
        We conclude the non-smooth locus is empty.

        To prove the converse, we note that due to the convexity of $\kappa(s)$ in $s$, we always have $k\leq \max\{\kappa(0),\kappa\left(\lfloor\frac{r-1}{2}\rfloor\right)\}$.
        This means we need to show the scheme is not smooth in the following four cases: $k\leq \kappa(0)$, $r$ is odd and $k\leq \kappa\left(\frac{r-3}{2}\right)$, $r$ is even and $\max\{\kappa(0),\kappa\left(\frac{r-4}{2}\right)\}<k\leq \kappa\left(\frac{r-2}{2}\right)$, $r$ is even and $k\leq \kappa\left(\frac{r-4}{2}\right)$.
        Table~\ref{tab:summary-of-cases} summarizes why $\fano{k}{\SD[n]{r}}$ is not smooth in each case.
        Recall that by Corollary~\ref{cor:s-s'-on-distinct-comp}, for distinct $s,s'$, the subschemes $\compk{k}{s}$, $\compk{k}{s'}$  lie on distinct components of $\fano{k}{\SD[r]{n}}$.

        \begin{table}
            \scalebox{0.9}{
        \begin{tabular}{p{.37\textwidth} p{.63\textwidth}}
            \toprule
            \textbf{Case} & \textbf{Reason why ${\fano{k}{\SD[r]{n}}}$ is not smooth}\\
            \midrule
            $\kappa(0)-(r-2)<k\leq \kappa(0)$,\newline $r$ even or odd  & ${\fano{k}{\SD[r]{n}}}$ is generically non-reduced on component $\compk{k}{0}$ (Corollary~\ref{cor:component-s=0} and Proposition~\ref{prop:when-gen-non-reduced}).\\
            \midrule
            $k\leq \kappa(0)-(r-2)$,\newline $r$ even or odd & $\compk{k}{0}\cap\compk{k}{1}\neq\emptyset$\newline because $\kappa(0)-(r-2)=g(\{0,1\})$ (Remark~\ref{rem:interpret-labels}).\\
            \midrule
            $k=\kappa(\frac{r-3}{2})$,\newline $r$ odd  & ${\fano{k}{\SD[r]{n}}}$ is generically non-reduced on the component $\compk{k}{\frac{r-3}{2}}$ (Corollary~\ref{cor:component-s>=1} and Proposition~\ref{prop:when-gen-non-reduced}).\\
            \midrule
            $k\leq \kappa(\frac{r-3}{2})-1$,\newline $r$ odd  & $\compk{k}{\frac{r-1}{2}}\cap\compk{k}{\frac{r-3}{2}}\neq 0$ \newline because $\kappa(\frac{r-3}{2})-1=g(\{\frac{r-3}{2},\frac{r-1}{2}\})$ (Remark~\ref{rem:interpret-labels}).\\
            \midrule
            $\Scale[0.9]{\max\{\kappa(0),\kappa(\frac{r-4}{2})\}<k\leq \kappa(\frac{r-2}{2})}$, \newline $r$ even & ${\fano{k}{\SD[r]{n}}}$ in its reduced structure equals $\compk{k}{\frac{r-2}{2}}$ (Theorem~\ref{thm:pazzis}) and is generically non-reduced (Proposition~\ref{prop:when-gen-non-reduced}).\\
            \midrule
            $k=\kappa(\frac{r-4}{2})$,\newline $r\geq 6$ even & ${\fano{k}{\SD[r]{n}}}$ is generically non-reduced on the component $\compk{k}{\frac{r-4}{2}}$ (Corollary~\ref{cor:component-s>=1} and Proposition~\ref{prop:when-gen-non-reduced}).\\
            \midrule
            $k=\kappa(\frac{r-4}{2})-1$,\newline $r=n=6$ & $\compk{k}{2}\cap\compk{k}{1}\neq\emptyset$ \newline because $\kappa(1)-1=g(\{1,2\})$ (Remark~\ref{rem:interpret-labels}).\\
            \midrule
            $k=\kappa(\frac{r-4}{2})-1$, \newline $r\geq 6$ even, $(n,r)\neq (6,6)$ & ${\fano{k}{\SD[r]{n}}}$ is generically non-reduced on the component $\compk{k}{\frac{r-4}{2}}$  (see Corollary~\ref{cor:component-s>=1} and Proposition~\ref{prop:s=1-special-case-component} for being a component and Proposition~\ref{prop:when-gen-non-reduced} for non-reducedness).\\
            \midrule
            $k\leq \kappa(\frac{r-4}{2})-2$,\newline $r\geq 6$ even & $\compk{k}{\frac{r-2}{2}}\cap\compk{k}{\frac{r-4}{2}}\neq \emptyset$ \newline because $\kappa(\frac{r-4}{2})-2=g(\{\frac{r-4}{2},\frac{r-2}{2}\})$ (Remark~\ref{rem:interpret-labels}).\\
            \bottomrule
        \end{tabular}}
        \caption{Non-smoothness of $\fano{k}{\SD[r]{n}}$ in Theorem~\ref{thm:characterize-smoothness}}
        \label{tab:summary-of-cases}
    \end{table}

    \end{proof}

    We used the following lemma in the proof of Theorem~\ref{thm:characterize-smoothness}.

    \begin{lemma}\label{lem:smoothness-of-middle-s}
    Let $r$ be odd and assume $[Q]$ is a Borel fixed point of $\compk{k}{\frac{r-1}{2}}$ such that $Q$ is a nested $s$-compression space only for $s=\frac{r-1}{2}$.
    Then $[Q]$ is a smooth point of $\fano{k}{\SD[r]{n}}$.
    \end{lemma}
    \begin{proof}
        The point $[Q]$ can be written in the form
        \begin{equation}
            \label{eq:s-mid-comp-space-block-form-r}
            \begin{pNiceArray}{ccccc}[first-row,first-col]
                &&\outmat{s}&&\outmat{s+1}&\outmat{n-2s-1}\\
                \outmat{s}&&B&&D_1&D_2 \\
                & &  &&& \\
                \outmat{s+1}&&D_1^t&&0&0\\
                & &&  && \\
                \outmat{n-2s-1}&&D_2^t&&0&0\\
            \end{pNiceArray},
        \end{equation}
        where the blocks $B,D_1,$ and $D_2$ have entries in $\KK[z_0,\ldots,z_k]_1$ and each entry is either an independent variable or zero.
        Since $Q$ is only a nested $(\frac{r-1}{2})$-compression space, the entries in $B$ and those on and above the anti-diagonal of $D_1$ are distinct independent variables (see Figure~\ref{fig:D1}).
        This is because $Q$ is a Borel fixed point, and by Proposition~\ref{prop:characterize-borel-fixed-points-grassmannian} if an entry in $B$ or on or above the anti-diagonal of $D_1$ is zero, then every entry to the right and bottom will be zero, and we will have a nested $s$-compression space for some $s<\frac{r-1}{2}$.

        The dimension of the tangent space at $[Q]$ does not change if we act on $[Q]$ by $\gl{n}$.
        Thus, after reordering the columns and rows in $Q$ so that the columns in $D_1$ are reordered from the last column to the first (see Figure~\ref{fig:D1-reordered}), and a linear change of variables, we may assume $D_1=D_s(z_0,z_1,z_2)+D'$ where $D_s(z_0,z_1,z_2)$ is as in~\eqref{eq:block-D} and $D'$ has entries in $\KK[z_3,\ldots,z_k]_1$.
        \begin{figure}
            \centering
            \begin{subfigure}[b]{.4\textwidth}
                \centering
                $
                \begin{pmatrix}
                    &     &     &        &  z_s   \\
                    & \text{\Large *} &  & \iddots &     \\
                    &        & z_2    &     &  \\
                    z_0 & z_1  &        &     &
                \end{pmatrix}_{s\times (s+1)}
                $
                \caption{The block $D_1$}\label{fig:D1}
            \end{subfigure}
            \quad\quad
            \begin{subfigure}[b]{.4\textwidth}
                \centering
                $
                \begin{pmatrix}
                    z_s &     &     &        &     \\
                    & \ddots &  & \text{\Large *} &     \\
                    &        & z_2    &     &  \\
                    &        &        & z_1    & z_0
                \end{pmatrix}_{s\times (s+1)}
                $
                \caption{$D_1$ after reordering columns}\label{fig:D1-reordered}
            \end{subfigure}
            \caption{Reordering columns in $D_1$; entries in $*$ are non-zero}
            \label{fig:reordering-columns}
        \end{figure}

        By Lemma~\ref{lem:account-of-tangent-dimension}, we have  $\dim T_{[Q]}\fano{k}{\SD[r]{n}}=a_{\det} + (\kappa(s)-k)(k+1)$,
        where $a_{\det}$ is the dimension of the $\KK$-vector space of all symmetric matrices $A$ of the form
        \[\begin{pNiceArray}{ccc}[first-row,first-col]
              && \outmat{s+1} & \outmat{n-2s-1}\\
              \outmat{s+1}&& A_1 & A_2 \\
              \outmat{n-2s-1}&& A_2^t & A_3
        \end{pNiceArray}_{(n-s)\times(n-s)},
        \]
        with entries in $\KK[z_1,
        \ldots,z_k]_1$ such that the $s$-anchored $(2s+1)\times(2s+1)$ minors of
        \begin{equation}\label{eq:mid-comp-three-A-blocks}
        \begin{pNiceArray}{ccccc}[first-row,first-col]
            &&\outmat{s}&&\outmat{s+1}&\outmat{n-2s-1}\\
            \outmat{s}&&0&&D_1&D_2 \\
            & &  &&& \\
            \outmat{s+1}&&D_1^t&&A_1&A_2\\
            & &&  && \\
            \outmat{n-2s-1}&&D_2^t&&A_2^t&A_3\\
        \end{pNiceArray}
        \end{equation}
        vanish.

        This vector space contains a subspace of dimension $s(n-s)$, namely all matrices $A=\tilde{D}+\tilde{D}^t$ where $\tilde{D}$ is a square matrix of size $(n-s)$ whose rows are in the row span of $(D_1\quad D_2)$.
        To show the scheme is smooth at $[Q]$, it is enough to prove $a_{\det}=s(n-s)$, since by Corollary~\ref{cor:dim-compk}, we will have $T_{[Q]}\fano{k}{\SD[r]{n}}=\dim\compk{k}{\frac{r-1}{2}}$.
        Equivalently, we need to show if the $s$-anchored $(2s+1)\times(2s+1)$ minors of~\eqref{eq:mid-comp-three-A-blocks} vanish, then by \rowcol{} operations using the first $s$ rows and columns we can reduce the block $A$ to zero.
        We do this in three steps.

        Step 1 (reducing $A_1$ to zero):
        Consider the submatrix in~\eqref{eq:mid-comp-three-A-blocks} given by the first $2s+1$ rows and columns.
        This submatrix is exactly the matrix in the proof of Proposition~\ref{prop:tangent-dim-smooth} (case $k\geq 2$).
        By our argument there we can use the blocks $D_1,D_1^t$ and \rowcol{} operations on~\eqref{eq:mid-comp-three-A-blocks} to reduce the block $A_1$ to zero.

        Step 2 (reducing $A_2$ to zero):
        By row-column operations and the entries on the diagonal of $D_1{^t}$, we eliminate $z_0$ from the entries in $A_2$ except possibly those in the last row of $A_2$.
        We claim this reduces $A_2$ to zero.
        For each $2s+2\leq \ell\leq n$, let $P_{\ell}$ be the submatrix of~\eqref{eq:mid-comp-three-A-blocks} given by rows $1,\ldots, 2s+1$ and columns $1,\ldots,2s,\ell$:
        \begin{equation*}
            P_{\ell}=\left(\begin{array}{ccc|c|c}
                               &&&&d_{1,\ell}\\
                               &0&&\text{the matrix }D_1&\vdots\\
                               &&&\text{with last column removed}& d_{s,\ell}\\
                               \hline
                               &&&& a_{1,\ell}\\
                               &D_1^t&&0&\vdots\\
                               &&&&a_{s+1,\ell}
            \end{array}
            \right).
        \end{equation*}
        We have $\det P_{\ell}=0$ since this is an $s$-anchored minor.
        In the expansion of $\det P_{\ell}$, no $d_{i,\ell}$ appears.
        This is because each term of $\det P_{\ell}$ picks one entry from each row.
        If a term picks one of the entries $d_{1,\ell},\ldots, d_{s,\ell}$, then because $D_1^t$ is of size $(s+1)\times s$, one entry in the zero block below $D_1$ must be chosen, which makes the term zero.

        Expanding $\det P_{\ell}$ along the last column gives
        \[
            \det P_{\ell}=\sum_{i=1}^{s+1} (-1)^{s+i}a_{i,\ell} p_i p_{s+1},
        \]
        where $p_i$ is the minor of $D_1$ after deleting column $i$.
        Since $p_{s+1}\neq 0$, we have $\sum_{i=1}^{s+1} (-1)^{s+i}a_{i,\ell}p_i=0$.

        Modulo $z_3,\ldots,z_k$, we have $D_2=D_s(z_0,z_1,z_2)$ and by Lemma~\ref{lem:lin-indep-minors}, the smallest $z_0$-degree in $p_i$ is $i-1$.
        By assumption, $a_{1,\ell},\ldots,a_{s,\ell}$ do not involve $z_0$.
        Inspecting the $z_0$-degrees in $\sum_{i=1}^{s+1} (-1)^{s+i}a_{i,\ell}p_i=0$ yields $a_{i,\ell}=0$ for all $i$, hence $a_{1,\ell},\ldots,a_{s+1,\ell}$ do not involve $z_0, z_1,z_2$.

        Pick $3\leq j\leq k$.
        Modulo $z_3,\ldots,z_{j-1},z_{j+1},\ldots,z_k$, we have $a_{i,\ell}=\alpha_{i,j,\ell} z_j$ for some $\alpha_{i,j,\ell}\in \KK$.
        Therefore
        \[
            0=\sum_{i=1}^{s+1} (-1)^{s+i}a_{i,\ell} p_i=z_j\sum_{i=1}^{s+1} (-1)^{s+i}\alpha_{i,j,\ell} p_i,
        \]
        which implies $\sum_{i=1}^{s+1}(-1)^{s+i}\alpha_{i,j,\ell} p_i=0$.
        By inspecting the degree of $z_0$, we have $\alpha_{i,j,\ell}=0$ for all $i,\ell$, i.e.\ $a_{i,\ell}$ does not involve $z_j$.
        As $i,j,\ell$ are arbitrary, we conclude the entries in $A_2$ are all zero.

        Step 3 ($A_3$ is now reduced to zero):
        To show $A_3$ is reduced to zero, pick an entry $a_{\ell,\ell'}$ at position $(\ell,\ell')$ in~\eqref{eq:mid-comp-three-A-blocks} located in the block $A_3$.
        The minor corresponding to rows $1,\ldots,2s,\ell$ and columns $1,\ldots,2s,\ell'$ is $p_{s+1}^2 a_{\ell,\ell'}$.
        Since this minor vanishes and $p_{s+1}\neq 0$, we get $a_{\ell,\ell'}=0$.
        This completes the proof of $a_{\det}=s(n-s)$.
    \end{proof}

    \begin{example}[Complete description of $\fano{k}{\SD[3]{3}}$]\label{ex:3by3matrix}
       We fully describe the schemes $\fano{k}{\SD[3]{3}}$, parametrizing linear spaces of singular plane conics.
        We have $0\leq s\leq 1$, and $1\leq k\leq 2$ (by Proposition~\ref{prop:max-dim-non-empty}).
        When $k=2$, because $\kappa(0)=\kappa(1)=2$, Corollary~\ref{cor:kappa-s-component} gives that the scheme $\fano{2}{\SD[3]{3}}$ is the disjoint union of two irreducible components $\compk{2}{0}$, $\compk{2}{1}$ each of dimension 2 and each a single $\gl{3}$ orbit (see also~\cite[§4]{nets-of-conics}).
        The component $\compk{2}{0}$ is generically non-reduced and is isomorphic in its reduced structure to $\PP^2$, and $\compk{2}{1}$ is isomorphic to $\PP^2$ (see Proposition~\ref{prop:when-gen-non-reduced} and Corollary~\ref{cor:kappa-s-component}).

        When $k=1$, the scheme $\fano{1}{\SD[3]{3}}$ has two intersecting components $\compk{1}{0}$, $\compk{1}{1}$, each of dimension 4 (see Theorem~\ref{thm:fano-lines-components}).
        The component $\compk{1}{0}$ is generically non-reduced and $\compk{1}{1}$ is generically reduced (see Proposition~\ref{prop:when-gen-non-reduced}).
        For a computation of the classes of the orbit closures in the Chow ring of $\Gr(2,\sym_3)$, see~\cite{dhruv2023chow}.
    \end{example}

    \section{Further questions}\label{sec:further-questions}
    In our investigation of the Fano schemes $\fano{k}{\SD[r]{n}}$, we proved in particular cases, listed in Table~\ref{tab:summary-symmetric-components}, that the closed subscheme $\compk{k}{s}$ consisting of nested $s$-compression spaces gives an irreducible component of $\fano{k}{\SD[n]{r}}$.
    Based on these results, we make the following conjecture.
    \begin{conjecture}\label{conj:symmetric-ck-component}
        For every $n,r,k,s$, the subscheme $\compk{k}{s}$ is an irreducible component of the Fano scheme $\fano{k}{\SD[r]{n}}$ taken with the reduced subscheme structure.
    \end{conjecture}

    Our next question involves Theorem~\ref{thm:pazzis} stating that the Fano scheme consists only of compression spaces when $k>\max\{\kappa(1),\kappa(s_{\max}-1)\}$.
    Here $s_{\max}=\lfloor\frac{r-1}{2}\rfloor$.
    It is of interest to know if this bound can be improved.
    One reason this is important is that when the closed points of $\fano{k}{\SD[r]{n}}$ are all compression spaces, we can classify the irreducible components of the scheme.
    \begin{question}
        Given $n$ and $r$, what is the smallest and largest value $k$ for which $\fano{k}{\SD[r]{n}}$ has a non-compression point?
    \end{question}

    The Kronecker canonical form classifies pencils of symmetric matrices of bounded rank up to a linear change of coordinates.
    It may be too much to ask to have a complete classification for linear spaces in higher dimensions but there may be hope for enumerating the irreducible components of $\fano{k}{\SD[r]{n}}$.

    \begin{problem}
        Enumerate the irreducible components of $\fano{k}{\SD[r]{n}}$ for $k>1$.
    \end{problem}

    \printbibliography

@book {3264,
    AUTHOR = {Eisenbud, David and Harris, Joe},
    TITLE = {3264 and all that---a second course in algebraic geometry},
    PUBLISHER = {Cambridge University Press, Cambridge},
    YEAR = {2016},
    PAGES = {xiv+616},
    ISBN = {978-1-107-60272-4; 978-1-107-01708-5},
    MRCLASS = {14-01 (14C15 14M15 14N10)},
    MRNUMBER = {3617981},
    MRREVIEWER = {Arnaud Beauville},
    DOI = {10.1017/CBO9781139062046},
}

@article {horrocks,
    AUTHOR = {Horrocks, G.},
    TITLE = {Fixed point schemes of additive group actions},
    JOURNAL = {Topology},
    FJOURNAL = {Topology. An International Journal of Mathematics},
    VOLUME = {8},
    YEAR = {1969},
    PAGES = {233--242},
    ISSN = {0040-9383},
    MRCLASS = {14.50},
    MRNUMBER = {244261},
    MRREVIEWER = {James\ Milne},
    DOI = {10.1016/0040-9383(69)90013-5},
}

@article {hartshorne-connectedness,
    AUTHOR = {Hartshorne, Robin},
    TITLE = {Connectedness of the {H}ilbert scheme},
    JOURNAL = {Inst. Hautes \'{E}tudes Sci. Publ. Math.},
    FJOURNAL = {Institut des Hautes \'{E}tudes Scientifiques. Publications
              Math\'{e}matiques},
    NUMBER = {29},
    YEAR = {1966},
    PAGES = {5--48},
    ISSN = {0073-8301},
    MRCLASS = {14.55},
    MRNUMBER = {213368},
    MRREVIEWER = {H. R\"{o}hrl},
    URL = {http://www.numdam.org/item?id=PMIHES_1966__29__5_0},
}

@article {atkinson-lloyd,
    AUTHOR = {Atkinson, M. D. and Lloyd, S.},
    TITLE = {Large spaces of matrices of bounded rank},
    JOURNAL = {Quart. J. Math. Oxford Ser. (2)},
    FJOURNAL = {The Quarterly Journal of Mathematics. Oxford. Second Series},
    VOLUME = {31},
    YEAR = {1980},
    NUMBER = {123},
    PAGES = {253--262},
    ISSN = {0033-5606},
    MRCLASS = {15A30},
    MRNUMBER = {587090},
    MRREVIEWER = {A. D. Porter},
    DOI = {10.1093/qmath/31.3.253},
}

@article {atkinson-rank3,
    AUTHOR = {Atkinson, M. D.},
    TITLE = {Primitive spaces of matrices of bounded rank. {II}},
    JOURNAL = {J. Austral. Math. Soc. Ser. A},
    FJOURNAL = {Australian Mathematical Society. Journal. Series A. Pure
              Mathematics and Statistics},
    VOLUME = {34},
    YEAR = {1983},
    NUMBER = {3},
    PAGES = {306--315},
    ISSN = {0263-6115},
    MRCLASS = {15A30},
    MRNUMBER = {695915},
    MRREVIEWER = {J.\ L.\ Brenner},
}

@misc{landsberg-rank4,
    title={On Linear spaces of of matrices bounded rank},
    author={Hang Huang and J. M. Landsberg},
    year={2023},
    eprint={2306.14428},
    archivePrefix={arXiv},
    primaryClass={math.AG}
}

@article {barth-fano,
    AUTHOR = {Barth, W. and Van de Ven, A.},
    TITLE = {Fano varieties of lines on hypersurfaces},
    JOURNAL = {Arch. Math. (Basel)},
    FJOURNAL = {Archiv der Mathematik},
    VOLUME = {31},
    YEAR = {1978/79},
    NUMBER = {1},
    PAGES = {96--104},
    ISSN = {0003-889X},
    MRCLASS = {14C05 (14M99)},
    MRNUMBER = {510081},
    MRREVIEWER = {Allen B. Altman},
    DOI = {10.1007/BF01226420},
}

@article {langer-fano,
    AUTHOR = {Langer, Adrian},
    TITLE = {Fano schemes of linear spaces on hypersurfaces},
    JOURNAL = {Manuscripta Math.},
    FJOURNAL = {Manuscripta Mathematica},
    VOLUME = {93},
    YEAR = {1997},
    NUMBER = {1},
    PAGES = {21--28},
    ISSN = {0025-2611},
    MRCLASS = {14J45 (14J70)},
    MRNUMBER = {1446187},
    MRREVIEWER = {P. Abellanas},
    DOI = {10.1007/BF02677454},
}

@article {nets-of-conics,
    AUTHOR = {Abdallah, Nancy and Emsalem, Jacques and Iarrobino, Anthony},
    TITLE = {Nets of conics and associated {A}rtinian algebras of length 7},
    JOURNAL = {Eur. J. Math.},
    FJOURNAL = {European Journal of Mathematics},
    VOLUME = {9},
    YEAR = {2023},
    NUMBER = {2},
    PAGES = {[Paper No. 22], 71},
    ISSN = {2199-675X,2199-6768},
    MRCLASS = {14C21 (01A75 13E10 14-03)},
    MRNUMBER = {4566789},
    DOI = {10.1007/s40879-023-00600-9},
    URL = {https://doi.org/10.1007/s40879-023-00600-9},
}

@article {debarre-fano,
    AUTHOR = {Debarre, Olivier and Manivel, Laurent},
    TITLE = {Sur la vari\'{e}t\'{e} des espaces lin\'{e}aires contenus dans une
              intersection compl\`ete},
    JOURNAL = {Math. Ann.},
    FJOURNAL = {Mathematische Annalen},
    VOLUME = {312},
    YEAR = {1998},
    NUMBER = {3},
    PAGES = {549--574},
    ISSN = {0025-5831},
    MRCLASS = {14M15 (14J40)},
    MRNUMBER = {1654757},
    MRREVIEWER = {Yuri G. Prokhorov},
    DOI = {10.1007/s002080050235},
}

@article {ilten-kelly,
    AUTHOR = {Ilten, Nathan and Kelly, Tyler L.},
    TITLE = {Fano schemes of complete intersections in toric varieties},
    JOURNAL = {Math. Z.},
    FJOURNAL = {Mathematische Zeitschrift},
    VOLUME = {300},
    YEAR = {2022},
    NUMBER = {2},
    PAGES = {1529--1556},
    ISSN = {0025-5874},
    MRCLASS = {14M25 (14D20 14M10)},
    MRNUMBER = {4363787},
    DOI = {10.1007/s00209-021-02809-4},
}

@article {altman-fano,
    AUTHOR = {Altman, Allen B. and Kleiman, Steven L.},
    TITLE = {Foundations of the theory of {F}ano schemes},
    JOURNAL = {Compositio Math.},
    FJOURNAL = {Compositio Mathematica},
    VOLUME = {34},
    YEAR = {1977},
    NUMBER = {1},
    PAGES = {3--47},
    ISSN = {0010-437X},
    MRCLASS = {14C05 (14C20 14M15 14N10)},
    MRNUMBER = {569043},
    MRREVIEWER = {Joel Roberts},
    URL = {http://www.numdam.org/item?id=CM_1977__34_1_3_0},
}

@article{fano1904,
    title={Sulle superficie algebriche contenute in una varieta cubica dello spazio a quattro dimensioni},
    author={Fano, Gino},
    journal={Atti Ace. Torino},
    volume={39},
    pages={597--615},
    year={1904}
}

@article {pazzi-revisited,
    AUTHOR = {de Seguins Pazzis, Cl\'{e}ment},
    TITLE = {Large spaces of bounded rank matrices revisited},
    JOURNAL = {Linear Algebra Appl.},
    FJOURNAL = {Linear Algebra and its Applications},
    VOLUME = {504},
    YEAR = {2016},
    PAGES = {124--189},
    ISSN = {0024-3795},
    MRCLASS = {15A30 (15A03)},
    MRNUMBER = {3502533},
    MRREVIEWER = {LeRoy B. Beasley},
    DOI = {10.1016/j.laa.2016.03.051},
}

@incollection {beasley2,
    AUTHOR = {Beasley, LeRoy B.},
    TITLE = {Null spaces of spaces of matrices of bounded rank},
    BOOKTITLE = {Current trends in matrix theory ({A}uburn, {A}la., 1986)},
    PAGES = {45--50},
    PUBLISHER = {North-Holland, New York},
    YEAR = {1987},
    MRCLASS = {15A21},
    MRNUMBER = {898892},
    MRREVIEWER = {M. Pearl},
}

@article {buchsbaum,
    AUTHOR = {Buchsbaum, David A. and Eisenbud, David},
    TITLE = {Algebra structures for finite free resolutions, and some
              structure theorems for ideals of codimension {$3$}},
    JOURNAL = {Amer. J. Math.},
    FJOURNAL = {American Journal of Mathematics},
    VOLUME = {99},
    YEAR = {1977},
    NUMBER = {3},
    PAGES = {447--485},
    ISSN = {0002-9327},
    MRCLASS = {13D15},
    MRNUMBER = {453723},
    MRREVIEWER = {M. Nagata},
    DOI = {10.2307/2373926},
}

@article {kleppe,
    AUTHOR = {Kleppe, H. and Laksov, D.},
    TITLE = {The algebraic structure and deformation of {P}faffian schemes},
    JOURNAL = {J. Algebra},
    FJOURNAL = {Journal of Algebra},
    VOLUME = {64},
    YEAR = {1980},
    NUMBER = {1},
    PAGES = {167--189},
    ISSN = {0021-8693},
    MRCLASS = {14M12 (14B07)},
    MRNUMBER = {575789},
    MRREVIEWER = {Luchezar L. Avramov},
    DOI = {10.1016/0021-8693(80)90140-4},
}

@article {flanders,
    AUTHOR = {Flanders, H.},
    TITLE = {On spaces of linear transformations with bounded rank},
    JOURNAL = {J. London Math. Soc.},
    FJOURNAL = {The Journal of the London Mathematical Society},
    VOLUME = {37},
    YEAR = {1962},
    PAGES = {10--16},
    ISSN = {0024-6107},
    MRCLASS = {15.40},
    MRNUMBER = {136618},
    MRREVIEWER = {Marvin Marcus},
    DOI = {10.1112/jlms/s1-37.1.10},
}

@article {parfenov-orbits,
    AUTHOR = {Parfenov, P. G.},
    TITLE = {Orbits and their closures in the spaces {${\Bbb
              C}^{k_1}\otimes\dots\otimes{\Bbb C}^{k_r}$}},
    JOURNAL = {Mat. Sb.},
    FJOURNAL = {Matematicheski\u{\i} Sbornik},
    VOLUME = {192},
    YEAR = {2001},
    NUMBER = {1},
    PAGES = {89--112},
    ISSN = {0368-8666},
    MRCLASS = {14L30},
    MRNUMBER = {1830474},
    MRREVIEWER = {Dmitri I. Panyushev},
    DOI = {10.1070/SM2001v192n01ABEH000537},
}

@misc{mezzetti,
    title={Pencils of singular quadrics of constant rank and their orbits},
    author={Ada Boralevi and Emilia Mezzetti},
    year={2022},
    eprint={2206.00503},
    archivePrefix={arXiv},
    primaryClass={math.AG}
}

@misc{dhruv2023chow,
    title={The Chow Ring Classes of $\mathrm{PGL}_3$ Orbit Closures in $\mathbb{G}(1, 5)$},
    author={Gaurav and Goel},
    year={2023},
    eprint={2310.18571},
    archivePrefix={arXiv},
    primaryClass={math.AG}
}

@article {degen-landsberg,
    AUTHOR = {Ilic, Bo and Landsberg, J. M.},
    TITLE = {On symmetric degeneracy loci, spaces of symmetric matrices of
              constant rank and dual varieties},
    JOURNAL = {Math. Ann.},
    FJOURNAL = {Mathematische Annalen},
    VOLUME = {314},
    YEAR = {1999},
    NUMBER = {1},
    PAGES = {159--174},
    ISSN = {0025-5831},
    MRCLASS = {14N05 (14J60 14M99 15A30)},
    MRNUMBER = {1689267},
    MRREVIEWER = {Fyodor L. Zak},
    DOI = {10.1007/s002080050290},
}

@book {gantmacher,
    AUTHOR = {Gantmacher, F. R.},
    TITLE = {The theory of matrices. {V}ols. 1, 2},
    NOTE = {Translated by K. A. Hirsch},
    PUBLISHER = {Chelsea Publishing Co., New York},
    YEAR = {1959},
    PAGES = {Vol. 1, x+374 pp. Vol. 2, ix+276},
    MRCLASS = {15.00 (65.00)},
    MRNUMBER = {0107649},
}

@article {sturmfels-pencil,
    AUTHOR = {Fevola, C. and Mandelshtam, Y. and Sturmfels, B.},
    TITLE = {Pencils of quadrics: old and new},
    JOURNAL = {Matematiche (Catania)},
    FJOURNAL = {Le Matematiche},
    VOLUME = {76},
    YEAR = {2021},
    NUMBER = {2},
    PAGES = {319--335},
    ISSN = {0373-3505},
    MRCLASS = {14Q15 (14M15 14N20 62F10 62R01)},
    MRNUMBER = {4334890},
    MRREVIEWER = {Guillermo Matera},
    DOI = {10.4418/2021.76.2.2},
}

@article {moduli-333,
    AUTHOR = {Ng, Kok Onn},
    TITLE = {The moduli space of {$(3,3,3)$} trilinear forms},
    JOURNAL = {Manuscripta Math.},
    FJOURNAL = {Manuscripta Mathematica},
    VOLUME = {88},
    YEAR = {1995},
    NUMBER = {1},
    PAGES = {87--107},
    ISSN = {0025-2611},
    MRCLASS = {14D25},
    MRNUMBER = {1348793},
    MRREVIEWER = {Yi Hu},
    DOI = {10.1007/BF02567808},
}

@article {dieudonne,
    AUTHOR = {Dieudonn\'{e}, Jean},
    TITLE = {Sur une g\'{e}n\'{e}ralisation du groupe orthogonal \`a quatre
              variables},
    JOURNAL = {Arch. Math.},
    FJOURNAL = {Archiv der Mathematik},
    VOLUME = {1},
    YEAR = {1949},
    PAGES = {282--287},
    ISSN = {0003-889X},
    MRCLASS = {09.0X},
    MRNUMBER = {29360},
    MRREVIEWER = {L. K. Hua},
    DOI = {10.1007/BF02038756},
    URL = {https://doi.org/10.1007/BF02038756},
}

@article {goto,
    AUTHOR = {Goto, Shiro},
    TITLE = {The divisor class group of a certain {K}rull domain},
    JOURNAL = {J. Math. Kyoto Univ.},
    FJOURNAL = {Journal of Mathematics of Kyoto University},
    VOLUME = {17},
    YEAR = {1977},
    NUMBER = {1},
    PAGES = {47--50},
    ISSN = {0023-608X},
    MRCLASS = {13G05},
    MRNUMBER = {435063},
    MRREVIEWER = {David\ W.\ Sharpe},
    DOI = {10.1215/kjm/1250522811},
    URL = {https://doi.org/10.1215/kjm/1250522811},
}

@book {vakil-rising-sea,
    AUTHOR = {Vakil, Ravi},
    TITLE = {The rising sea---foundations of algebraic geometry},
    PUBLISHER = {Princeton University Press, Princeton, NJ},
    YEAR = {[2025] \copyright 2025},
    PAGES = {xxiii+662},
    ISBN = {978-0-691-26866-8; 978-0-691-26867-5; 978-0-691-26868-2},
    MRCLASS = {14-01},
    MRNUMBER = {4942570},
}

@book {borel,
    AUTHOR = {Borel, Armand},
    TITLE = {Linear algebraic groups},
    SERIES = {Graduate Texts in Mathematics},
    VOLUME = {126},
    EDITION = {Second},
    PUBLISHER = {Springer-Verlag, New York},
    YEAR = {1991},
    PAGES = {xii+288},
    ISBN = {0-387-97370-2},
    MRCLASS = {20-01 (20Gxx)},
    MRNUMBER = {1102012},
    MRREVIEWER = {F. D. Veldkamp},
    DOI = {10.1007/978-1-4612-0941-6},
}

@book {hartshorne-deformation,
    AUTHOR = {Hartshorne, Robin},
    TITLE = {Deformation theory},
    SERIES = {Graduate Texts in Mathematics},
    VOLUME = {257},
    PUBLISHER = {Springer, New York},
    YEAR = {2010},
    PAGES = {viii+234},
    ISBN = {978-1-4419-1595-5},
    MRCLASS = {14D15 (13D10 14B07 14B12)},
    MRNUMBER = {2583634},
    MRREVIEWER = {Arvid\ Siqveland},
    DOI = {10.1007/978-1-4419-1596-2},
    URL = {https://doi.org/10.1007/978-1-4419-1596-2},
}

@article {kutz74,
    AUTHOR = {Kutz, Ronald E.},
    TITLE = {Cohen-{M}acaulay rings and ideal theory in rings of invariants
              of algebraic groups},
    JOURNAL = {Trans. Amer. Math. Soc.},
    FJOURNAL = {Transactions of the American Mathematical Society},
    VOLUME = {194},
    YEAR = {1974},
    PAGES = {115--129},
    ISSN = {0002-9947},
    MRCLASS = {13C15 (13H10 14M15 15A72 20G15)},
    MRNUMBER = {352082},
    MRREVIEWER = {Melvin Hochster},
    DOI = {10.2307/1996796},
}

@article {fano-nathan-chan,
    AUTHOR = {Chan, Melody and Ilten, Nathan},
    TITLE = {Fano schemes of determinants and permanents},
    JOURNAL = {Algebra Number Theory},
    FJOURNAL = {Algebra \& Number Theory},
    VOLUME = {9},
    YEAR = {2015},
    NUMBER = {3},
    PAGES = {629--679},
    ISSN = {1937-0652},
    MRCLASS = {14M12 (14J45 14M15)},
    MRNUMBER = {3340547},
    MRREVIEWER = {Ulrich Derenthal},
    DOI = {10.2140/ant.2015.9.629},
}

@article {meshulam-max-dim,
    AUTHOR = {Meshulam, Roy},
    TITLE = {On two extremal matrix problems},
    JOURNAL = {Linear Algebra Appl.},
    FJOURNAL = {Linear Algebra and its Applications},
    VOLUME = {114/115},
    YEAR = {1989},
    PAGES = {261--271},
    ISSN = {0024-3795},
    MRCLASS = {15A03 (15A15)},
    MRNUMBER = {986879},
    MRREVIEWER = {Lajos L\'{a}szl\'{o}},
    DOI = {10.1016/0024-3795(89)90465-5},
}

@article {jensen,
    AUTHOR = {Jensen, J. L. W. V.},
    TITLE = {Sur une identit\'{e} d'{A}bel et sur d'autres formules analogues},
    JOURNAL = {Acta Math.},
    FJOURNAL = {Acta Mathematica},
    VOLUME = {26},
    YEAR = {1902},
    NUMBER = {1},
    PAGES = {307--318},
    ISSN = {0001-5962},
    MRCLASS = {DML},
    MRNUMBER = {1554966},
    DOI = {10.1007/BF02415499},
}

@incollection {loewy-radwan,
    AUTHOR = {Loewy, Raphael and Radwan, Nizar},
    TITLE = {Spaces of symmetric matrices of bounded rank},
    NOTE = {Second Conference of the International Linear Algebra Society
              (ILAS) (Lisbon, 1992)},
    JOURNAL = {Linear Algebra Appl.},
    FJOURNAL = {Linear Algebra and its Applications},
    VOLUME = {197/198},
    YEAR = {1994},
    PAGES = {189--215},
    ISSN = {0024-3795},
    MRCLASS = {15A03},
    MRNUMBER = {1275615},
    MRREVIEWER = {Donald W. Robinson},
    DOI = {10.1016/0024-3795(94)90488-X},
}

@article {pazzis-classification,
    AUTHOR = {de Seguins Pazzis, Cl\'{e}ment},
    TITLE = {Large spaces of symmetric or alternating matrices with bounded
              rank},
    JOURNAL = {Linear Algebra Appl.},
    FJOURNAL = {Linear Algebra and its Applications},
    VOLUME = {508},
    YEAR = {2016},
    PAGES = {146--189},
    ISSN = {0024-3795},
    MRCLASS = {15A30 (15A03)},
    MRNUMBER = {3542986},
    MRREVIEWER = {Tatjana Petek},
    DOI = {10.1016/j.laa.2016.07.005},
}

@article {dopico,
    AUTHOR = {De Ter\'{a}n, Fernando and Dmytryshyn, Andrii and Dopico,
              Froil\'{a}n M.},
    TITLE = {Generic symmetric matrix pencils with bounded rank},
    JOURNAL = {J. Spectr. Theory},
    FJOURNAL = {Journal of Spectral Theory},
    VOLUME = {10},
    YEAR = {2020},
    NUMBER = {3},
    PAGES = {905--926},
    ISSN = {1664-039X,1664-0403},
    MRCLASS = {15A22 (15A18 15A21)},
    MRNUMBER = {4164010},
    MRREVIEWER = {Sanne\ ter Horst},
    DOI = {10.4171/jst/316},
    URL = {https://doi.org/10.4171/jst/316},
}

@article {waterhouse,
    AUTHOR = {Waterhouse, William C.},
    TITLE = {The codimension of singular matrix pairs},
    JOURNAL = {Linear Algebra Appl.},
    FJOURNAL = {Linear Algebra and its Applications},
    VOLUME = {57},
    YEAR = {1984},
    PAGES = {227--245},
    ISSN = {0024-3795,1873-1856},
    MRCLASS = {15A21 (15A18)},
    MRNUMBER = {729275},
    MRREVIEWER = {George\ P.\ Barker},
    DOI = {10.1016/0024-3795(84)90190-3},
    URL = {https://doi.org/10.1016/0024-3795(84)90190-3},
}
\end{document}